\documentclass[13pt, reqno]{amsart}

\usepackage{graphicx}
\usepackage[normalem]{ulem}
\usepackage{mathrsfs}
\usepackage{amssymb}
\usepackage{amsfonts, mathtools}
\usepackage{amssymb,amsmath,amsthm,tikz-cd, tikz}
\usepackage{stix}

\tikzcdset{arrow style=tikz, diagrams={>=stealth'}}

\def\ds{{\delta, \sigma}}

\usepackage{hyperref}
\usepackage{palatino}
\usepackage{wrapfig}
\usepackage{subcaption}

\usetikzlibrary{patterns, matrix,shapes,arrows,calc,topaths,intersections,positioning,decorations.pathreplacing,decorations.pathmorphing,fit, 3d}

\newif\ifpaper
\papertrue

\usepackage{color}

\title {Variation of GIT and Variation of Lagrangian Skeletons II: Quasi-Symmetric Case}
\author{Jesse Huang and Peng Zhou}
\date{\today}

\usepackage{color}

\renewcommand{\ss}{\subsection}
\newcommand{\sss}{\subsubsection}

\DeclareMathOperator{\Aff}{Aff}

\DeclareMathOperator{\cone}{cone}

\DeclareMathOperator{\conv}{conv}

\DeclareMathOperator{\Coh}{Coh}

\DeclareMathOperator{\Ext}{Ext}

\DeclareMathOperator{\Fun}{Fun}

\DeclareMathOperator{\Hom}{Hom}

\newcommand{\uhom}{{\hcal  om}}
\DeclareMathOperator{\holim}{holim}

\DeclareMathOperator{\Int}{Int}

\renewcommand{\Re}{\text{Re}}

\renewcommand{\span}{\z{span}}

\DeclareMathOperator{\nbhd}{nbhd}

\DeclareMathOperator{\supp}{Supp}
\DeclareMathOperator{\sign}{sign}

\DeclareMathOperator{\Vect}{Vect}

\newcommand{\NN}{\mathsf{N}}
\newcommand{\MM}{\mathsf{M}}

\newcommand\ceil[1]{\lceil #1 \rceil}
\newcommand\floor[1]{\lfloor #1 \rfloor}

\newcommand{\dm}{\diamondsuit}

\newcommand{\z}{\text}

\newcommand{\signset}{\{+, -, 0\}}

\newcommand{\pa}{\partial}

\newcommand{\la}{\langle}
\newcommand{\ra}{\rangle}
\newcommand{\ot}{\otimes}

\newcommand{\bu}{\bullet}

\newcommand{\half}{\frac{1}{2}}

\newcommand{\wt}{\widetilde}
\newcommand{\wh}{\widehat}

\newcommand{\wb}{\overline}

\newcommand{\acal}{\mathcal{A}}
\newcommand{\bcal}{\mathcal{B}}
\newcommand{\ccal}{\mathcal{C}}
\newcommand{\dcal}{\mathcal{D}}
\newcommand{\ecal}{\mathcal{E}}
\newcommand{\fcal}{\mathcal{F}}

\newcommand{\hcal}{\mathcal{H}}
\newcommand{\ical}{\mathfrak{I}} 

\newcommand{\lcal}{\mathcal{L}}
\newcommand{\mcal}{\mathcal{M}}

\newcommand{\pcal}{\mathcal{P}}

\newcommand{\rcal}{\mathcal{R}}
\newcommand{\scal}{\mathcal{S}}

\newcommand{\wcal}{\mathcal{W}}
\newcommand{\xcal}{\mathcal{X}}

\newcommand{\gfrak}{\mathfrak{g}}

\newcommand{\B}{\mathbb{B}}

\newcommand{\C}{\mathbb{C}}

\newcommand{\R}{\mathbb{R}}
\newcommand{\bS}{\mathbb{S}}

\newcommand{\Z}{\mathbb{Z}}

\newcommand{\T}{\mathbb{T}}
\renewcommand{\P}{\mathbb{P}}

\newcommand{\xto}{\xrightarrow}

\newcommand{\LRA}{\Leftrightarrow}
\newcommand{\LA}{\Leftarrow}
\newcommand{\RA}{\Rightarrow}
\newcommand{\RM}{\backslash}

\newcommand{\into}{\hookrightarrow}

\newcommand{\bea}{\begin{eqnarray*} }
\newcommand{\eea}{\end{eqnarray*} }
\newcommand{\be}{\begin{equation} }
\newcommand{\ee}{\end{equation} }
\newcommand{\bp}{\begin{proposition}}
\newcommand{\ep}{\end{proposition}}
\newcommand{\bt}{\begin{theo}}
\newcommand{\et}{\end{theo}}

\newcommand{\btu}{\begin{theou}}
\newcommand{\etu}{\end{theou}}
\newcommand{\bpf}{\begin{proof}}
\newcommand{\epf}{\end{proof}}
\newcommand{\bl}{\begin{lemma}}
\newcommand{\el}{\end{lemma}}
\newcommand{\bc}{\begin{corollary}}
\newcommand{\ec}{\end{corollary}}
\newcommand{\bd}{\begin{definition}}
\newcommand{\ed}{\end{definition}}
\newcommand{\bex}{\begin{example}}
\newcommand{\eex}{\end{example}}
\newcommand{\bA}{\left(\begin{array}}
\newcommand{\eA}{\end{array}\right)}
\newcommand{\bma}{\begin{bmatrix}}
\newcommand{\ema}{\end{bmatrix}}
\newcommand{\bcd}{\begin{tikzcd}}
\newcommand{\ecd}{\end{tikzcd}}
\newcommand{\bcs}{\begin{cases}}
\newcommand{\ecs}{\end{cases}}
\newcommand{\bee}{\begin{eqnarray} }
\newcommand{\eee}{\end{eqnarray} }
\newcommand{\brem}{\begin{remark}}
\newcommand{\erem}{\end{remark}}
\newcommand{\bnum}{\begin{enumerate}}
\newcommand{\enum}{\end{enumerate}}

\newtheorem*{theou}{Main Theorem}
\newtheorem{theo}{Theorem}[section]
\newtheorem{lemma}[theo]{Lemma}
\newtheorem{corollary}[theo]{Corollary}
\newtheorem{proposition}[theo]{Proposition}

\newtheorem{definition}[theo]{Definition}
\newtheorem{remark}[theo]{Remark}

\newtheorem{examplex}[theo]{Example}
\newenvironment{example}
  {\pushQED{\qed}\examplex}
  {\popQED\endexamplex}

\theoremstyle{plain}
\numberwithin{equation}{section}

\newcommand{\CS}{{\C^*}}

\newcommand{\La}{\Lambda}

\newcommand{\In}{\subset}

\newcommand{\congto}{\xto{\sim}}
\renewcommand{\cong}{\simeq}

\newcommand*{\hrlen}{5}
\newcommand*{\hramp}{3}

\tikzset{
asdstyle/.style={blue,thick},
righthairs/.style={postaction={decorate,draw,decoration={border,amplitude=\hramp,segment length=\hrlen,angle=-90,pre=moveto,pre length=\hrlen/2}}},
lefthairs/.style={postaction={decorate,draw,decoration={border,amplitude=\hramp,segment length=\hrlen,angle=90,pre=moveto,pre length=\hrlen/2}}},
righthairsnogap/.style={postaction={decorate,draw,decoration={border,amplitude=\hramp,segment length=\hrlen,angle=-90}}},
lefthairsnogap/.style={postaction={decorate,draw,decoration={border,amplitude=\hramp,segment length=\hrlen,angle=90}}},
graphstyle/.style={thick},
arrowstyle/.style={thick,decorate,decoration={snake,amplitude=1.7,segment length=10pt,post length=.5mm,pre length=0}},
genmapstyle/.style={thick,-stealth'},
arrhdstyle/.style={thick},
exceptarcstyle/.style={red, ultra thick},
dualquiverstyle/.style={thick,->},
patstyle/.style={pattern color = gray, pattern = north east lines, opacity=0.3}
}

\begin{document}

\begin{abstract}
Consider $(\C^*)^k$ acting on $\C^N$ satisfying certain 'quasi-symmetric' condition which produces a class of toric Calabi-Yau GIT quotient stacks. Using subcategories of $\Coh([\C^N / (\C^*)^k])$ generated by line bundles whose weights are inside certain zonotope called the 'magic window', Halpern-Leistner and Sam give a combinatorial construction of equivalences between derived categories of coherent sheaves for various GIT quotients. We apply the coherent-constructible correspondence for toric varieties to the magic windows and obtain a non-characteristic deformation of Lagrangian skeletons in $\R^{N-k}$ parameterized by $\R^k$, exhibiting derived equivalences between A-models of the various phases. Moreover, by translating the magic window zonotope in $\R^k$, we obtain a universal skeleton over $\R^k \times \R^k \RM \dcal$ for some fattening of hyperplane arrangements $\dcal$, and we show that the the universal skeleton induces a local system of categories over $\R^k \times \R^k \RM \dcal$. We also connect our results to the perverse schober structure identified by {\v{S}}penko and Van den Bergh.

\end{abstract}

\maketitle

\section{Introduction}


This is the second paper in a series \cite{zhou2020} where we construct deformation family of Lagrangian skeletons in the toric variation of GIT setting. Consider a complex torus $\T = (\C^*)^k$ acting on a vector space $V = \C^N$. There are possibly many GIT quotient stacks $[V /_\theta \T]$ depending on a character $\theta$ of $\T$. If the $\T$ action preserves the standard volume form on $\C^N$, then the resulting quotients are (non-compact) toric Calabi-Yau stacks whose derived categories of coherent sheaves are equivalent, although the equivalences are not canonical. The lack of uniqueness of GIT quotients and their equivalences is a feature instead of a bug, and one expects a local system of categories over the complexified K\"{a}hler moduli space minus some discriminant, such that away from the discriminant, various asymptotic regions in the parameter space correspond to various GIT quotients; parallel transports along paths between the regions correspond to derived equivalences. 

In general, it is difficult to construct such a local system of categories explicitly, since the discriminant can be complicated. However, under the 'quasi-symmetric' condition (see Definition \ref{d:qsym} below), the discriminant is an affine hyperplane arrangement, and we are able to construct the local system of categories using microlocal sheaf theory and Lagrangian skeletons. 

The Coherent-Constructible-Correspondence (CCC) for toric varieties was initiated by Bondal \cite{bondal2006derived} and Fang-Liu-Treumann-Zaslow \cite{FLTZ-morelli, FLTZ-hms}, and finalized in the paper by Kuwagaki \cite{Ku16} using wrapped constructible sheaves developed in \cite{nadler2016wrapped}. These results relate A-model to B-model within each GIT chamber, and our work here interpolates the various A-models across different chambers. 





\subsection{Toric GIT}
Let $n = N -k$. Let $\Z^N = \Hom( (\C^*)^N, \C^*)$ and $\Z^k = \Hom((\C^*)^k, \C^*)$ be the character lattices, and let $(\Z^N)^\vee$ and $(\Z^k)^\vee$ be their duals. Let $\{e_i\}_{i=1}^N$ be the standard basis of $\Z^N$ and $\{e_i^\vee\}_{i=1}^N$ be the dual basis of $(\Z^N)^\vee$. We assume that the torus action $(\C^*)^k$ on $\C^N$ induces short exact sequences of lattices
\be 0 \to \MM \to \Z^N \xto{\mu_\Z} \Z^k \to 0, \label{eq:Mseq} \ee
\be 0 \to (\Z^k)^\vee \to (\Z^N)^\vee \xto{\nu_\Z} \NN  \to 0,  \label{eq:Nseq}  \ee
where $\MM$ and $\NN$ are dual lattices of rank $n$. 
Let $\mu_\R: \R^N \to \R^k$ be induced from $\mu_\Z$ by $-\ot \R$. We let $T = \R / \Z \cong S^1$, and let $T^n$ be the real n-dimensional torus. We use $\MM_\R, \MM_T, \NN_\R, \NN_T$ to denote the result of tensoring the lattice with $T$ or $\R$. For simplicity, we will choose a basis of $\MM$ and identify $\MM_\R \cong \R^n$ and $\MM_T \cong T^n$ by an abuse of notation. 

Let $\beta_i \in (\Z^k)^\vee$ be the image of $e_i$, which are also the weights of the action $\T$ on $V$. We assume that $\beta_i$ are all nonzero and span $\R^k$.

\bd[Quasi-Symmetric Condition]\label{d:qsym}
Let $\lcal = \{L_1, \cdots, L_m\}$ be the collection of lines in $\R^k$ where each line contains at least one $\beta_i$.  If $$ \sum_{\beta_i \in L} \beta_i = 0 $$
for any line $L \in \lcal$, then we say the action of $\T$ on $V$ satisfies the quasi-symmetric condition. 
\ed

\brem 
This condition may seem unnatural at first sight, since there is no natural meaning of the lines $\lcal$. It may seem less imposing if one considers the following factorization of $\mu_\Z$
$$ \prod_{L \in \lcal} \Z^{[N]_L} \xto{\pi_\Z} \prod_{L \in \lcal} \Z \xto{q} \Z^k, $$
where we denote
$$ [N] = \{1,\cdots, N\}, \quad [N]_L = \{i \in [N]: \beta_i \in L\}. $$
\erem 

\subsection{Zonotopes and Windows}
We take the zonotope
$$ \nabla  := \half \sum_{i=1}^N [0, \beta_i], $$
where $[0, v]$ is the line segment connecting $0$ and $v$, and the sum is the Minkowski sum. The quasi-symmetric condition implies $\nabla = -\nabla$ and that $\nabla$ can be translated to a lattice zonotope. 


For any $\delta \in \R^k$, we define the shifted zonotope and window
$$ \nabla_\delta = \delta + \nabla, \quad W_\delta = \nabla_\delta \cap \Z^k. $$
We call $\delta$ a shift parameter, and we say $\delta$ is generic if $\pa \nabla_\delta$ contains no lattice points. In general, there is a stratification of $\R^k$ induced by the function $\delta \mapsto W_\delta$. 

For any $w \in \Z^N$, we consider a translated open positive quadrant $w + \R_{>0}^N$, and  the constructible sheaf 
\be \label{e:Lw}  L_w := \C_{w + \R_{>0}^N} \ee
locally constant with stalk $\C$ on $w + \R_{>0}^N$ and $0$ elsewhere.
We define a conical Lagrangian in $T^*\R^N$, 
\be \label{e:Law} \La_w = SS(L_w) \quad = \begin{tikzpicture}[anchor=base, baseline, scale=0.5,  shift={(0,-1)}]
\draw [righthairs] (0,2) -- (0,0) -- (2,0);
\draw [fill=gray, opacity=0.3] (0,0) rectangle (2,2);
\node at (-0.5,-0.5) {$w$};
\end{tikzpicture}. \ee
Concretely, for $w = 0$, we have
$$ \La_{ 0} = SS(\C_{\R_{>0}^N}) = \{ (x_1, \cdots, x_N, \xi_1, \cdots, \xi_N) \in T^* \R^N \mid \z{ for each $i$, either $(x_i > 0, \xi_i =0)$, or $(x_i=0, \xi_i \leq 0)$ } \}.$$

\bd [Window skeletons]
Let $\delta \in \R^k$ be a shift parameter and $W_\delta$ be a window. Then the window skeleton of $W_\delta$ is defined as the union
$$ \La_{W_\delta} = \bigcup_{w \in \wt W_\delta} \La_w  \In T^* \R^N, $$
where $ \wt W_\delta =  \mu_\Z^{-1}(W_\delta)$. We also define the non-equivariant window skeleton 
$$ \wb \La_{W_\delta} = \La_{W_\delta} / \MM \In T^*(\R^N / \MM) \cong T^*(T^n \times \R^k), $$ where $\MM$ acts on $\R^N$ and $T^* \R^N$ by translation. \ed

Recall $\mu_\R: \R^N \to \R^k$ is the induced character map from the torus action $(\C^*)^k$ on $\C^N$. For any $l \in \R^k$, we define the restrction of $\La_{W_\delta}$ to the fiber $X_l =  \mu_\R^{-1}(l)$ to be the symplectic reduction
$$ \La_{\delta, l} := \La_{W_\delta}|_{X_l} \In T^* X_l \cong T^* \R^n, $$
where we use the notion of restricting a conical Lagrangian $\La \In T^*M$ to the cotangent bundle $T^*S$ of a smooth submanifold $S \In M$, 
$$ \La|_S := (\La  \cap T^* M|_{S}) / T^*_{S} M \In  (T^* M|_S)  / T_S^* M \cong T^* S. $$
Similarly, we define the non-equivariant version of fiber skeleton
$$ \wb \La_{\delta, l} = \La_{\delta, l} / \MM \In T^* (X_l / \MM) \cong T^* T^n. $$


\subsection{Constructible sheaf of categories}
We first present our result about $\La_{W_\delta}$ for a fixed $\delta$. 

We first state a result saying $\La_{W_\delta}$ over various asymptotic region are the mirror skeletons to the various GIT quotients. Recall that $\R^k$ admits stratification according to how the fiber of $\mu_\R$ intersects with the positive quadrant $\R_{\geq 0}^N$. This stratification fits into a GKZ fan (secondary fan), where each maximal cone's interior $C$ is called a GKZ chamber, and labels the GIT quotient $X_C$. In Section \ref{s:geom-window},  we prove that for any $\delta \in \R^k$ and chamber $C$, there is an stable region $C_{\nabla_\delta}$ in the direction $C$, 
such that $\La_{W_\delta}$ restricted over $C_{\nabla_\delta}$ is the mirror to the GIT quotient $X_C$.

\bt [Theorem \ref{t:LaWLaC}]
For any $\delta \in \R^k$ and any GKZ chamber $C$, if we define 
$$ C_{\nabla_\delta} = \bigcap_{p \in \nabla_{\delta}} p + C $$
and the skeleton for GIT quotient in chamber $C$, 
$$ \La_C = \bigcup_{I \in \ical_C} (\Z^N + \R^I) \times (- (\R^\vee_{\geq 0})^{I^c}) \In T^* \R^N, \quad \ical_C = \{ I \In [N] \mid \mu_\R^{-1}(C) \cap \R_{>0}^I \neq \emptyset \}, $$
then we have
$$ \La_{W_\delta} = \La_C ,\quad \z{ inside } \mu_\R^{-1} (C_{\nabla_\delta}). $$
In particular, if $l \in C_{\nabla_\delta} \cap \Z^k$, then $\wb \La_{\delta, l} \In T^* T^n$ is the (non-equivariant) FLTZ skeleton mirror to the GIT quotient $\C^N //_l (\C^*)^k$. 
\et 
In fact, Theorem \ref{t:LaWLaC} holds for general toric GIT, without toric CY or quasi-symmetric condition.

Next, we prove a result showing that, for $\delta$ generic, $\La_{W_\delta}$ induces an equivalence among all skeletons $\La_{\delta, l}$ for $l \in \R^k$. We recall the notion of non-characteristic deformation of Lagrangian. For any smooth manifold $M$ and conical Lagrangian $\La \In T^*M$, we follow the notation of Nadler \cite{nadler2016wrapped} and consider the cocomplete dg derived category of sheaves on $M$ that are cohomologically constructible and have singular support contained in $\La$, and denote it by $Sh^\dm (M, \La)$. The compact objects in $Sh^\dm (M, \La)$ are called wrapped constructible sheaves, and constitute the full subcategory $Sh^w(M, \La)$. It is proven recently \cite{GPS3, nadler2020sheaf} that the wrapped constructible sheaves (and more generally wrapped microlocal sheaves) and the corresponding wrapped Fukaya categories $\wcal(T^*M, \La)$ are equivalent. We follow \cite{nadler2015non}, and say a family of Lagrangian skeletons $\{\La_t \In T^*M \}$ parametrized by $t$ is a non-characteristic deformation family, if the corresponding categories $Sh^\dm(M, \La_t)$ are independent of $t$. 

\bt [Theorem \ref{t:l-indep}]
If $\delta \in \R^k$ is generic, then $\{\La_{\delta, l}\}_{l \in \R^k}$ and $\{\wb \La_{\delta, l}\}_{l \in \R^k}$ are non-characteristic deformations of Lagrangian skeletons parametrized by $l \in \R^k$. 
More precisely, for any $l \in \R^k$ the restriction functors are equivalence of categories 
$$ \rho_{\delta, l}: Sh^\dm(\R^N, \La_{W_\delta}) \congto Sh^\dm(\R^n, \La_{\delta, l}), \quad \wb \rho_{\delta, l}: Sh^\dm(T^n \times \R^k, \wb \La_{W_\delta}) \congto Sh^\dm(T^n, \wb \La_{\delta, l}). $$
\et

Hence given a choice of the window $W_\delta$ for $\delta$ generic, there is a simultaneous non-characteristic interpolation among the universal FLTZ skeletons $\La_C$ in all chambers. Consequently, the A-model categories are all equivalent through parallel transport.

Our next result is about what happens when $\delta$ is non-generic.  On the B-side, the story is explained by \cite{halpern2020combinatorial} and \cite{vspenko2019class}, and we have the 'window inclusion functors' and their left and right adjoints. For each $\delta \in \R^k$, we have $\nabla_\delta, W_\delta$ as before, and we have the window category $\bcal_\delta$ (denoted as $\mcal(\nabla_\delta)$ in \cite{halpern2020combinatorial} and $\ecal_C$ in \cite{vspenko2019class} for the strata in stratification containing $\delta$) as full subcategory of $Coh([\C^N / (\C^*)^k])$ generated by line bundles $\lcal_w$ with weights $w \in W_\delta$, 
\be \label{e:bcal} \bcal_\delta := \la \{\lcal_w \mid   w \in W_\delta \} \ra. \ee
If $\delta, \delta' \in \R^k$,  and $\delta$ is a specialization of $\delta'$, i.e $W_{\delta} \supset W_{\delta'}$, then we have the window inclusion functor 
$$ \iota_{\delta', \delta}: \bcal_{\delta'} \to \bcal_{\delta}, \quad \z{ if } W_{\delta'} \In W_{\delta}. $$
The adjoint functors are studied in \cite[Lemma 6.7]{halpern2020combinatorial}, and the right adjoints are used to establish a perverse schober structure on certain periodic affine hyperplane arrangements \cite{vspenko2019class}. 

On the A-side, we give a complete Lagrangian skeletal translation of the above story. We have the same construction of the window skeleton $\La_{W_\delta}$ in $\R^N$, however the variation of skeleton $\{\La_{\delta, l}\}_{l \in \R^k}$ is no longer non-characteristic, there will be jumps in the constructible sheaf category  $Sh^\dm(\R^n, \La_{\delta, l})$ on the fiber over $l$. In Figure \ref{fig:N6k2}, the jumps are shown with blue hairy lines, and dashed lines indicate there are no jumps. Hence it is useful to adopt the language of a constructible sheaf of categories and use its singular support and microlocal stalks to indicate the 'location and amount' of the jumps\footnote{Admittedly, the use of 'constructible sheaf of categories and singular supports' is rather informal and naive in this paper. However, since the situation is combinatorial and explicit, we hope the language can be formalized without affecting the final result. }. More precisely, we consider a sheaf of categories $Sh^\dm_{\La_{W_\delta}}$ on $\R^N$ defined in \cite[Section 3.6]{nadler2016wrapped} and its push-forward along $\mu$ to $\R^k$. We denote these two sheaves of categories as
\be \label{e:def-C-delta}
\wt \ccal_\delta := Sh^\dm_{\La_{W_\delta}}, \quad \ccal_\delta := \mu_* Sh^\dm_{\La_{W_\delta}}. 
\ee 
Hence for any convex open sets $\wt U \In \R^N$ and $U \In \R^k$, we have
\be \label{e:C-delta-U} \wt \ccal_\delta(\wt U) := Sh^\dm_{\La_{W_\delta}} (\wt U),  \quad \z{ and } \quad  \ccal_\delta(U) := Sh^\dm_{\La_{W_\delta}} (\mu^{-1}(U)), \ee
where we require convexity of the open sets so that value of the sheaf and pre-sheaf agree. The same can be defined when $\La_{W_\delta}$ is replaced by the non-equivariant skeleton $\wb \La_{W_\delta}$, and/or the large constructible sheaf is replaced by the wrapped constructible sheaves
\be \ccal_\delta^w:= \mu_* Sh^w_{\La_{W_\delta}}, \quad \wb \ccal_\delta := \mu_* Sh^\dm_{\wb \La_{W_\delta}}, \quad \wb \ccal^w_\delta := \mu_* Sh^w_{\wb \La_{W_\delta}}.  \ee
Note however, that for the wrapped constructible sheaf, one in general only get a cosheaf of categories, where as for $Sh^\dm$ one get both a sheaf and co-sheaf. Hence, we will prove results for $Sh^\dm$ first, and pass to $Sh^w$ only when needed. 

The definition of singular support of $\wt \ccal_\delta$ and $\ccal_\delta$ is given in Section \ref{ss:SS-cat}. It is a straightforward generalization of the one for constructible sheaves \cite[Chapter 5]{KS}, and measures the failure of existence of unique extensions of objects and their hom spaces from a smaller open set to a bigger open set. For any given non-zero vector $(x,\xi) \in T^* \R^k$, we have the {\it microlocal restriction functor} 
\be \label{e:rho-x-xi} \rho_{x,\xi}: \ccal_\delta(B_x) \to \ccal_\delta(B_{x, \xi, -}) \ee
from a small ball $B_x$ centered at $x$ to a half ball $B_{x, \xi, -}$ obtained by 'retreating' in the $\xi$ direction. Then $\rho_{x, \xi}$ is an equivalence if and only if $\rho_{x,\xi}$ and its left-adjoint, the {\it microlocal co-restriction functor} 
\be \label{e:rho-L-x-xi} \rho_{x,\xi}^L, \ccal_\delta(B_{x, \xi, -}) \to \ccal_\delta(B_x) \ee are both fully-faithful. Hence, we define the usual singular support $SS(\ccal_\delta)$ measuring the failure of $\rho_{x,\xi}$ being an equivalence, and two smaller versions, $SS_{Hom}(\ccal_\delta)$  and $SS_{Hom}^L(\ccal_\delta)$, measuring when $\rho_{x,\xi}$ and $\rho_{x,\xi}^L$ fails to be fully-faithful, i.e. quasi-isomorphism on the hom space. By construction, the various versions of singular supports satisfy
$$ SS(\ccal_\delta) = SS_{Hom}(\ccal_\delta) \cup SS_{Hom}^L(\ccal_\delta). $$
If one thinks in terms of extending a constructible sheaf in $\ccal_\delta(U) = Sh^\dm(\pi^{-1}(U), \La_{W_\delta})$ to one in $\ccal_\delta(V)$ over a slightly larger convex open set $V \supset U$, then $SS_{Hom}^L(\ccal_\delta)$ is the obstruction for the existence of an extension, and $SS_{Hom}(\ccal_\delta)$ is the obstruction for the uniqueness of the extension.


Our next theorem describes the location of the jumping loci for $\ccal_\delta$ on $\R^k$. Let $\Sigma_\nabla$ be the exterior conormal fan to the zonotope $\nabla$, then faces of $\nabla$ are in one-to-one correspondence to the cones in $\Sigma_\nabla$. 

\bt \label{intro-t: jump-loci} [Theorem \ref{t:SS-Hom} and \ref{t:SSL-Hom-window} and ]
For any $\delta \in \R^k$, let $\ccal_\delta$ denote the constructible sheaf of categories on $\R^k$ defined in Eq. \eqref{e:def-C-delta}, then $SS_{Hom}^L(\ccal_\delta)$ is the zero section of $T^* \R^k$, and 
$$ SS(\ccal_\delta) = SS_{Hom}(\ccal_\delta) = \bigcup_{F_{\delta, \sigma}: F_{\delta, \sigma} \cap \Z^k \neq \emptyset} \Aff(F_{\delta, \sigma}) \times (-\sigma) \In \R^k \times (\R^k)^\vee \cong T^* \R^k  $$
where $F_{\delta, \sigma}$ is the face of the shifted zonotope $ \nabla_\delta = \delta + \nabla$ labelled by a cone $\sigma \In \Sigma_\nabla$. 
\et

\begin{remark}
We do not have a good way to prove that $SS_{Hom}^L(\ccal_\delta)$ is the zero section of $T^* \R^k$ except by verifying by hand that the hom spaces between microlocal skyscrapers are invariant under microlocal co-restriction. This is analogous to asking certain covariant sectorial inclusion functor on the wrapped Fukaya categories being fully-faithful. It would be interesting to develope a general method to check this, even just for piecewise linear conical Lagrangians in $\R^{2n}$. 
\end{remark}

For any $\La \In T^*M$ conical Lagrangian, we define the contact sphere bundle $T^\infty M = (T^*M \RM T_M^*M) / \R_{>0}$, and the corresponding Legendrian $\La^\infty = (\La\RM T_M^*M) / \R_{>0}$.  We define the jumping loci $\scal(\ccal_\delta)$ for $\ccal_\delta$ as the projection image of the Legendrian $SS^\infty(\ccal_\delta) \In T^\infty \R^k$ to the base $\R^k$. By the above theorem,  
$$\scal(\ccal_\delta)= \bigcup_{F_{\delta, \sigma}: F_{\delta, \sigma} \cap \Z^k \neq \emptyset} \Aff(F_{\delta, \sigma}).  $$
If $l \in \R^k$ is not in the jumping loci $\scal(\ccal_\delta)$, we say $l$ is a regular value for $\ccal_\delta$.

With the knowledge of singular support, we may deform a convex open set $U \In \R^k$ without changing the categorical output $\ccal_\delta(U)$ as along as the exterior unit conormal of $U$ is disjoint from $SS^\infty(\ccal_\delta)$. The following corolloary is immediate, noting the connected components in the complement of the jumping loci are convex hence contractible. 

\bc 
If $l_1, l_2 \in \R^k$ are in the same connected component in the complement of the jumping loci for  $\ccal_\delta$, then we have a canonical isomorphism of stalks $\ccal_\delta|_{l_1} \cong \ccal_\delta|_{l_2}$. 
\ec 

More usefully, we have the following result identifying the stalk of $\ccal_\delta|_l$ at a regular value $l$ with the constructible sheaf category on the fiber over $l$. 
\bp 
If $l$ is a regular value of $C_\delta$, then we have 
$$ \ccal_\delta|_l \cong Sh^\dm(\R^n, \La_{\delta, l}).$$
Similar statement holds if we change $Sh^\dm$ to $Sh^w$ and/or $\La_{\delta, l}$ with $\wb \La_{\delta, l}$. 
\ep
\bpf 
By Proposition \ref{p:coFF-fiber} and the triviality of $SS_{Hom}^L(\ccal_\delta)$, we may extend a sheaf in  $Sh^\dm(\mu_\R^{-1}(l), \La_{\delta, l})$ to a tubular neighborhood of $\mu_\R^{-1}(l)$, e.g. $\mu_\R^{-1}(B_x)$. By regularity of $l$, the restriction from $\mu_\R^{-1}(B_x)$ back to $\mu_\R^{-1}(l)$ is non-characteristic with respect to $SS_{Hom}(\ccal_\delta)$, hence the extension is unique, and we have proven that the restriction from $\mu^{-1}(B_x)$ to the fiber $\mu_\R^{-1}(l)$ is an equivalence of categories. 
\epf

In the following example, we illustrate what the singular support looks like. 
\begin{example} [$N=6, k=2$]
Consider the example of $(\C^*)^2$ acting on $\C^6$ with weight vectors $\beta_i$ (as column vectors) given by
$$ (\beta_1, \beta_2, \cdots, \beta_6) = \begin{pmatrix}
1 & -1 & 0 & 0  & 1 & -1\\
0 & 0 & 1  & -1 & 1 & -1
\end{pmatrix}
$$
There are 6 GKZ chambers, separated by the 6 rays generated by $\beta_i$. 

The stratification of the shift parameter space $\R^k_\delta$ (subscript is used to indicate the name of the coordinate) is shown in Figure \ref{fig:strata-delta}. We consider three sample choices of $\delta$ as shown above, with $\delta_1$ being the most non-generic and $\delta_3$ being generic. For each $\delta_i$, we illustrate in Figure \ref{fig:N6k2} the zonotope, window points, and the singular support of $\ccal_{\delta}$. Note that in the first figure, over the vertices the zonotope, we have Lagrangian cones in the cotangent fiber, marked by the blue arcs, and over other intersections of the blue hairy lines, we don't have anything extra in the cotangent fiber. 
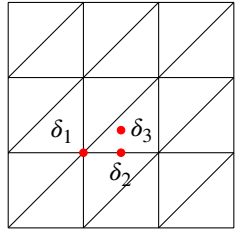
\begin{figure} [h]
    \centering
\begin{tikzpicture}
\draw (-1,-1) grid (2,2); 
\clip (-1, -1) rectangle (2,2); 
\foreach \c in {-3, -2, -1, 0, 1, 2} \draw [shift={(\c, 0)} ] (-1,-1) -- (2,2); 
\filldraw[red] (0,0) circle (0.05) node[anchor=south east, black] {$\delta_1$};
\filldraw[red] (0.5,0) circle (0.05) node[anchor=north, black] {$\delta_2$};
\filldraw[red] (0.5, 0.3) circle (0.05) node[anchor= west, black] {$\delta_3$};
\end{tikzpicture}
    \caption{Stratification of the shift parameter space $\R^k_\delta$.}
    \label{fig:strata-delta}
\end{figure}

\begin{figure}[h]
    \centering
      \begin{subfigure}[b]{0.3\textwidth}
         \centering
\begin{tikzpicture}[scale=0.5]
\clip (-4,-4) rectangle (4,4); 
\draw [gray, opacity=0.2] (-5,-5) grid (5,5); 
\foreach \x/\y in {0/0, -1/-1, 0/-1, 1/0, 1/1, 0/1, -1/0} {\draw [fill, red] (\x, \y) circle (0.15);}
\def\sx{0};
\def\sy{0};
\begin{scope} [shift={(\sx,\sy)}] 
\draw [fill = yellow, opacity=0.3] (1,1)  -- (0,1) -- (-1,0) -- (-1,-1) -- (0, -1) -- (1,0) -- cycle; 
\draw [lefthairs, blue] (-4, -5) -- (5, 4) (4,5) -- (-5, -4); 

\draw [lefthairs, blue] (-5, -1) -- (5, -1) (5, 1) -- (-5, 1); 

\draw [lefthairs, blue] (1, -5) -- (1, 5) (-1,5) -- (-1, -5); 

\def\r{0.3}
\draw [thick, blue] (-1,-1) -- +(0: \r) arc (0: 90: \r) -- cycle; 

\draw [thick, blue, rotate=180] (-1,-1) -- +(0: \r) arc (0: 90: \r) -- cycle;

\draw [thick, blue] (0,-1) -- +(45: \r) arc (45: 180: \r) -- cycle; 

\draw [thick, blue , rotate=180] (0,-1) -- +(45: \r) arc (45: 180: \r) -- cycle; 

\draw [thick, blue] (1,0) -- +(90: \r) arc (90: 225: \r) -- cycle; 

\draw [thick, blue, rotate=180] (1,0) -- +(90: \r) arc (90: 225: \r) -- cycle; 

\end{scope}
\end{tikzpicture}
         \caption{$\delta = \delta_1$}
         \label{fig:d1}
     \end{subfigure}
\begin{subfigure}[b]{0.3\textwidth}
         \centering
\begin{tikzpicture}[scale=0.5]
\clip (-4,-4) rectangle (4,4); 
\draw [gray, opacity=0.2] (-5,-5) grid (5,5); 
\foreach \x/\y in {0/0,  0/-1, 1/0, 1/1} {\draw [fill, red] (\x, \y) circle (0.15);}

\def\sx{0.3};
\def\sy{0};
\begin{scope} [shift={(\sx,\sy)}] 
\draw [fill = yellow, opacity=0.3] (1,1)  -- (0,1) -- (-1,0) -- (-1,-1) -- (0, -1) -- (1,0) -- cycle; 

\draw [dashed, blue] (-4, -5) -- (5, 4) (4,5) -- (-5, -4); 

\draw [lefthairs, blue] (-5, -1) -- (5, -1) (5, 1) -- (-5, 1); 

\draw [dashed, blue] (1, -5) -- (1, 5) (-1,5) -- (-1, -5); 

\end{scope}
\end{tikzpicture}
         \caption{$\delta = \delta_2$}
         \label{fig:d2}
     \end{subfigure}
         \begin{subfigure}[b]{0.3\textwidth}
         \centering
\begin{tikzpicture}[scale=0.5]
\clip (-4,-4) rectangle (4,4); 
\draw [gray, opacity=0.2] (-5,-5) grid (5,5); 
\foreach \x/\y in {0/0,   1/0, 1/1} {\draw [fill, red] (\x, \y) circle (0.15);}

\def\sx{0.4};
\def\sy{0.1};
\begin{scope} [shift={(\sx,\sy)}] 

\draw [fill = yellow, opacity=0.3] (1,1)  -- (0,1) -- (-1,0) -- (-1,-1) -- (0, -1) -- (1,0) -- cycle; 

\draw [dashed, blue] (-4, -5) -- (5, 4) (4,5) -- (-5, -4); 

\draw [dashed, blue] (-5, -1) -- (5, -1) (5, 1) -- (-5, 1); 

\draw [dashed, blue] (1, -5) -- (1, 5) (-1,5) -- (-1, -5); 

\end{scope}
\end{tikzpicture}
         \caption{$\delta = \delta_3$}
         \label{fig:d3}
     \end{subfigure} 

    \caption{The zonotope $\nabla_\delta$ (yellow), the windows points $W_\delta$ (red) and the singular supports (blue hairy lines and blue arcs) of the sheaf of categories on $\R^k_l \cong \R^2$. }
    \label{fig:N6k2}
\end{figure}
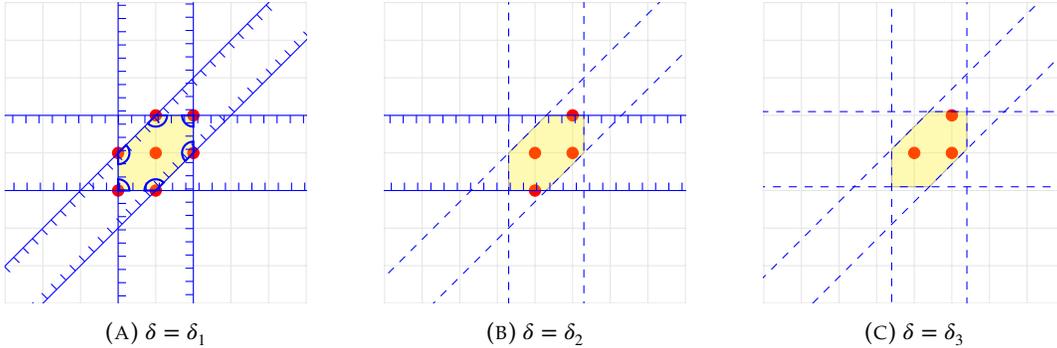
\end{example}

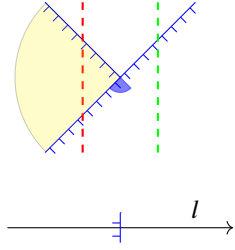
\begin{figure}[h]
    \centering
\begin{tikzpicture}
\draw [righthairs, blue] (-1,-1) -- (1, 1)  (-1,1) -- (0,0); 
\draw [fill, opacity=0.5, blue] (0,0) -- + (-45:0.2) arc (-45:-135:0.2) -- cycle; 
\draw [green, thick, dashed] (0.5, -1) -- (0.5, 1); 
\draw [red, thick, dashed] (-0.5, -1) -- (-0.5, 1); 
\draw [fill=yellow, opacity=0.2] (0,0) -- +(135: 1.4) arc (135:225:1.4); 
\begin{scope} [shift={(0,-2)}]
\draw [->] (-1.5,0) -- (1.5,0) node at (1, 0) [anchor=south] {$l$}; 
\draw [lefthairs, blue] (0, -0.2) -- (0, +0.2); 
\end{scope} 
\end{tikzpicture}       
\caption{This is a local picture of a skeleton $\La_{W_\delta}$ in $\R^2$ such that $\ccal_\delta$ has jumping loci in $\R^1_l$. The fiber skeleton over the green dashed line has one stop, and the one over the red line has two stops, a zero-dimensional Legendrian vanishing sphere $S^0$. The yellow region is the support of a sheaf that vanishes when restricted to the right of the green line, hence $SS_{Hom}(\ccal)$ has a covector pointing to the left.}
         \label{fig:1d-slices}
\end{figure}

Next, we describe what is causing the jump. Intuitively, as $l$ moves pass the jumping loci in the direction of the singular support (i.e. pass a blue hairy line in the direction of the hair), the fiber skeleton $\La_{\delta, l}$ will pick up something extra, it turns out the Legendrian boundary  $\La^\infty_{\delta, l}$ will obtain certain 'Legendrian vanishing spheres'. The prototypical local behaviors are shown in Figure \ref{fig:1d-slices}, where the fiber skeleton is 1-dimensional, and the base $\R^k_l$ is 1-dimensional. 

In general, for higher dimensional base $\R^k_l$, the jumping loci $V_\ds = \Aff(F_\ds)$ are affine linear extrapolations of the faces $F_{\delta, \sigma}$ in $\nabla_\delta$ that contain lattice points. There exist global vanishing cycles $\{ L_{\sigma, \wt v} \mid \wt v \in \Z^N, \mu(\wt v) \in V_\ds\}$ (Eq. \eqref{e:L-sigma-v}) associated to $F_\ds$, whose support in $\R^N$ has projection image in $\R^k$ contained in the affine cone
$$ C_\ds := V_\ds - \sigma^\vee = V_\ds + \R_{\geq 0} \cdot \{ y - x \mid y \in \nabla, x \in F_\sigma \}. $$  
For each $x \in V_\ds$, we can define a category of the local vanishing cycles, namely those sheaves in $\ccal_\delta(B_x)$ that restrict to 0 in $B_{x,\xi, -}$ for $\xi \in \Int(-\sigma)$.  This collection of categories again forms a sheaf of categories on $V_\ds$, denoted as $\ccal_\ds$. 
\bd 
For each $(\ds)$ such that the face $F_\ds$ contains lattice points, we define {\bf the sheaf of vanishing cycles $\ccal_\ds$ on $V_\ds$} as a sheaf of full subcategories of $\ccal_\delta$ supported on $V_\ds$, such that,    for any convex open subset $U \In V_\ds$, 
$$ \ccal_\ds(U) = \varinjlim_{\Omega \supset U} \{ F \in \ccal_\delta(\Omega) \mid \mu(\supp(F)) \In C_\ds\},  $$
where $\Omega \In \R^k$ is open and $ \Omega \cap V_\ds = U$. For small enough $\Omega$, $\ccal_\delta(\Omega)$ are all equivalent. 
\ed 
We also abuse notation and view $\ccal_\ds$ as a sheaf on $\R^k$ supported on $V_\ds$, so that for any convex open $\Omega \In \R^k$, we have $\ccal_\ds(\Omega) := \ccal_\ds(\Omega \cap V_\ds)$.

In Theorem \ref{t:decomposition}, we relate the sheaf of vanishing cycles $\ccal_\ds$ with certain window subcategories on $V_\ds$. For details, see Section \ref{s:microlocal-stalk}. Here we present a simplified version. 

\bt 
Let $F_\ds$ be a face of $\nabla_\delta$ with exterior conormal $\sigma$ which contains lattice points in $\Z^k$. Then for any $x \in \Aff(F_\ds)$ and $\xi \in \Int(-\sigma)$, we have a semi-orthogonal decomposition
$$ \ccal_\delta(B_x) = \la \ccal_{\ds}(B_x), \ccal_\delta(B_{x,\xi,-}) \ra. $$
Furthermore, $\ccal_{\ds}(B_x)$ is generated by the restriction of the following sheaves to $\mu^{-1}_\R(B_x)$, 
$$ \{ L_{\sigma, \wt v} \mid \wt v \in \Z^N, \mu_\Z(\wt v) \in F_\ds \}
$$
$$ L_{\sigma, \wt v} = \uwave{L_{\wt v}} \to \bigoplus_{I \In I_\sigma, |I|=1} L_{\wt v - e_I} \to \bigoplus_{I \In I_\sigma, |I|=2} L_{\wt v - e_I} \to \cdots, $$
where $L_{w}$ for $w\in \Z^N$ is given in Eq. \eqref{e:Lw}, and the morphisms are induced by inclusions of open sets. 
\et

\begin{example}[$N=2, k=1$] \label{ex:N2k1}
In this example, we consider the simplest possible case, where $N=2, k=1$, and $\C^*$ acts on $\C^2$ with weight $(1,-1)$. The example already illustrates many of the notions mentioned above. The zonotope $\nabla = [-1/2, 1/2]$, and the stratification of $\delta$-space consisting of $(1/2) + \Z$ in $\R$ and its complement. For a non-generic choice of $\delta$, say $\delta=1/2$, we can draw the window skeleton in $\R^N = \R^2$, as Figure \ref{fig:N2k1}.

\begin{figure}[h]
    \centering
    
\begin{tikzpicture}[scale=0.3]
\begin{scope}[rotate=45, scale=1.414]
\def\s{2} 
\draw [opacity=0.2] (-\s, -\s) grid (\s, \s); 
\clip (-\s-0.1,-\s-0.1) rectangle (\s,\s);
\foreach \c in {-2, ..., 2} {
\draw [righthairs] (\c, \s) -- (\c, \c) -- (\s, \c); 
};
\end{scope}
\begin{scope} [shift={(0,-5)}]
\draw [->] (-5,0) -- (5, 0);
\draw [dashed, red] (0, 0) -- (0, 8); 
\draw [dashed, red] (1, 0) -- (1, 8); 
\draw [thick] (0, 0.2) -- (0, -0.2) node at (0,-0.2) [anchor=north] {$0$};
\draw [thick] (1, -0.2) -- (1, 0.2) node at (1,-0.2) [anchor=north] {$1$}; 
\end{scope}
\end{tikzpicture}   
$\quad$
\begin{tikzpicture}[scale=0.5]
\begin{scope}[rotate=45, scale=1.414]
\def\s{2} 
\draw [opacity=0.2] (-\s, -\s) grid (\s, \s); 
\clip (-\s-0.1,-\s-0.1) rectangle (\s,\s);
\foreach \c in {-2, ..., 2} {
\draw [righthairs] (\c, \s) -- (\c, \c) -- (\s, \c); 
};
\draw [fill,  opacity=0.2, yellow] (-1, \s) rectangle (0, -1);
\foreach \c in {-2, ..., 2} {
\draw [righthairs, shift={(0,-1)}] (\c, \s) -- (\c, \c) -- (\s, \c); 
}
\draw [fill, opacity=0.2, green] (-1, -2) rectangle (2, -1);
\draw [very thick, red] (-1,-1) -- (0,-1);
\end{scope}
\begin{scope} [shift={(0,-5)}]
\draw [->] (-5,0) -- (5, 0);
\draw [dashed, red] (0, 0) -- (0, 8); 
\draw [dashed, red] (1, 0) -- (1, 8); 
\draw [thick, lefthairs] (0, 0.2) -- (0, -0.2) node at (0,-0.2) [anchor=north] {$0$};
\draw [thick, lefthairs] (1, -0.2) -- (1, 0.2) node at (1,-0.2) [anchor=north] {$1$}; 
\end{scope}
\end{tikzpicture}  
$\quad$
\begin{tikzpicture}[scale=0.3]
\begin{scope}[rotate=45, scale=1.414]
\def\s{2} 
\draw [opacity=0.2] (-\s, -\s) grid (\s+1, \s); 
\foreach \c in {-2, ..., 2} {
\draw [righthairs, shift={(1,0)}] (\c, \s) -- (\c, \c) -- (\s, \c); 
}
\end{scope}
\begin{scope} [shift={(0,-5)}]
\draw [->] (-5,0) -- (5, 0);
\draw [dashed, red] (0, 0) -- (0, 8); 
\draw [dashed, red] (1, 0) -- (1, 8); 
\draw [thick] (0, 0.2) -- (0, -0.2) node at (0,-0.2) [anchor=north] {$0$};
\draw [thick] (1, -0.2) -- (1, 0.2) node at (1,-0.2) [anchor=north] {$1$}; 
\end{scope}
\end{tikzpicture}    
\caption{Window skeleton of $\La_{W_\delta} \In T^* \R^N$ for $N=2$. From left to right, we have $W_\delta=\{0\}, \{0,1\}, \{1\}$ respectively. In the middle figure,  the yellow region is an example of a 'vanishing cycle sheaf' whose tip has $\mu$ image at $l=1$. Similarly, the green region is the support of a 'vanishing cycle sheaf' whose tip has $\mu$ image at $l=0$. They are the images of the structure sheaf of the unstable loci in the action of $\C^*$ on $\C^2$ with weight $(1,-1)$. See Example \ref{ex:N2k1} for details. }
    \label{fig:N2k1}
\end{figure}
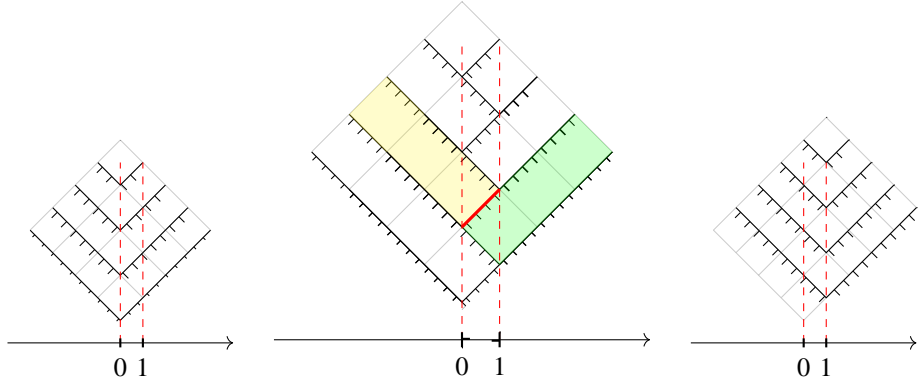

The non-generic window $W_{\delta}=\{0,1\}$ has two sub-windows, $W_{\delta-\epsilon}=\{0\}$ and $W_{\delta+\epsilon} = \{1\}$. We note that the three window skeletons only differ in the 'transition region', shown between the two red dashed lines. We denote $\wh C_- = (-\infty, 0)$  and $\wh C_+ = (1, \infty)$, the quantized GKZ chambers, and we also use $\La_\pm$ to denote the universal FLTZ skeleton, which can be determined by restrict (any of) the window skeleton $\La_{W_\delta}$ to $\mu_\R^{-1}(\wh C_\pm)$ and extend periodically in the obvious way.

Consider the window inclusion functor induced by window inclusion $\{1\} \into \{0,1\}$:
$$ \iota_+ : Sh^\dm (\R^2, \La_{\{1\}}) \into Sh^\dm (\R^2, \La_{\{0,1\}}), $$ 
the skeletons are shown in the right and middle panel in Figure \ref{fig:N2k1}.

The inclusion functors has two adjoints, the left adjoint $\iota_+^L$ and right adjoint $\iota_+^R$. The left-adjoint is the stop removal functor, where the stop removed is the difference $ \La_{\{0,1\}} \RM \La_{\{1\}}$. One such stop is shown as a red segment, whose microlocal skyscraper sheaf in $\La_{\{0,1\}}$ is shown in green. Hence, stop removal will kill the green sheaf (and its various shifts by $\ker(\mu_\Z)$ up and down), but will leave the sheaf in the 'left half space' $\mu^{-1}(\wh C_-)$ unchanged. Since if we know a sheaf in the left stable region, and the sheaf is admissble for the skeleton $\La_{\{1\}}$, then the sheaf is uniquely determined, hence the stop removal functor is the composition of restriction and co-restriction with respect to a new skeleton $\La_{\{1\}}$: 
$$ \iota^L_+: \ccal_{\delta} (\R) \to  \ccal_{\delta} (\wh C_-) \to  \ccal_{\delta+\epsilon} (\R). $$

The right adjoint $\iota_+^R$ is going to kill objects invisible by the image $\iota_+$. Since the image of $\iota_+$ is generated by the microlocal skyscrapers for cotangent fiber in open cells $\{ w + (0,1)^2 \mid w \in \mu_\Z^{-1}(1) \}$, e.g. $w + (0,1)^2$ with $w=(0,1)$,  and the stalks of the yellow sheaf in these cells are zero, the yellow sheaf and its translates by $\ker(\mu_\Z)$ are killed by the right-adjoint. \footnote{ Another way of seeing that the yellow sheaf needs to be killed is that, it is the representing object in $Sh^\dm(\R^, \La_{\{0, 1\}})$ of the following contravariant functor, 
$$ F \mapsto \Hom(F,  \tikz[anchor=base, baseline, scale=0.6, shift={(0,-0.3)}]{\draw [very thick, red] (0,0) -- (1,1); 
\draw [righthairs, blue] (0.3, 0) -- (1, 0.7); 
\draw [righthairs, blue] (0.3, 0) arc (225:45:0.50); 
} ) $$
where the red segment is the same as in Figure \ref{fig:N2k1} middle panel, and the blue hollow half-dome represents the locally constant sheaf supported in its interior with the indicated boundary. 
Hence, we call the yellow sheaf the microlocal 'co-skyscraper' for this Lagrangian disk, just as the green sheaf is the microlocal skyscraper. This operation is not so natural in a general wrapped Fukaya category, since it is cosheaf-like and one only has stop removal; however it might work if the wrapping is fully-stopped.}
By a similar argument, we have
$$ \iota^R_+: \ccal_{\delta} (\R) \to  \ccal_{\delta} (\wh C_+) \to  \ccal_{\delta+\epsilon} (\R). $$

We note that this is precisely the construction of the adjoint functors given by \cite[Lemma 6.7]{halpern2020combinatorial}. 

\end{example}

Finally, we prove homological mirror symmetry for the window categories between $A$-model and $B$-model.
\bt \label{intro-t:gen}
For any $\delta \in \R^k$, we have an equivalence of categories
$$ \bcal_\delta \cong Sh^w( \R^N / \MM, \wb \La_{W_\delta}). $$
\et 
Note that this is not automatic by coherent-constructible correspondence, which only gives a fully-faithful embedding of $\bcal_\delta$ into the A-model category. One still need to prove essential surjectivity, which is done in Theorem \ref{t:generation}.

\subsection{Window Skeleton and Perverse Schober}
Perverse schober \cite{kapranov2014perverse} on a disk is a nice combinatorial book-keeping tool to organize several spherical functors, and to localize the categorical computation to small open sets. More generally, perverse schobers with complex hyperplane arrangements \cite{kapranov2016perverse} can be defined. 
One area of motivation and appliation is the Fukaya category with coefficients, the idea of which has been adopted by \cite{nadler2019mirror, nadler2018wall} for example. 

The connection with GIT has been made in \cite{donovan2019perverse, donovan-kuwagaki-2019mirror}, and the connection with quasi-symmetric window categories has been explained in \cite[Proposition 5.1] {vspenko2019class}. There it is stated that the window inclusion functors and their right-adjoints between various window categories $\bcal_\delta$ give rises to a perverse schober. 

By Theorem \ref{intro-t:gen}, we automatically have a skeletal realization of the perverse schober, using Spenko and Van den Bergh's result. 

A slightly different point of view is to recognize $\ccal_\delta$ itself as giving a perverse schober. 

\bp 
For any $\delta \in \R^k$, $\ccal_\delta$ defines a perverse schober on $\R^k_l$ that is isomorphic to the local schober structure on $\R^k_{\delta}$ near $\delta$.
\ep 
\bpf 
On the one hand, for $\delta$ non-generic, let $S_\delta$ denote the strata in $\R^k$ that contain $\delta$, and let $\scal_\delta = \{S \mid \wb S \supset S_\delta\}$ denote the neighboring strata. On the other hand, for each affine hyperplane $H$ passing through $\delta$ in the stratification in $\R^k_\delta$, we can assign a thickened hyperplane $\wh H = \Int(\nabla_\delta) + H$, and for each strata $S' \in \scal_\delta$, there is a connected component $R_{S'}$ in the complement of the jumping loci, made from the intersections of thickened hyperplanes. 

\begin{figure} [h]
    \centering
\begin{tikzpicture}[scale=0.5]
\clip (-4,-4) rectangle (4,4); 
\draw [gray, opacity=0.2] (-5,-5) grid (5,5); 
\foreach \x/\y in {0/0,  0/-1, 1/0, 1/1} {\draw [fill, red] (\x, \y) circle (0.15);}

\def\sx{0.3};
\def\sy{0};
\begin{scope} [shift={(\sx,\sy)}] 


\draw [lefthairs, blue] (-5, -1) -- (5, -1) (5, 1) -- (-5, 1); 


\draw [green, righthairs] (-2, 0) circle (0.5); 
\draw [green, righthairs] (-2, 0) circle (1.5); 
\draw [green, rounded corners, lefthairs] (-3.8, -3.5) rectangle (3, 3.5); 

\end{scope}

\end{tikzpicture}
    \caption{An increasing sequence of open sets $U_t$ (the interior regions of the green curves). See Figure \ref{fig:d2}, where we only keep the window $W_\delta$ (red dots) and the singular support $SS(\ccal_\delta)$ (blue hairy line).  This is a non-characteristic expansion of open sets (since when the green line become tangent to the blue lines, the hairs are in opposite direction),  hence the categories $\ccal_\delta(U_t)$ are  invariant. }
    \label{fig:C-delta-U-t}
\end{figure}
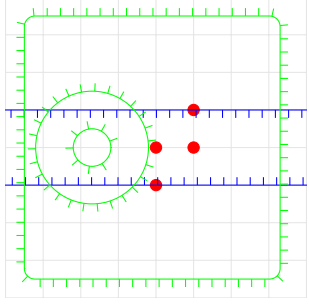

And one can verify that for any $S \in \scal_\delta$, we have 
$$ \ccal_\delta( R_S ) \cong \ccal_{S}(R_S) \cong \ccal_S(\R^k) $$
where  $\ccal_{S} = \ccal_{\delta'}$ for some $\delta' \in S$. Indeed, the first equivalence is seen by considering the window skeletons, we have $\La_{W_\delta} = \La_{W_{\delta'}}$ over $R_S$. And the second equivalence is by non-characteristic expansion from the open set $R_S$ to $\R^k$. See the Figure \ref{fig:C-delta-U-t} for an illustration of the expansion. 

And if $S_1 < S_2$ for $S_i \in \scal_\delta$, i.e. $S_1 \In \wb S_2$, we have a commuting diagram
$$ 
\begin{tikzcd}
\ccal_\delta( R_{S_2} ) \ar[r, "P_\gamma"] \ar[d, "\cong"]  & \ccal_\delta( R_{S_1}) \ar[d, "\cong"] \\
\ccal_{S_2} ( \R^k ) \ar[r, "\iota"] & \ccal_{S_1}(\R^k) \\
\end{tikzcd}
$$ 
where  the top row is realized by extension along a path $\gamma$ from a point $p_2 \in R_{S_2}$ to a point $p_1 \in R_{S_1}$, such that $\gamma$ is monotone in the sense that whenever $\gamma$ crosses a jumping locus, $\dot \gamma$ pairs positively with the covector in the singular support. We call such functors parallel transports along $\gamma$, denoted as $P_\gamma$. Then since $\ccal_\delta( R_{S_i} ) \cong \ccal_\delta|_{p_i}$, this gives the desired top arrow. The bottom row is the window inclusion functor. 
\epf

\subsection{Universal Skeleton and Local System of Categories} Finally, we discuss how to connect to the local system on the 'stringy K\"{a}hler moduli' space.  Let 
$$ \bcal = \R^k_\delta \times \R^k_l, \quad \xcal = \R^k_\delta \times \R_x^N $$
where the subscripts denote the names of the coordinate variables. Consider the trivial fibration 
$$ \wt \mu: \xcal \to \bcal, \quad (\delta, x) \mapsto (\delta, \mu(x)). $$

For any $\delta \in \R^k$, let 
\be \dcal_{\nabla_\delta} = \bigcup \{ \nabla_\delta + TF \mid \z{F is a facet of $\nabla_\delta$ containing lattice point} \}  \label{e:dcal-delta} \ee
where $TF$ is the tangent space of $F$, a hyperplane in $\R^k$. Note for generic $\delta$, $\pa \nabla_\delta$ does not contain any lattice point, and $\dcal_{\nabla_\delta} = \emptyset$. 
Consider all $\delta$ together, we define 
$$ \dcal = \bigcup_{\delta \in \R^k} \{\delta\} \times \dcal_{\nabla_\delta}. $$
Define
$$ \bcal^o = \bcal \RM \dcal, \quad \xcal^o = \wt \mu^{-1}(\bcal^o). $$

\bp 
There is a universal skeleton $\La$ define over $\xcal^o$, such that for any $(\delta, l) \in \bcal^o$, 
$$ \La|_{(\delta, l)} = \La_{\delta, l}. $$ 
\ep
\bpf 
For any $(\delta, l) \in \bcal^o$, there is a small product neighborhood $U_\delta \times U'_l$ of $(\delta, l)$, such that $l \in U'_l \In \R^k_l \RM \dcal_{\nabla_\delta}$ and $U_\delta$ is small enough that it is contained in the link of the strata that $\delta$ belongs to. One can show that, as $\delta'$ moves in $U_\delta$, $\La_{W_{\delta'}}$ remains constant over $\R^k_l \RM \dcal_{\nabla_\delta}$, hence one can define skeleton on $U_\delta \times \mu^{-1}(U'_l)$ by $T^*_{U_\delta} U_\delta \times (\La_{W_\delta}|_{U'_l})$. This universal skeleton is the unique skeleton such that when restricted to $B^o|_\delta$, it agrees with $\La_{W_\delta}$. 
\epf



\bt 
The variation of Lagrangian skeletons $\La$ for $\wt \mu: \xcal^o \to \bcal^o$ is non-characteristic, i.e., the assignment
$$ (\delta, l) \mapsto Sh^\dm(\R^n, \La_{(\delta, l)}), \quad \forall \delta, l \in \bcal^o $$
is a local system of categories. The same is true for the wrapped sheaf version and the non-equivariant version.  
\et 

Let $\Z^k$ acts on $\bcal$ diagonally, then $\dcal$ is preserved under this action. The universal skeleton is also $\Z^k$-invariant. We can consider the quotient base $\wb \bcal^o = \bcal^o / \Z^k$, with trivial $\R^n$ fibration over it,  and the family of skeleton $\La_{\delta, l}$ over the fiber $\wt \mu^{-1}(\delta, l) = \mu^{-1}(l)$. 

If $C$ is a GKZ chamber and $\wh C_\delta \In \R^k_l$ is the quantized GKZ chamber for $\nabla_\delta$, then if $l \in \wh C_\delta \cap \Z^k$, we may identify $\wb \La_{\delta, l}$ with the FLTZ skeleton for the toric GIT quotient corresponding to $C$. For any $q \in \Z^k$, the segment from $(\delta,l)$ to $(\delta+q, l+q)$ does not meet any discrimant locus. $\La_{\delta, l}$ is locally independent of $\delta$ (unless one passes through the discriminant locus), hence the local variation is the same as keeping $\delta$ and changing $l$ to $l + q$. By inspecting the universal FLTZ skeleton $\La_C$, we see this is mirror to tensoring by a line bundle labelled by the Chern character $q \in \Z^k \cong H^2(X_C)$, thus monodromies of the local system of categories are automatic. 

\bex  [{(1,-1) action of $\C^*$}]
We draw the base for the universal skeleton $\Lambda$ for the $\C^*$ action of weight (1,-1). $\Lambda$ is a 2 dimensional skeleton over the complement of the red colored slits in $\R_\delta\times \R_l$ drawn below. 
\begin{figure}[htb]
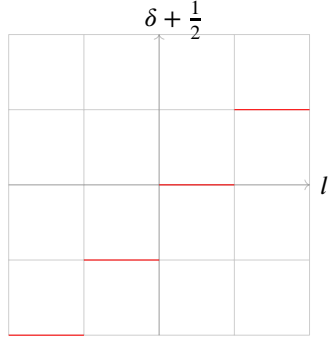

    \centering
    \tikz{
    \begin{scope};
\draw[opacity=0.2,ultra thin] (0,0) grid (4,4);
\draw[->, opacity=0.5, ultra thin, color=gray] (0,2) -- (4,2);
\draw[->, opacity=0.5, ultra thin, color=gray] (2,0) -- (2,4);
\foreach\i/\j in {0/0,1/1,2/2,3/3}
\draw [thin, color=red] (\i,\j)-- (\i+1,\j);
\node at (4.2, 2) {$l$};
\node at (2.2, 4.2) {$\delta+\frac{1}{2}$};
    \end{scope};
}

    \caption{$\R_l\times \R_\delta \RM \{[k,k+1]\times \{k\}\}_{k\in \Z}$ for the (1,-1) action of $\C^*$. If we quotient out the diagonal $\Z$-action on $\R_\delta \times \R_l$, we get a cylinder with a slit, which is homotopic with 3 punctured sphere.}
    \label{favex}
\end{figure}
\eex 

\subsection{Related works}

\subsubsection{Mirror symmetry for Toric Calabi-Yau} 
Let $P \In \{1\} \times \R^{k-1}$ be a bounded convex polytope with vertices in $\{1\} \times \Z^{k-1}$, and let $\wh P = \conv(\{0 \in \R^k\} \cup P)$. From this data we can construct a collection of toric Calabi-Yau $X_{P,T}$ depending on subdivision $T$ of $P$ as B-model, and a mirror Landau-Ginzburg A-model $W_P: (\C^*)^k \to \C$. 

On the B-side, $X_P$ and $X_{P,T}$ is defined as following. Let $\Sigma \In \R^k$ be a fan consisting of a single top dimensional cone $\sigma_P = \R_{\geq 0} \cdot P$ generated by $P$ and its faces, then $X_P$ is the toric stack associated to $\Sigma$. More generally, for any polyhedral subdivision $T$ of $P$, we can get a refinement $\Sigma_T$ of $\Sigma$, and correspondingly, we get a (partial) resolution of toric Calabi-Yau $X_{P, T} \to X_P$. The stack $X_P$ and its resolutions $X_{P,T}$ have equivalent derived categories of coherent sheaves. 

On the A-side, we may consider a $k$-variable Laurent polynomial $W_P(z_1, \cdots, z_k)$ associated with $P \cap \Z^k$ with generic coefficients of the monomials.  
Then, we may consider the wrapped Fukaya category on $(\C^*)^k$ with stop given by $W_P^{-1} (\{ w \in \C: \Re(w) \geq 1 \}).$ Thus, one get a local system of category as the coefficients of $W_P$ moves away from some discriminant locus. We may compactify $(\C^*)^k$ to a smooth projective toric variety $Y_P$ whose moment polytope is $P$, then $W_P$ compactify to a pencil $w_P: Y_P \dashrightarrow \P^1$, with base loci $B$ in the divisor $D_P$ corresponding to the facet $P \In \wh P$. The genericity  condition of $W_P$ is that we require $B$ to be smooth, and the intersection of $B$ with all the boundary strata of $D_P$ to be smooth. 
The wrapped Fukaya category is defined by removing a tubular neighborhood of all the toric boundary divisor in $Y_P$, the stop is the subset of $\arg(w_P)=0$ in the boundary of the tubular neighborhood of $D_P$. 

\bex 
(1) Let $k=1$ and $P =\{ 1 \}$. Then $X_P = \C$ and $W_P = z: \C^* \to \C$. 

(2) Let $k=2$ and $P = \conv \{ (0,1), (1,1) \}$. Then $X_P = \C^2$ and $W_P = z_1(1+z_2): (\C^*)^2 \to \C$, which after reparamerization $x=z_1, y=z_1z_2$ of $(\C^*)^2$, we have $W_P = x + y$. $Y_P = \C\P^2$, and $w_P([T_0: T_1: T_2]) = [T_0: T_1 + T_2 ]$. 

(3) Let $k=2$ and $P = \conv \{ (0,1), (n,1)\}$, then $X_P =  [C^2 / \Z_n]$ as a stack, and $W_P = z_2 ( a_n z_1^n + a_{n-1} z_1^{n-1} + \cdots + a_0)$ for generic coefficients $a_i \in \C$.  Then $Y_P$ is the weighted projective space $\P^2_{n,1,1} = [\C^3 \RM \{(0,0,0)\} /_{(n,1,1)} \C^* ]$, and  $w_P([T_0: T_1: T_2]) = [T_0 : a_n T_1^n + a_{n-1} T_1^{n-1} T_2 + \cdots + a_0 T_2^n ]$. Here $(a_i)_i$ is generic if the polynomial $a_n z_1^n + a_{n-1} z_1^{n-1} + \cdots + a_0$ has no multiple roots. The moduli space of $a_i$ is the same as configuration of $n$ distinct points in $\C^*$. 
\eex 

Borisov and Horja studied on the level of $K$-theory, the shadow of the equivalences of categories between A-model and B-models, in particular, one perform period integral on the mirror Lagrangian cycles and get GKZ hypergeometric functions. \cite{horja1999hypergeometric, borisov2006mellin, borisov2015applications}. 

The $K$-theory still loses information. For example, Seidel and Thomas shows that given a chain of exceptional collection,  the braid group acts faithfully on certain derived category of coherent sheaves \cite{seidel2001braid}, but not so on the K-theory level. The window approach for braid group action is also considered by Segal and Donovan \cite{donovan-segal-2014window}. 

\subsubsection{Windows and Quasi-Symmetric GIT}
We first recall the notion of a B-side window subcategory in GIT problem. Rouhgly speaking, if $Y$ is a GIT quotient stack of $Y=[X^{ss}_\lcal / G]$ determined by some $G$-ample equivariant line bundle $\lcal$, then one has a restriction functor on the derived category of coherent sheaves $r: Coh[X / G] \to Coh[X^{ss}_\lcal / G]$. A window subcategory for $Y$ is the image of a fully faithful embedding $Coh[X^{ss}_\lcal / G]$ back into  $Coh[X / G]$. In the case of $\C^*$-action on linear space $\C^n$ with weights $(a_1, \cdots, a_n)$, such that $\eta = \sum_{i:a_i>0} a_i = \sum_{i:a_i<0} |a_i|$, the subcategory in $[\C^n / \C^*]$ is given by $[a,a+\eta] \cap \Z$ for some $a \in \R \RM \Z$, hence the terminology of 'window'. 

The notion of window, or 'grade restriction rules' was first discovered by physicists \cite{herbst2008phases}, then imported to math by the seminal work of Ed Segal \cite{segal2011equivalences}. The work in \cite{BFK} and \cite{halpern2015derived} then shows that window subcategories exist for a general GIT setup. There the window subcategories are constructed by first stratifying the unstable locus by the Kempf-Ness (KN) strata ordered Morse theoretically via the Hilbert-Mumford function and then eliminating the obstruction to extend coherent sheaves on $X^{ss}$ to the KN strata iteratively. A choice is made in each extension step, and the final result of the composition of embeddings depends on the many choices made thus can be complicated, even for a general toric Calabi-Yau setting where $(\C^*)^k$ acts on $\C^N$ preserving the volume form.

The quasi-symmetric case is much simpler and sufficiently rich to be interesting. It is simple because one can choose the window to be defined by a polytope; it is interesting since it includes the case of symplectic representation $T^*V$ of a reductive group $G$, and the adjoint representation $\gfrak$ of $G$ is also quasi-symmetric. The quasi-symmetric case is studied by \cite{vspenko2017non, vspenko2019class} and \cite{halpern2020combinatorial}. For quasi-symmetric torus actions, a mirror calculation to the window equivalence in \cite{halpern2020combinatorial} is carried out by the first author in her Ph.D thesis \cite{jesse-thesis} using periodic FLTZ skeletons.

Another important feature of quasi-symmetry pertains to mirror symmetry expectations. It was shown in \cite{kite2017discriminants} that the complex moduli space of the A-model mirror of a quasi-symmetric GIT quotient admits a combinatorial description as the complement of the complexified $\Z^k$-periodic hyperplane arrangement parallel to facets of $\nabla$, which agrees topologically with $\bcal^o$ where the universal skeleton $\La$ lives over.


\subsection{Outline}
In section \ref{s:geom-window}, we study how the geometry of a general window skeleton in the asymptotic regions relate to the periodic skeletons for the various GIT quotients. 

In section \ref{s:rectilinear}, we study the local behavior of the window skeleton near any point $x \in \R^N$ (lattice point or not). This give rises to a class of bi-conic skeleton in $T^* \R^N$, termed rectlinear skeleton. We proved a few lemmas in preparation for later. 

In section \ref{s:SS-Hom}, we study the local obstruction for the uniqueness of extending an object in $\ccal_\delta$. This boils down to an local calculation done in \ref{s:rectilinear}. 

The section \ref{s:SS-Hom-L} is the heart of this paper, in which we proved the existence of extension of object (or the co-restriction functor is fully-faithful) for the sheaf of category $\ccal_\delta$ on $\R^k$, for any $\delta$. 

Finally, in section \ref{s:global} we study the parallel transport, the jumping loci and microlocal stalk for $\ccal_\delta$, and we show that any sheaf admissible for the window skeleton are generated by the window objects. 
\subsection{Acknowledgements}
PZ wishes to thank the KSU M-seminar organizers for inviting him to give a talk, where this work is a fruitful 'by-product'. PZ also wishes to thank M. Kontsevich and C. Diemer for generously sharing the idea of using window skeleton as interpolations. JH would like to thank her thesis advisor J. Pascaleff for useful discussions on the topic and B. Webster for inviting her to speak at Perimeter Institute where the main results in this paper were presented.


\section{Notation and Conventions}
Let $Z \In \R^N$ be a locally closed subset $Z = A \cap U$ for $U$ open and $A$ closed. We let $\C_Z$ denote the locally constant sheaf $\C_Z = j_! i_* \underline{\C}_{A \cap U}$, where $A \cap U \xto{i} U \xto{j} \R^N$. We will further abuse notation, and write $Z$ for $\C_Z$.

Let $(\CS)^k$ acts on $\C^N$ quasi-symmetrically. In other words,   the cocharacter map $\mu_\Z: \Z^N \to \Z^k$ factorizes as
$$ \bigoplus_{i=1}^m \Z^{N_i} \xto{\pi_Z = (\pi_{1,\Z}, \cdots, \pi_{m, \Z})} \Z^m \to \Z^k, $$
where $\pi_{i,\Z}$ is balanced, i.e, $\pi_{i,\Z}(1,\cdots, 1) = 0$. 

For each $i \in [m]$, let $e_{i,1},\cdots, e_{i, N_i}$ denote the basis of $\R^{N_i}$, and $\beta_{i,1}, \cdots, \beta_{i, N_i} \in \Z$ the corresponding image under $\pi_{i,\Z}$. By assumption, $\sum_{j=1}^{N_i} \beta_{i,j} = 0$. We denote $\eta_i = \sum_{j :\beta_{i,j}>0} \beta_{i,j}$ the 'window size' of $\pi_{i,\Z}$. 

We define the zonotope for $\pi_\Z: \R^N \to \R^m$, 
$$ \B = \half \pi_\R([0,1]^N) = \half \prod_{i=1}^m \pi_{i,\R}([0,1]^{N_i})  = \prod_{i=1}^m [-\eta_i/2, \eta_i/2] $$
and it a 'box'. Then $\pi'_\R(\B) = (1/2) \mu_\R([0,1]^N) = \nabla$.

For any $I \In [N]$, we define the following open regions in $\R^N$, 
$$ Q_I = \{ x \in \R^N \mid x_i < 0 \z{ if } i \in I ; x_i > 0 \z{ if } i \In I^c \} $$
$$ P_I = \{ x \in \R^N \mid x_i \in \R \z{ if } i \in I ; x_i > 0 \z{ if } i \In I^c \}. $$
And we define a Lagrangian
$$ \La_I = SS(P_I). $$

\begin{enumerate}
    \item Let $\{e_i : i \in [N]\}$ be the standard basis of $\Z^N$. For any $I \In [N]$, let $\tau_I = \cone(e_i, i \in I)$. Dually, let $\{e_i^\vee : i \in [N]\}$ be the dual basis in $(\Z^N)^\vee$. For any $I \In [N]$, let $\sigma_I = \cone(e_i^\vee, i \in I)$. 

\item For any $I \In [N]$, let $\R^I$ denote the subspace of $\R^N$ generated by $e_i, i \in I$, and $[0,1]^I$ be the unit cube in $\R^I$. 

\item For any subset $A \In \R^k$, we denote $\wt A = \mu_\R^{-1}(A)$; and if $A \In \Z^k$, we overload the notation and denote $\wt A = \mu_\Z^{-1}(A)$.


\item If $\La \In T^*M$ is a Lagrangian skeleton, $U \In M$ an open set, we let $\La|_U = \La \cap T^*U$ be the restriction of $\La$ to $U$. If we are given a map $f: M \to N$, and $V \In N$ an open set, we sometimes abuse notation and write $\La|_V  := \La|_{f^{-1}(V)}$. 

\item A sign is an element in $\{+, -, 0\}$. 
An $n$-dimensional sign vector is an $n$-tuple $s=(s_1, \cdots, s_n) \in \{+,-,0\}^n$. We endow $\{+,-,0\}$ with a partial ordering
$ 0 < +, \quad 0 < -. $
\end{enumerate}

If $s \in \signset$, we define cone $\sigma_s \In \R^\vee$ by
$$ \sigma_+ = \R_{\geq 0}^\vee, \quad \sigma_- = \R_{\leq 0}^\vee, \quad \sigma_0 = \{0\}. $$
For $s=(s_1, \cdots, s_n) \in \signset^n$, we have cone 
$$ \sigma_s=\sigma_{s_1} \times \cdots \times \sigma_{s_n} \In (\R^\vee)^n.$$ 
We define $\tau_s$ in $\R^n$ similarly. 
We define the obvious map
$$ \sign: \R \to \{+,-,0\}. $$
And we generalize it to $\R^\vee$ and $\R^n, (\R^n)^\vee$. 

All categories used dg derived categories, and functors are derived functors.






\section{Geometry of the Window Skeleton: Stable Region\label{s:geom-window}} 
Here we consider $(\CS)^k$ acting on $\C^N$ such that the cocharacter maps fits into short exact sequences \eqref{eq:Mseq} and \eqref{eq:Nseq}, without any further assumptions.

For any $I \In [N]$, we define the $I$-box as the unit cube $[0,1]^I \In \R^I$, and its image under $\mu_\R$ the $I$-zonotope
$$ \nabla_I =  \mu_\R([0,1]^I) = \sum_{i \in I} [0, \beta_i]. $$

\subsection{GKZ fan and GIT skeleton}

The closed positive quadrant $\R_{\geq 0}^N$ and its faces (in various dimensions) project to $\R^k$, and their intersections induces a fan structure on $\R^k$, called {\it GKZ fan} $\Sigma_{GKZ}$. We call the interior of the top dimensional cones in $\Sigma_{GKZ}$ as chambers. 

Let $C$ be a GKZ chamber, and $p \in C$ any point. Then the fiber $\mu_\R^{-1}(p)$ intersects all the faces of $\R_{\geq 0}^N$ transversally if the intersection is non-empty.  We define
$$ \ical_C = \{ I \In [N] \mid \tau_{I} \cap \mu_\R^{-1}(p) \neq \emptyset \}. $$ 
Then $\ical_C$ satisfies the following property, if $I \in \ical_C$ and $J \supset I$, then $J \in \ical_C$. Recall that $C_I = \mu_\R(\tau_I)$, hence $I \in \ical_C$ if and only if $C_I \supset C$. 

There is a sub-fan $\wt \Sigma_C \In \Sigma_{\C^N}$ in the fan of $\C^N$, given by 
$$ \wt \Sigma_\C = \{ \sigma_I: I^c \in \ical_N \}.$$

The image of $\wt \Sigma_C$ under the map $\nu_\R$ defines a simplicial stacky fan $\mathbf{\Sigma}_C$ in $\NN_\R$. More precisely, the stacky fan is given by 
\begin{enumerate}
    \item a simplicial fan $\Sigma_C = \{\wb \sigma_I: I \in S_C\}$ where $\wb \sigma_I = \nu_\R(\sigma_I)$; 
    \item and for each ray $\rho_i \in \Sigma_C$, an integral vector $\alpha_i \in \rho_i$. 
\end{enumerate}

We then have corresponding (universal) FLTZ-skeletons $\La_{\wt \Sigma_C} \In T^* \R^N$ and $\La_{\Sigma_C} \In T^* \MM_\R \cong T^* \R^n$, where
\be \label{eq:LaC} \La_C = \La_{\wt \Sigma_C} = \bigcup_{I \in S_C} \{x \in \R^N \mid x_i = \la x, e_i^\vee \ra  \in \Z, \,  \forall i \in I \} \times (-\sigma_I) \In \R^N \times (\R^N)^\vee = T^* \R^N \ee
and 
\be \La_{\mathbf{\Sigma}_C} = \bigcup_{I \in S_C} \{\wb x \in \MM_\R \mid \la \wb x, \alpha_i \ra \in \Z, \forall i \in I\}  \times (-\wb \sigma_I) \In \MM_\R \times \NN_\R = T^* \MM_\R. \ee

We recall the definition of the specialization of a conical Lagrangian $L \In T^*X$ to a submanifold $S \In X$ transverse to $L$, i.e,  $T_S^* X \cap L \In T_X^* X$: 
$$ L|_S = (T^*X|_S \cap L) / T_{S}^* X \In T^*X|_{S}  / T_{S}^* X \cong T^* S. $$
If $f: X \to Y$ is a smooth submersion, and $L$ is transverse to all the fibers of $f$, then we say $L$ is transverse to $f$, or $f$ is non-characteristic to $L$. 

\bp
$\La_{\wt \Sigma_C}$ is transverse to $\mu_\R: \R^N \to \R^k$. And for any $v \in \Z^k$ and $\wt v \in \mu_\Z^{-1}(v)$, if we identify $(\mu_\R^{-1}(v), \wt v)$ with $(\MM_\R, 0)$, then the specialization $\La_{\wt \Sigma_C}$ to the fiber $\mu_\R^{-1}(v)$ is the same as $\La_{\mathbf \Sigma_C}$ on $T^*\MM_\R$. 
\ep
\bpf
The fiber of conormal bundle to $\mu_\R^{-1}(p)$ can be identified with cotangent fiber $T_p^*\R^k \cong (\R^k)^\vee$. Since by construction, $\sigma_I$ maps to $\wb \sigma_I$ bijectively, hence $\sigma_I \cap (\R^k)^\vee = 0$. 

To see the matching of Lagrangians, it suffices to note that if $x \in \mu_\R^{-1}(v)$, then for any $i \in [N]$ $$ \la x, e_i \ra \in \Z \LRA \la x-\wt v, e_i \ra \in \Z \LRA \la x-\wt v, \alpha_i \ra \in \Z$$
where in the last step $x - \wt v \in \MM_\R$, and  $\la x-\wt v, e_i \ra = \la x-\wt v, \alpha_i \ra $. 
\epf

\subsection{Window skeleton over stable regions}

Let $P \In \R^k$ be a bounded convex polytope, $W = P \cap \Z^k$.  
For each GKZ chamber $C$, we will compare
$$ \La_W = \bigcup_{w \in \wt W} \La_w, \quad \z{ and }  \La_C = \La_{\wt \Sigma_C}. $$
We will show that if $P$ is large enough (see Definition \ref{d:P-large}), then the two Lagrangians agrees over certain asymptotic (i.e. stable) regions  in the direction of $C$.

First we define the notion of a 'large enough' polytope. Let $H$ be a generic (closed) half space in $\R^k$, 'generic' meaning  $\pa H$ does not contain any $\beta_i$. Then we define the zonotope for $H$ by 
$$ \nabla_H = \sum_{\beta_i \in H} [0, \beta_i]. $$
There are only finitely many possible such zonotopes as $H$ varies.

\bl 
In the quasi-symmetric case, there is only one $\nabla_H$ up to translation, which is exactly the zonotope $\nabla$ defined before. 
\el 
\bpf 
Let $\rcal$ be the collection of rays in $\R^k$ containing some $\beta_i$. For each $\rho \in \rcal$, let $\beta_\rho = \sum_{\beta_i} \beta_i$. By quasi-symmetric condition, if $\rho \in \rcal$, then $-\rho \in \rcal$, and $\beta_\rho = -\beta_{-\rho}$. Since for any generic choice of half space $H$, one and exactly one of $\pm \rho$ is in $H$, hence the resulting Minkowski sum is the same up to translation. 
\epf

Let $A, B$ are subsets of $\R^N$, we say '$A$ fits in $B$' and  '$B$ contains a translate of $A$', if there exists $s \in \R^N$ such that $s + A \subset B$. We denote this by $A \sqsubset B$.

\bd \label{d:P-large}
Let $P$ be a polytope. We say $P$ is large enough, if for all generic closed half space $H$, $\nabla_H$ fits in $P$. 

\ed


Next, we define the asymptotic region in the direction of $C$.
$$ C_P := \bigcap_{p \in P} p + C = \{ x \in \R^k \mid P \subset x - C \}. $$

The rest of the section is devoted to prove the following theorem
\bt \label{t:LaWLaC}
Let $P$ a large enough bounded convex polytope, and $W = P \cap \Z^k$. Then for any GKZ chamber $C$,  $\La_W$ and $\La_C$ agree over $C_P$, i.e. $\La_W|_{C_P} = \La_C|_{C_P}$. 
\et 

This has an immediate corollary. 
\bc 
Assume we have the quasi-symmetric condition. For any $\delta \in \R^k$ and any GKZ chamber $C$, the window skeleton $\La_{W_\delta}$ and the GIT skeleton $\La_C$ agree over $C_{\nabla_\delta}$. 
\ec

\ss{Proof of Theorem \ref{t:LaWLaC}}
Recall the definitions of $\La_W$ and $\La_C$: 
$$ \La_W = \bigcup_{w \in \wt W} \La_w = \bigcup_{w \in \wt W} \left( w + (\bigcup_{I \In [N]} \tau_I \times (-\sigma_{I^c})) \right) $$
and 
$$ \La_C = \Z^N + \bigcup_{I \in \ical_C} ( \R^I) \times (-\sigma_{I^c}), $$
where the addition is done in the base direction $\R^N \In T^* \R^N$.

\sss{} We first show that $\La_W \In \La_C$ over  $C_P$. This follows from the following lemma. 
\bl
For any $w \in \Z^N$, $\La_w \subset \La_C$ over $ \wb w + C$. \el
\bpf
Since $\La_C$ is translation invariant by $\Z^N$, it suffices to check the case $w=0$. Since $\tau_I \cap \wt C \neq \emptyset$ if and only if $I \in \ical_C$, and since $\tau_I \In \R^I$, we have
$$ \La_0 \cap \wt C \In \bigcup_{I \in \ical_C} ( \R^I) \times (-\sigma_{I^c}) \In \La_C. $$
\epf
Since for each $w \in \wt W$, $\wb w + C \supset C_P$, hence $\La_W \In \La_C$ over $C_P$. 

\sss{} 
We then show that $\La_W \supset \La_C$ over $C_P$.


\bl\label{l:P_nabla}
If $P$ is large enough , then for any $I \In [N]$ such that $C_I$ is a proper cone containing $C$,  $P$ contains a translate for all the zonotope $\nabla_I$. 
\el 
\bpf
If $C_I$ is a proper cone, then by definition it is contained in a generic half space $H$. And we have $C \In C_I \In H$. Let $J = \{ i \in [N]: \beta_i \in H \}$, then $C_J$ is the maximal proper cone in $H$, and $\nabla_H = \nabla_J$.  Since $I \In J$, $\nabla_J \supset \nabla_I$. Hence, up to a translate, $P$ contains $\nabla_J$ hence $\nabla_I$. 
\epf

\bpf [Proof of the theorem \ref{t:LaWLaC}].  For any $x \in \wt C_P$, we have
$$ \La_C|_x = \cup\{ -\sigma_{I^c} \mid I \in \ical_C, x \in \Z^N + \R^I \}, \quad \La_W|_x = \cup\{ -\sigma_{I^c} \mid I \In [N], x \in \wt W + \tau_I \} $$

Hence, to show $\La_C|_x \In \La_W|_x$, it suffices to show that for any $I \in \ical_C$ such that $x \in \Z^N + \R^I$, we have $w \in \wt W$, such that $x \in w + \tau_I^o$.

{\bf case (a):} $C_I$ is a proper cone. Let $x \in a + \R^I$ for some $a \in \Z^N$. Define the affine space $\R_I = a + \R^I$, and affine lattice $\Z_I = a + \Z^I$. Since $C_I$ is a proper cone, by Lemma \ref{l:P_nabla},  $P$ contains a translate of $\nabla_I$. Since $\nabla_I$ and $-\nabla_I$ differ by a translation, $P \sqsupset -\nabla_I$, i.e there exists $s \in \R^k$, such that 
$$ s - \nabla_I \In P \In \wb x - C_I. $$
Lift $s$ to $\wt s \in x - \tau_I^o$, and since $- [0,1]^I \In - \tau_I$,  we have 
$$ \wt s - [0,1]^I \In x - \tau_I^o. $$
Intersecting both with $\Z_I$, we have $(\wt s - [0,1]^I) \cap \Z_I \neq \emptyset$. Let $w \in (\wt s - [0,1]^I) \cap \Z_I \In \wt P \cap \Z^N = \wt W$, we have $w \in x - \tau_I^o$, or $x \in w + \tau_I^o$. 

{\bf case (b):} $C_I$ is not a proper cone. Let $J \In I$ be a maximal subset such that $C_J$ contains $C$ and is proper. Then for any $i \in I \RM J$, $\beta_i \in - C_J$, since otherwise, one can add $i$ to $J$, and keep $C_{J \cup \{i\}}$ as a proper cone, contradicting with the maximality of $J$. Let $l = \sum_{i \in I \RM J} \beta_i$, then $l \in -C_J$ hence $l + C_J \supset C_J$. 

Recall that we need to find a $w \in \wt W$, such that $w \in x - \tau_I$. Define $y = x - \sum_{I \RM J} e_i$, then $\wb y = \wb x - l \in \wb x + C_J$. Hence, suffice to find a $w \in \wt W$, such that $w \in y - \tau_J^o$, which reduces to case (a). 
\epf 

\section{Rectilinear Skeletons \label{s:rectilinear}}
In this section, we setup some notations and establish some basic properties of constructible sheaves in Euclidean space adapted to the coordinate hyperplane stratifications. 

Let $N$ be any positive integer, $[N]=\{1,2,\cdots, N\}$. Let $\{e_i\}$ be the standard basis of $\R^N$. For $I \In [N]$, we define $e_I = \sum_{i \in I} e_i$, $\R^I = \span_\R \la e_i, i \in I \ra$, $\tau_I = \cone\la e_i, i \in I \ra$. Let $(\R^N)^\vee$ be the dual space of $\R^N$, with dual basis $e_i^\vee$. For $I \In [N]$, we define $\sigma_I = \cone(e_i^\vee ,i \in I)$. 

Consider a stratification $\scal_N$ of $\R^N$ generated by the coordinate hyperplanes and their intersections, where the largest strata are $n$-dimensional open quadrants in $\R^N$. We call $\scal_N$ {\em quandrants stratification} of $\R^N$. Let $\sign: \R \to \{+,0,-\}$ be the sign map. For each strata $s \in \scal_N$, we have a sign vector $\sign(s) \in \{+,0,-\}^N$, and we sometimes write $s_i = \sign(s)_i$ the i-th component of the sign vector. 

Recall that we defined two types of open sets in $\R^N$:  quadrants $Q_I$ and wedge $P_I$
$$ Q_I: \bcs x_i < 0 & \z{ if } i \in I \\
x_i > 0 & \z{ if } i \notin I \ecs, \quad  \z{ and } \quad 
P_I: \bcs x_i \z{ free } & \z{ if } i \in I \\
x_i > 0 & \z{ if } i \notin I \ecs.  $$ 
We also use $P_I$ and $Q_I$ to denote constructible sheaves in $Sh_{\scal_N}(\R^N)$ with non-zero stalk $\C$ over $P_I$ and $Q_I$ respectively. 

We define a skeleton $\La_N \In T^*\R^N$ by 
$$ \La_N = \bigcup_{I \In [N]} \R^I \times (-\sigma_{I^c}). $$
Then $\La_N$ is adapted to the stratification $\scal_N$. Then the sheaf $P_I$ has $\La_I = SS(P_I) \In \La_N$, and in fact $\La_N$ can be written as 
\be \label{e:Lan} \La_N = \bigcup_{I \In [N]} \La_I. \ee 

Let $\pcal_n = \{I \In [N]\}$ be the power set of $[N]$. For any index subset $\ical \In \pcal_n$, we have a sub-skeleton 
$$ \La_\ical = \bigcup_{I \in \ical} \La_I \In \La_N. $$
We call $\La_\ical$ a {\bf rectilinear skeleton in $\R^N$}. 
We also define the sheaf categories $\ccal_\ical = Sh^\dm(\R^N, \La_\ical)$, and $\ccal^w_\ical$ the wrapped constructible sheaves, i.e. compact objects in $\ccal_\ical$. If $\ical = \pcal_n$, we omit the subscript $\ical$. 


\subsection{Probe sheaves for stalks on quadrants \label{ss:probe-general} }
Fix an index set $\ical$, and consider category $\ccal_\ical$.  For any $I \In [N]$, we consider the stalk functor $$ \Phi_{I, \ical}: \ccal_\ical \to \Vect, \quad G \mapsto \Hom(Q_I, G) = \z{ stalk of $G$ at $x_I \in Q_I$ }. $$
As explained in \cite{nadler2016wrapped}, the functor $\Phi_{I, \ical}$ admits left adjoint, hence has co-representing objects $F_{I, \ical} \in \ccal^w_\ical$, i.e. 
$$ \Hom(F_{I, \ical}, -) = \Phi_{I, \ical}(-): \ccal_\ical \to \Vect. $$
We $F_{I, \ical}$ the {\em probe sheaf} for quadrants $Q_I$ with repect to skeleton $\La_{\ical}$. 

We first study the category $\ccal$. 
\bp \label{p:LaFull}
(1) For each $I \in \pcal_n$, the probe sheaf 
$$ \fcal_{I, \pcal_n} = P_I. $$
That is, for any $G \in \ccal$, we have $\Hom(Q_I, G) \cong \Hom(P_I, G)$. 

(2) For each $I \in \pcal_n$, $P_I$ is a projective object in $\ccal$, in the sense that $\Hom(P_I, -)$ is an exact functor.

(3) $\{P_I, I \in \pcal_n\}$ forms a compact generator in $\ccal$, in the sense that if $G \in \ccal$ such that $\Hom(P_I, G) \cong 0$ for all $I \in \pcal_n$, then $G \cong 0$. 
\ep 
\bpf 
(1) We can define a non-characteristic deformation of open sets from $Q_I$ to $P_I$ with respect to $\La_N$. For example, for $ t \geq 0$, let $P_I(t) = Q_I - t e_I$, then $P_I(0) = Q_I$, $P_I = \cup_{t =0}^\infty P_I(t)$, and $\Hom(P_I(t), G)$ is independent of $t$. 

(2) This follows since the stalk functor is an exact functor. 

(3) Suppose $\Hom(P_I, G) = 0$ for all $I \in \pcal_n$, then $G$ has vanishing stalk in all open quadrants. Then $\supp(G)$ is in the union of coordinate hyperplane, hence $SS(G)$ contains conormal of some strata $s \in \scal_N$ which cannot be contained in $\La_N$. Hence $\supp(G) = \emptyset$ and $G=0$. 
\epf 

\bc \label{c:quiver}
The category $\ccal$ is equivalent to the category of presheaves on $\pcal_n$
$$ \Phi: \ccal \congto Psh(\pcal_n) := \Fun(\pcal_n^{op}, \Vect), \quad G \mapsto (I \mapsto \Hom(P_I, G)), $$
where $\pcal_n$ is viewed as a partially ordered set and hence a category.
\ec 
\bpf 
The full subcategory of $\ccal_n$ with objects $\{P_I, I \in \pcal_n\}$ is equivalent to the $\C$-linearization of $\pcal_n$ where if $I < J$ in $\pcal$, we define $\Hom_{\pcal_n}(I, J) = \C$. Then the result follows since $\{P_I, I \in \pcal_n\}$ compactly generates $\ccal$. 
\epf 


Next, for any $\La_\ical$, we consider the closed subskeleton $\La_\ical \In \La_N$ and full subcategory $\iota_{\ical}: \ccal_\ical \into \ccal$. We have analogous results. 

For any $I \In [N]$, if there is $J \in \ical$ such that $I \In J$, we say $I$ is dominated by $\ical$. 

\bp \label{p:LaPartial}
(1) For any $I \in \ical$, the probe sheaf 
$$ F_{I, \ical} = P_I. $$
That is, for any $G \in \ccal_\ical$, we have $\Hom(Q_I, G) \cong \Hom(P_I, G)$. 

(2) For any $I \In [N]$, if $I$ is not dominated by $\ical$, then $F_{I, \ical} = 0$.

(3) For any $I \In [N]$,  $I \notin \ical$ and $I$ is dominated by $\ical$, then the following complex is acyclic
$$ F_{I, \ical}  \to  \bigoplus_{J \In I^c, |J|=1} F_{I \sqcup J, \ical} \to \bigoplus_{J \In I^c, |J|=2} F_{I \sqcup J, \ical} \to \cdots \to F_{[N], \ical} $$ 

(4) $\{P_I, I \in \ical\}$ forms a compact generator in $\ccal_\ical$, in the sense that if $G \in \ccal_\ical$ such that $\Hom(P_I, G) \cong 0$ for all $I \in \ical$, then $G \cong 0$. 

(5) The category $\ccal_\ical$ is equivalent to the category of presheaves on $\ical$
$$ \Phi: \ccal_\ical \congto Psh(\ical) := \Fun(\ical^{op}, \Vect), \quad G \mapsto (I \mapsto \Hom(P_I, G)). $$
\ep 
\bpf 
(1) The proof is the same as Proposition \ref{p:LaFull}(1). 

(2) If $I$ is not dominated by $\ical$, then the skeleton $\La_\ical$ is not supported on $Q_I$, hence $F_{I, \ical} = 0$.

(3) If  $s = -\Int(\tau_I)$ is a strata of $\scal_N$, then $\Int(\La_N|_{s}) \cap (\La_{\ical})_s = \emptyset$, since only $\La_I$ contribute to $\Int(\La_N|_{s})$. Hence we have an exact sequence of stalk functors
$$ \Phi_{I, \ical}  \gets  \bigoplus_{J \In I^c, |J|=1} \Phi_{I \sqcup J, \ical} \gets \bigoplus_{J \In I^c, |J|=2} \Phi_{I \sqcup J, \ical} \gets \cdots \gets \Phi_{[N], \ical}. $$ 
Translate the above relation to co-representing objects, we get the desired claim. 

(4) By (3) and induction on $|I|$ from $n$ to $0$, we see for each  $I \In [N]$,   $F_{I, \ical}$ can be expressed as a finite complex using $\{P_J, J \in \ical\}$. Since $\{F_{I, \ical}, I \in [N]\}$ forms a set of compact generators of $\ccal_\ical$, hence $\{P_I, I \in \ical\}$ also generates. 

(5) The proof is the same as in Corollary \ref{c:quiver}. 
\epf 

We have seen that $\{ \fcal_{I, \ical}, I \In [N]\}$ can be expressed using $\{P_I, I \in \ical\}$. The following is a more explicit (but maybe huge) formula. Let $\ical_{\geq I} = \{ J: J \in \ical, J \geq I \}$. We recall the following definition. 

\bd \label{d:order-complex}
For any poset $P$, we define a simplicial complex $\Delta(P)$, where a $k$-chain is  of the form $[p] = (p_0 < p_1 < \cdots < p_k)$ in $P$. We call $\Delta(P)$ the  order complex of $P$. 
\ed 

\bp \label{p:probe-as-holim}
Let $\Delta = \Delta(\ical_{\geq I})$ be the order complex of $\ical_{\geq I}$ and let $[J]=(J_0<\cdots<J_k)$ denote an element in $\Delta_k$.  Then the probe sheaf $F_{I, \ical}$ for quadrant $Q_I$ in $\La_{\ical}$ can be resolved as
\be F_{I, \ical} \cong  \bigoplus_{[J] \in \Delta_0} P_{J_0} \to \bigoplus_{[J] \in \Delta_1 } P_{J_1} \to \cdots. \label{e:probe-as-holim} \ee
\ep 
\bpf 
We will prove that $$ F_{I, \ical} \cong \int_{J \in \ical } \Hom(Q_I, P_J)^\vee \ot P_J  \cong \holim_{J \in \ical_{\geq I}} P_J $$
then the desired result follows as a concrete realization of homotopy limit over a poset.

Let $\wh \ical = Psh(\ical)$ and $h: \ical \into \wh \ical$ be the Yoneda embedding. Then $h(J) = P_J$ under identification of $Psh(\ical) \cong \ccal_\ical$. For any $I \In [N]$, we have a functor $\varphi_{\geq I}: \ical \to \Vect$, sending $J \in \ical$ to $\C$ if $J \supset I$ and to $0$ if $J \not\supset I$. Then, we can realize the left-adjoint $\Phi_{I, \ical}^L$ as a right Kan extension (see the first diagram below). Since $\ical$ generates $\wh \ical$, one can compute the right Kan extension using the second diagram below. 
\[ 
\begin{tikzcd}
& \Vect \ar[rd, dashed, "\Phi_{I, \ical}^L"] & \\
\wh \ical \ar[rr, "Id"] \ar[ru, "\Phi_{I, \ical}"] & & \wh \ical
\end{tikzcd}
\quad 
\begin{tikzcd}
& \Vect \ar[rd, dashed, "\Phi_{I, \ical}^L"] & \\
\ical \ar[rr, "h"] \ar[ru, "\varphi_{\geq I}"] & & \wh \ical
\end{tikzcd}
\]
Hence, we get the co-representing object as
\[ F_{I, \ical} = \Phi_{I, \ical}^L(\C)= \int_{J \in \ical} (\Hom(\C, \varphi_{\geq I}(J))^\vee \otimes h(J) =  \holim_{J \in \ical_{\geq I}} P_J. \]

\epf 

Recall that $\iota_\ical: \ccal_\ical \to \ccal$ is the inclusion functor induced by the closed embedding of skeleton $\La_\ical \In \La_N$. It has a left adjoint 
$$ \iota_\ical^L: \ccal \to \ccal_\ical. $$
Since for any $G \in \ccal_{\ical}$, we have 
$$ \Phi_{I, \ical}(G) = \Hom(Q_I, \iota_{\ical}(G)) \cong \Hom(P_I, \iota_{\ical}(G)) = \Hom(\iota_{\ical}^L P_L, G) $$
we have
$$ F_{I, \ical} = \iota_{\ical}^L(P_I). $$
\begin{remark}
Geometrically, $\iota_{\ical}^L$ is obtained by wrapping using geodesic flow stopped by $\La_\ical$. Here we uses the (covariant) equivalence between wrapped Fukaya category on $T^*\R^N$ and constructible sheaf on $\R^N$, where the sign convention of $\omega$ identfies the negative Reeb flow on the Fukaya side with the geodesic flow on the constructible sheaf side. 
\end{remark}

\subsection{Support of probe sheaves and topology of poset}
In this section, we study the support of the probe sheaf, or equivalently, the hom between probe sheaves. This turns out to be related to the topology of the poset $\ical$ and its various subsets $\ical_{\geq I}$. 

Let $P$ be a finite poset,  $\Delta(P)$ be its order complex (Definition \ref{d:order-complex}) and $|\Delta(P)|$ be its geometric realization. 

We define the cochain complex of a partially ordered set $P$ as
$$ C^0(P) \to C^1(P) \to \cdots $$ 
where $C^k(P) = \C^{\Delta(P)_k} = \z{Map}(\Delta(P)_k, \C)$ is the $\C$-vector space with basis in $\Delta(P)_k$, and the differential is pull-back along the face map. And we define the cohomology $H^\bu(P)$ as the cohomology of the above cochain complex. Clearly, we have
$$ H^\bu(P, \C) = H^\bu(|\Delta(P)|, \C). $$

If $P \to Q$ is a map of posets, then we have the corresponding maps of its geometric realizations $|\Delta(P)| \to |\Delta(Q)|$. In particular, if $P \into Q$ is an inclusion, then $|\Delta(P)| \to |\Delta(Q)|$ is a closed embedding. 

\bp \label{p:hom-of-subposets}
Let $P$ be a poset and $P_1, P_2 \In P$ be sub-posets. Then the following definition of $\Hom(P_1, P_2) \in \Vect$ are equivalent. 
\begin{enumerate}
    \item Consider the following cochain complex 
    $$ C^k(P_1, P_2) = \bigoplus_{[p] \in \Delta(P_2)_k} \delta_{P_1}([p]), \quad \delta_{P_1}([p]) = \bcs \C & [p] \cap P_1 \neq \emptyset \\ 0 & \z{ otherwise } \ecs $$
    The differential is defined in the usual way. We define $$ \Hom_\ical(P_1, P_2) = C^*(P_1, P_2). $$
    \item For $i=1,2$, consider the locally constant sheaf $\C_{Z_i}$ supported on the closed subset $Z_i = |\Delta(P_i)|$ in $Z = |\Delta(P)|$. We define 
    $$ \Hom_\ical(P_1, P_2) := \Hom^*(\C_{Z_1}, \C_{Z_2}) = C^*( Z_2, Z_2 \cap (Z_1)^c) $$
\end{enumerate}
\ep 
\bpf 
Let $U_1 = Z \RM Z_1$.  Since we have short exact sequence
$$ 0 \to \C_{U_1} \to \C_{Z} \to \C_{Z_1} \to 0 $$
hence we have 
$$ \Hom(\C_{Z_1}, \C_{Z_2}) \cong \cone( \Hom(\C_Z, \C_{Z_2}) \to \Hom(\C_{U}, \C_{Z_2}) $$
We may realize $\Hom(\C_Z, \C_{Z_2})$ as $C^*(P_2)$ and $\Hom(\C_{U}, \C_{Z_2})$ as $C^*(P_2 \RM P_1)$ i.e cochain in $P_2$ that avoids $P_1$. Then we have the following short exact sequence of cochain complexes
$$ 0 \to C^*(P_1, P_2) \to C^*(P_2) \to C^*(P_2 \RM P_1) \to 0. $$
\epf 

Fix an index subset $\ical$. Let $I, J \In [N]$. We have
\bp \label{p:hom-FI-FJ}
The stalk of $F_{J, \ical}$ in quadrant $Q_I$ is given by 
$$ \Hom(F_{I, \ical}, F_{J, \ical}) \cong \Hom_{\ical}(\ical_{\geq I}, \ical_{\geq J}). $$
\ep 
\bpf
$$ \Hom(F_{I, \ical}, F_{J, \ical}) = \Hom(Q_I, F_{J, \ical}) \cong \bigoplus_{[L] \in \Delta(\ical_{\geq J})_0} \Hom(Q_I, P_{L_0}) \to \bigoplus_{[L] \in \in \Delta(\ical_{\geq J})_1 } \Hom(Q_I, P_{L_1}) \to \cdots $$
Since $\Hom(Q_I, P_{L_k}) = \delta_{ \ical_{\geq I}}([L])$, the chain complex is the same as $C^*(\ical_{\geq I}, \ical_{\geq J})$. 
\epf

\bc \label{c:skeleton-with-floor}
If $[N] \in \ical$, then $\Hom(F_{I, \ical}, F_{I, \ical}) = \C$ for all $I \In [N]$. 
\ec 
\bpf 
Since $\ical_{\geq I}$ is non-empty and has a final element $[N]$, its geometric realization is contractible. Hence $H^*(|\Delta(\ical_{\geq I})|)=\C$. 
\epf

\subsection{Hourglass sheaf and microlocal skyscraper\label{ss:hourglass}} 
Let $\ical = \pcal_N \RM \{ \emptyset \}$, then the skeleton $$ \La_{\hourglass, N} := \La_\ical = \bigcup_{\emptyset \neq I \In [N]} \R^I \times (-\sigma_{I^c}) \In T^* \R^N. $$ For any $\emptyset \neq I \In [N]$, the wedge sheaf $P_I$ has its singular support $\La_I \In  \La_{\hourglass,N}$, hence the hourglass sheaf $\hourglass_N$ is admissible for $\La_{\hourglass,N}$. However, over the origin $0 \in \R^N$, the open cone $-\sigma_{[N]}^o \In T_0^* \R^N$ is absent from $\La_{\hourglass_N}|_0$, thus the sheaf $P_{\emptyset}$ on the open quadrant is not admissible.

For $[N] \supset I \neq \emptyset$, we have $F_{I, \ical} = P_I$. For $I = \emptyset$, we have the probe sheaf
\[ \hourglass_N: = F_{\emptyset, \ical} = \uwave{\bigoplus_{I \In [N], |I|=1} P_{I}} \to \bigoplus_{I \In [N], |I|=2} P_{I} \to \cdots \to P_{[N]}. \]
This sheaf has its support only in two quadrants, the totally positive open quadrant $Q_\emptyset$ and the totally negative closed quadrants $\wb{Q_{[N]}}$ (hence the name 'hourglass'). More precisely, we have the following exact triangle 
$$ Q_\emptyset \to \hourglass_N \to \wb{Q_{[N]}}[-N+1] \xto{+1} $$ 
where the connection morphism corresponds to the extension $\Ext^{N}(\wb{ Q_{[N]}}, Q_\emptyset) \cong \C$.


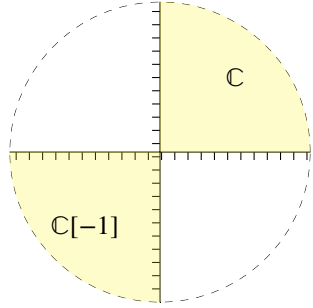
\begin{figure}[h]
    \centering
\begin{tikzpicture}
\draw [righthairs] (-2,0) -- (2,0) (0,2) -- (0,-2);
\clip (0,0) circle (2);
\draw [dashed] (0,0) circle (2);
\filldraw [yellow, opacity=0.2] (0,0) rectangle (2,2) (0,0) rectangle (-2,-2);
\node at (1,1) {$\C$};
\node at (-1,-1) {$\C[-1]$};
\end{tikzpicture}
    \caption{An hour-glass sheaf supported in the first and third quadrant, namely $Q_\emptyset$ and $Q_{12}$. }
    \label{fig:hour-glass}
\end{figure}

\begin{figure}[h]
    \centering
    \begin{subfigure}[b]{0.3\textwidth}
            \centering
            \begin{tikzpicture}
\draw [righthairs] (-2,0) -- (2,0) (0,2) -- (0,-2);
\draw [fill, opacity=0.3] (0,0) -- (-0.2, 0) arc (180:270:0.2); 
\clip (0,0) circle (2);
\draw [dashed] (0,0) circle (2);
\begin{scope}[scale=0.6,shift={(0.3, 0.3)}]
\draw [lefthairs] (1,-1) -- (-1,1);
\draw [lefthairs] (1,-1) arc (-45:-225:1.41);
\filldraw [blue, opacity=0.2] (1,-1) arc (-45:-225:1.41) -- cycle;  
\end{scope}
\end{tikzpicture}
    \end{subfigure}
\begin{subfigure}[b]{0.3\textwidth}
            \centering
            \begin{tikzpicture}
\draw [righthairs] (-2,0) -- (2,0) (0,2) -- (0,-2);
\draw [fill, opacity=0.3] (0,0) -- (-0.2, 0) arc (180:270:0.2); 
\clip (0,0) circle (2);
\draw [dashed] (0,0) circle (2);
\begin{scope}[scale=1,shift={(0.1, 0.1)}]
\draw [lefthairs, rounded corners] (0.2,-1) -- (0.2,0.2) -- (-1,0.2);
\draw [lefthairs] (0.2,-1) arc (-90:-180:1.2);
\fill [blue, opacity=0.2] (0.2,-1) -- (0.2,0.2) -- (-1,0.2) -- cycle;
\fill [blue, opacity=0.2] (0.2,-1) arc (-90:-180:1.2) -- cycle;  
\end{scope}
\end{tikzpicture}
    \end{subfigure}
    \begin{subfigure}[b]{0.3\textwidth}
            \centering
            \begin{tikzpicture}
\draw [righthairs] (-2,0) -- (2,0) (0,2) -- (0,-2);
\draw [fill, opacity=0.3] (0,0) -- (-0.2, 0) arc (180:270:0.2); 
\clip (0,0) circle (2);
\draw [dashed] (0,0) circle (2);
\begin{scope}[shift={(0.1, 0.1)}, scale=2]
\draw [lefthairs, rounded corners] (0,-1) -- (0,0) -- (-1,0);
\draw [lefthairs] (0,-1) arc (-90:-180:1);
\fill [blue, opacity=0.2] (0,-1) -- (0,0) -- (-1,0) -- cycle;
\fill [blue, opacity=0.2] (0,-1) arc (-90:-180:1) -- cycle;  
\end{scope}
\end{tikzpicture}
    \end{subfigure}
    \caption{The microlocal skyscraper sheaf $\wb{Q_{12}}$ for the partial cotangent fiber in $T_0^* \R^2$ for $\La_{[N]}$, $N=2$.}
    \label{fig:probe}
\end{figure}
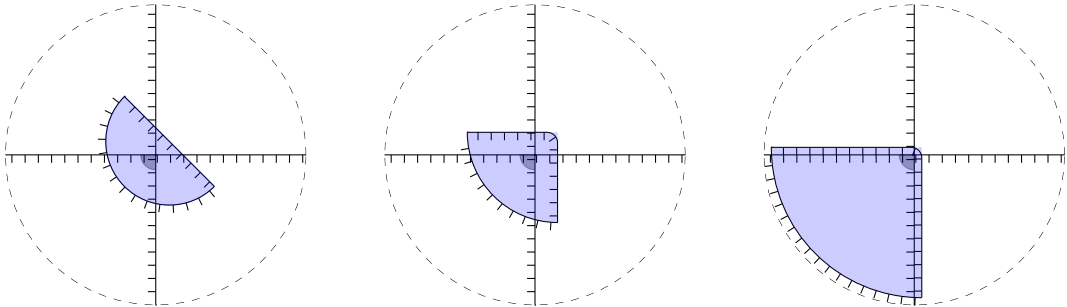

The hour-glass sheaf $\hourglass$ can also be obtained from $P_{\emptyset}$ by stop removal. For the full skeleton $\La_{[N]}$, the microlocal skyscraper for the smooth strata $\{0\} \times (- \sigma_{[N]})$ is $\wb{ Q_{[N]}}$ (up to a degree shift ambiguity). It can be proven either by a non-characteristic deformation, as shown in Figure \ref{fig:probe}, or by a computation using generators.

The probe sheaves can genenerate all other microlocal skyscrapers. Let $\scal_s \in \scal_N$ be a strata labelled by the sign vector $s$, which we abuse notation and also call $s$. Then the relative interior $\La_N|_s:=\La_N|_{\scal_s}$ is a smooth Lagrangian in the conormal to $\scal_s$. For each $s=+,-,0$, define $I_s = \{ i \in [N] : s_i = s\}$. We define the following complex of probe sheaves
\be \label{e:microstalk} F_{s, \pcal_n} \cong  \uwave{ P_{I_+ \sqcup I_-} } \to  \bigoplus_{J \in I_0, |J|=1} P_{I_+ \sqcup I_- \sqcup J } \to \bigoplus_{J \in I_0, |J|=2} P_{I_+ \sqcup I_- \sqcup J }\to \cdots \to P_{[N]}. \ee

\bp 
Let $G \in Sh^\dm(\R^N, \La_N)$, $S_s$ a strata of $\scal_N$. Then \[  \Int (\La_N|_s) \In SS(G) \quad \LRA \quad \Hom(F_{s, \pcal_n}, G) \neq 0. \] 
\ep 
\bpf 
If we construct a transverse disk to the strata $S_s$, then the only quadrants involved are those $Q_J$ with $I_- \In J \In I_- \sqcup I_0$. Hence if we represent the microlocal stalk functors using the stalk functors, and then co-represent the stalk functors in $Q_J$ as $P_J$, we get the desired expression.

\epf 
One can check that this microlocal stalk is supported in the quadrant $Q_{I_- \sqcup I_0}$ with the appropriate boundary faces.

\subsection{Leaks and flooded quadrants}
Let $\La_\ical$ be a rectilinear skeleton with zero section, i.e. $[N] \in \ical$. In this subsection, we try to understand $F_{I, \ical}$ as obtained from $P_I$ by 'leaking'. We will drop the subscript $\ical$ if the context is clear. 

Consider the closed embedding of skeleton $\La_\ical \into \La_N$, we have the corresponding fully faithful embedding of sheaves 
$$ \iota_\ical: \ccal_\ical \to \ccal $$ 
and its left adjoint is the 'stop removal' functor
$$ \iota_\ical^L: \ccal \to \ccal_\ical.$$
We have the co-unit $\eta: 1 \to \iota_\ical \circ \iota_\ical^L$. Then we have the probe sheaf
$$ F_I =  \iota_\ical \circ \iota_\ical^L(P_I) $$ 
and using the co-unit map, we define $G_I$ as the cone of $P_I \to F_I$
$$ P_I \to F_I \to G_I \xto{+1}. $$
By Corollary \ref{c:skeleton-with-floor}, the stalk of $F_I$ in $P_I$ is always $\C$, hence the support of $G_I$ is disjoint from the open set $P_I$. We say a quadrant $Q_J$ is {\bf flooded} by the probe $F_I$, if $J \not \In I$ but $F_I|_{Q_J} \neq 0$, or equivalently, $Q_J \In \supp(G_I)$.


\bd 
A {\bf leak} in $Q_I$ for the skeleton $\La_\ical$ is a face of $Q_I$ the form $L = -\tau_I + \tau_J$ for some $J \In I^c$, such that $\La_{\ical}|_{\Int(L)} \neq \La_N|_{\Int(L)}$. 
\ed 

\bl \label{l:leak2}
Let $Q_I$ be a quadrant, $J \In I^c$, the face $L = -\tau_I + \tau_J$ is a leak of $Q_I$ for the skeleton $\La_\ical$ if and only if for any $J' \In J$, $I \sqcup J' \notin \ical$. 
\el 
\bpf
The only possible contribution to the skeleton over the face $L$ are from $\La_{I \sqcup J'}$ for $J' \In J$, hence it is sufficient and necessary for these $I \sqcup J'$ to be absent in $\ical$. 
\epf

Using leaks, we can relate the probe sheaf $F_I$ to those $F_{I'}$ where $Q_{I'}$ is adjacent to $Q_I$ at the leaky edge. 
\bl \label{l:probe-from-leak}
If $L = -\tau_I + \tau_{I'}$ is a leak for $Q_I$ in $\La_\ical$, then 
$$ F_{I}  \to  \bigoplus_{J \In (I\cup I')^c, |J|=1} F_{I \sqcup J} \to \bigoplus_{J \In (I\cup I')^c, |J|=2} F_{I \sqcup J} \to \cdots \to F_{(I')^c} $$ 
\el 
\bpf 
If $L$ is a leak of codimension $k$, then by examining the $k$-dimensional transverse disk to $L$, we see the stalks on the $2^k$  quadrants on the disk are related by an acyclic Koszul complex, hence one can express $F_I$ using $F_{I'}$ for $Q_I'$ that are adjacent to $Q_I$ along the leak $L$. 
\epf



We may relate the leak and the flooded regions geometrically. The following proposition describe when the flooded region touches the source quadrant $Q_I$. 
\bp \label{p:leak-flood}
Let $\La_\ical$ be a rectilinear skeleton with zero-section. Let $I \notin \ical$, and define the following closed subset in $\R^N$, 
\begin{enumerate}
    \item $C_1 = \cup \{ -\tau_I + \tau_J \mid\z{$J \In  I^c$ and $-\tau_I + \tau_J$ is a leak for $Q_I$} \}$. 
    \item $C_2 = \supp(G_I) \cap \supp(P_I) \cap \wb{Q_I} = \supp(G_I) \cap  \wb{Q_I}$. 
    \item $C_3 = \supp (\uhom(G_I, P_I) ) \cap \wb{Q_I}$. 
    \item $C_4 = \supp(c) \cap \wb{Q_I}$, where $c: G_I \to P_I[1]$ is the connection homomorphism, and is viewed as a section of $\uhom(G_I, P_I)$. 
\end{enumerate}
Then we have 
$$ C_1 = C_4 \In C_3 = C_2. $$
\ep 
\bpf
Since all three sets can be written as a product form $C'_i \times (-\tau_I) \In \R^{I^c} \times \R^I$, hence we may quotient out the $\R^I$ factor. This is equivalent of replacing $[N]$ by $I^c$ and $\ical$ by $\ical_{\geq I}$ and  $I$ by $\emptyset$. 

Hence, we only consider the case $I = \emptyset$ from now on.  

First, we claim that $ \supp(G_I) \cap \supp(P_I) =  \supp (\uhom(G_I, P_I) )$, hence $C_2 = C_3$. It is clear that $\supp(G_I) \cap \supp(P_I) \supset \supp (\uhom(G_I, P_I) )$. Conversely, using $I = \emptyset$,  if $\tau_J$ is a maximal cone in $\supp(G_\emptyset) \cap \supp(P_\emptyset)$. For any $x_J \in \Int(\tau_J)$, consider $B(x_J)$ small open ball around $x_J$ such that $B(x_J) \cap \tau_J) \In \Int(\tau)J)$, then we see $Hom(G_\emptyset|_{B(x_J)}, P_\emptyset|_{B(x_J)}) \neq 0$, hence $\Int(\tau_J) \In \supp (\uhom(G_I, P_I) )$. By considering all such maximal $\tau_J$, we have $\supp(G_I) \cap \supp(P_I) = \wb {\cup_J(\Int(\tau_J)) }\In  \supp (\uhom(G_I, P_I) )$. Hence proving $C_2=C_3$. 

Next, we note $C_4 \In C_3$ automatically since $c$ is a section of $\uhom(G_I, P_I)$. 

Finally, suffice to show that $C_4 = C_1$. Let $\lcal = \{ J \In [N] \mid \z{ $\tau_J$ is a leak for $Q_\emptyset$} \}$, and let $J \in \lcal$ be a maximal element. Note that since $[N] \in \ical$, hence $J \neq [N]$, and $\tau_J$ is a proper face of $Q_\emptyset$. For any $x_J \in \Int(\tau_J)$ let $B(x_J)$ be a small enough open ball around $x_J$ as above, then the probe sheaf $F_\emptyset$ restricted to $B(x_J)$ has the form of an hourglass sheaf in the transverse direction of $\tau_J$ and constant sheaf along $\tau_J$, i.e. $F_\emptyset$ locally in $B(x_J)$ is an extension of $G_\emptyset$ and $P_\emptyset$, hence $\Int(\tau_J) \In C_4$. And considering all such maximal $J$, we have $C_1 = \wb{\cup_J \Int(\tau_J) } \In C_4$. Conversely, since $F_\emptyset$ is obtained from $P_\emptyset$ by stop removal along the leaks of $Q_\emptyset$, hence the connecting morphism $c$ has support in the leaks, i.e $C_4 \In C_1$. This proves $C_1 = C_4$. 
\epf 

\brem 
It is possible that the above inclusion $C_4 \In C_1$ is strict. Consider the example of $N=3$, $\ical = \{1, 2, 123\}$ (where we use an obvious abuse of notation for subsets of $[3]=\{1,2,3\}$). Then, probe sheaf for $Q_\emptyset$ is 
$$ F_\emptyset = P_1 \oplus P_2 \to P_{123} $$
and $G_\emptyset$ is supported in $Q_{12}$, $Q_{3}, Q_{13}, Q_{23}, Q_{123}$. We note that the quadrant $Q_{3}$ is adjacent to $Q_\emptyset$ through the maximal faces $\tau_{3}, \tau_{12}$, i.e
$$ C_2=C_3 = \tau_3 \cup \tau_{12} $$
However the leak is only at 
$$ C_1 = C_4 = \tau_3. $$
\erem

\subsection{Singular support of the constructible sheaf of categories} 
\label{ss:SS-cat}
Let $\La_\ical$ be a rectilinear skeleton in $\R^N$. For simplicity, we assume $\La_\ical$ contains the zero-section of $T^*\R^N$. For a open convex set $U \In \R^N$, we consider the category $Sh_{\La_\ical}^\dm(U)$. These data assemble into a constructible sheaf of categories (or Kashiwara-Schapira stack) denoted as $Sh_{\La_\ical}^\dm$, and we want to describe $SS(Sh_{\La_\ical}^\dm)$. 

We define three versions of singular supports. For any non-zero covector $(x, \xi) \in T^*\R^N$, we let $B_x$ be a small enough open ball around $x$, and $B_{x,\xi,-} = \{ x' \in B_x \mid \la x' - x, \xi \ra < 0 \}$, then we have the 'microlocal restriction functor' 
$$ \rho_{x,\xi}: Sh_{\La_\ical}^\dm(B_x) \to Sh_{\La_\ical}^\dm(B_{x,\xi,-}) $$
and the 'microlocal co-restriction functor'
$$ \rho_{x,\xi}^L: Sh_{\La_\ical}^\dm(B_{x,\xi,-}) \to Sh_{\La_\ical}^\dm(B_{x}). $$

Let $SS(Sh_{\La_\ical}^\dm)$ denote the complement of the set of covectors $(x,\xi)$ where there is an open neighborhood $U$ of $(x,\xi)$ such that $\rho_{x', \xi'}$ is an equivalence for any $(x', \xi') \in U$. Similarly, let $SS_{Hom}(Sh_{\La_\ical}^\dm)$ (resp. $SS_{Hom}^L(Sh_{\La_\ical}^\dm)$) be defined by replacing '$\rho_{x', \xi'}$ is an equivalence' with  '$\rho_{x', \xi'}$ (resp. $\rho^L_{x', \xi'}$) is fully-faithful'. Since the condition that $\rho_{x,\xi}$ is an equivalence is equivalent to $\rho_{x,\xi}$ and  $\rho_{x,\xi}^L$ are fully-faithful, we have 
$$ SS(Sh_{\La_\ical}^\dm) = SS_{Hom}(Sh_{\La_\ical}^\dm) \cup SS^L_{Hom}(Sh_{\La_\ical}^\dm). $$
For simplicity of notation, we will use $SS(\La_\ical)$ to denote $SS(Sh_{\La_\ical}^\dm)$ and similar for its variants.

\bp 
$$ SS_{Hom}(\La_\ical) =  \bigcup_{I, J \in \ical} SS(\uhom(P_I, P_J)) \In T^* \R^N.$$
\ep 
We denote the right hand side of the above equation as $\uhom(\La_\ical)$.
\bpf 
For any $x \in \R^N$, the category $Sh_{\La_\ical}^\dm(B_x)$ is generated by the restriction of $\{ P_I: I \in \ical\}$ to $B_x$, since the restriction of $\La_{\ical}$ to $B_x$ is another rectilinear skeleton centered at $x$ indexed by a subset of $\ical$. 
\epf

\bl \label{l:uhom-IJ}
For any $I, J \In [N]$, we have
$$ \uhom(P_I, P_J) = \R^{I \cap J} \times (\R_{>0})^{I \RM J} \times (\R_{\geq 0})^{I^c},$$ 
where we abuse notation and use a locally closed subset $A$ with the locally constant sheaf $\C_A$ supported on $A$. 
\el  
\bpf 
For any $I \In [N]$, we have 
$$ P_I = P_{I,1} \boxtimes P_{I,2} \boxtimes \cdots \boxtimes P_{I,N}, \quad \z{ where } P_{I,i} = \bcs \R_{>0} & i \notin I \\ \R & i \in I \ecs. $$
The hom sheaf in each factor can be computed using
\begin{eqnarray*}
\uhom(\R_{>0}, \R_{>0}) &\cong& \R_{\geq 0}, \\  
\uhom(\R_{>0}, \R) &\cong& \R_{\geq 0}, \\ 
\uhom(\R, \R_{>0}) &\cong& \R_{>0}, \\
\uhom(\R, \R) &\cong& \R
\end{eqnarray*}
\epf 

\bl \label{l:ss-uhom-IJ}
For any $I, J \In [N]$, we have
$$ SS(\uhom(P_I, P_J)) \cap T_0^* \R^N = \sigma_{s(I, J)} $$
where  $s(I, J)$ is a sign vector defined by
\[ s(I, J)_i = \bcs +1 & i \notin I \\
0 & i \in I, i \in J \\
-1 & i \in I, i \notin J
\ecs. \]
Moreover, for any $x \in \R^N$, we have
$$ SS(\uhom(P_I, P_J)) \cap T_x^* \R^N  \In \sigma_{s(I, J)}. $$
\el 
\bpf
It suffices to check we have the desired claim for each $i \in [N]$, which is a straight-forward calculation using Lemma \ref{l:uhom-IJ}. 
\epf 

\brem 
In \cite{zhou2020}, we give a different approach to compute the hom sheaves singular support $\uhom(\La_\ical)$ as 
$$ \uhom(\La_\ical) = \lim_{\epsilon \to 0^+ } -\La_\ical + G^\epsilon \La_\ical $$
where $G^\epsilon$ is translation by $-\epsilon \cdot (1,..., 1)$. For general skeleton $\La$, $G^\epsilon$ is some positive isotopy, and the right-hand side is only an upper bound for the singular support of the hom sheaf. 
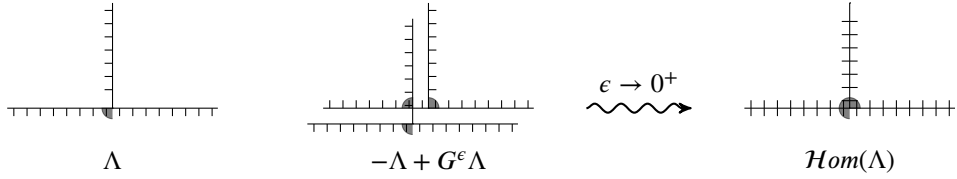
\begin{center}
\begin{tikzpicture} [scale=0.7]
\begin{scope} [shift={(-6,0)}]
\draw [righthairs] (-2,0) -- (2,0) (0,2) -- (0,0);
\draw [fill, opacity=0.5] (0,0) -- +(180:0.2) arc (180:270:0.2);
\node at (0, -1) {$\La$}; 
\end{scope}

\begin{scope}
\def\sf{0.3}
\draw [righthairs, shift={(-\sf, -\sf)}]  (-2,0) -- (2,0) (0,2) -- (0,0);
\draw [fill, opacity=0.5, shift={(-\sf, -\sf)}] (0,0) -- +(180:0.2) arc (180:270:0.2);

\draw [lefthairs] (-2,0) -- (2,0) (0,2) -- (0,0);
\draw [fill, opacity=0.5] (0,0) -- +(0:0.2) arc (0:90:0.2);

\draw [fill, opacity=0.5, shift={(-\sf, 0)}] (0,0) -- +(90:0.2) arc (90:180:0.2);
\node at (0, -1) {$-\La + G^\epsilon \La $}; 
\end{scope}

\draw [arrowstyle, -stealth'] (3,0) -- (5,0); 
\node at (4,0.5) {$\epsilon \to 0^+$};

\begin{scope}[shift={(8,0)}]
\draw [righthairs, lefthairs] (-2,0) -- (2,0) (0,0) -- (0,2);
\draw [fill, opacity=0.5] (0,0) -- +(0:0.2) arc (0:270:0.2);
\node at (0, -1) {$\uhom(\La)$}; 
\end{scope}

\end{tikzpicture}
\end{center}

\erem

\subsection{Rectilinear skeleton as family of skeletons}
We first recall how singular support behaves under push forward. Let $f: M \to N$ be a smooth submersion, and $\La \In T^*M$ be a conical Lagrangian. We define 
$$ f_* \La = \{ (x, \xi) \in T^*N \mid \z{ there exists $(\wt x, \wt \xi) \in \La$, such that $f(\wt x) = x, f^*(\xi) = \wt \xi.$} \} $$
We say $\La$ is {\bf $f$-non-characteristic} if $f_* \La$ is contained in the zero section of $T^*N$. 
In general, if $F$ is a constructible sheaf on $M$, then $SS(\pi_* F) \In \pi_*(SS F)$ due to possible cancellations in pushing forward a sheaf. However, if $F \in Sh(\R^N)$ is a conic constructible sheaf, and $\pi: \R^N \to \R^k$ is a linear surjective map, then $SS(\pi_* F) = \pi_* (SS F)$. If $b \in N$ and $\La$ is $f$-non-characteristic, then we define $ M_b = f^{-1}(b)$, and
$$\La|_b = q_b (\La \cap T^* M|_{M_b} ) \In T^* M_b $$
where $q_b$ is the quotient map in 
$$ 0 \to T_{M_b}^* M \to T^* M|_{M_b} \xto{q_b} T^* M_b \to 0. $$

Let $\La_\ical$ be a rectilinear skeleton, and $\pi: \R^N \to \R^k$ be a linear map. We always assume that $\La_N$ (hence the sub-skeleton $\La_\ical$) is $\pi$-non-characteristic.

\bl $\La_N$ is $\pi$-non-characteristic if and only if $\ker(\pi)$ intersects the open positive quadrant $\R_{>0}^N$. 
\el 
\bpf 
If $\emptyset \neq \R_{>0}^N \cap \ker(\pi)$ and $x=(x_1, \cdots, x_N)$ is in the intersection, then $x_i > 0$. For any $\xi=(\xi_1, \cdots, \xi_N)$ in some cone of $\Sigma_N$, we have $\xi_i \leq 0$. Thus $ \la \xi, x \ra = \sum_i \xi_i x_i < 0 $, hence $\xi$ is not conormal to $\ker(\pi)$. Thus, if   $\emptyset \neq \R_{>0}^N \cap \ker(\pi)$,  we have $\La_N$ is $\pi$-non-characteristic.

If $\emptyset = \R_{>0}^N \cap \ker(\pi)$, $\gamma = \pi(\R_{>0}^N)$ is an open convex cone in $\R^k$, hence we may lift an element $\bar \xi$ the exterior conormal to $\gamma$ at  $0\in \R^k$, to $\xi \in T_0^*\R^N$, which would be in the singular support of $\C_{\R_{>0}^N}$ at $0$, hence in the support of the fan $\Sigma_N$. Thus, $\La_N$ is not $\pi$-non-characteristic. 
\epf

We consider the sheaf of categories $\pi_* Sh_{\La_\ical}^\dm$ on $\R^k$, where for any $U$ convex open set in $\R^k$, we have 
$$ \pi_* Sh_{\La_\ical}^\dm(U) = Sh^\dm_{\La_\ical}(\pi^{-1}(U)). 
$$
We define the singular support $ SS^\bu( \pi_* Sh_{\La_\ical}^\dm ) $ as in section \ref{ss:SS-cat}, where $SS^\bu = SS, \, SS_{Hom}, \, SS_{Hom}^L$. 

We have the following relations for $SS$ and its variants,
$$ SS(\pi_* Sh_{\La_\ical}^\dm )) \In \pi_* SS(Sh_{\La_\ical}^\dm )),  $$

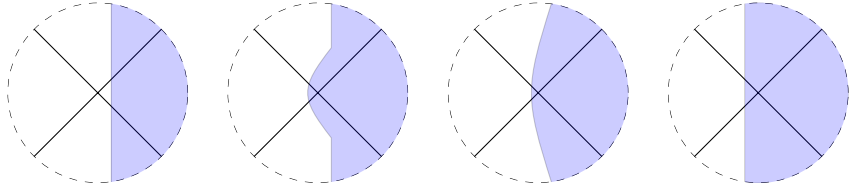
\begin{figure}[hh]
\begin{center}
    \begin{tikzpicture}[scale=0.6]

    \clip (0,0) circle (2);
    \draw [dashed] (0,0) circle (2);
    \draw (-2,2) -- (2, -2) (-2,-2)--(2,2); 
    \draw [fill=blue, opacity=0.2] (0.3,-2) rectangle (2,2); 

    \end{tikzpicture}
    $\quad$
    \begin{tikzpicture}[scale=0.6]

    \clip (0,0) circle (2);
    \draw [dashed] (0,0) circle (2);
    \draw (-2,2) -- (2, -2) (-2,-2)--(2,2); 
    \draw [fill=blue, opacity=0.2] (0.3,-2) -- (0.3,-1) .. controls (-0.4, -0.1) and (-0.4, 0.1) .. (0.3, 1) -- (0.3,2) -- (3,3) -- (3,-3) -- cycle; 

    \end{tikzpicture}
        $\quad$
    \begin{tikzpicture}[scale=0.6]

    \clip (0,0) circle (2);
    \draw [dashed] (0,0) circle (2);
    \draw (-2,2) -- (2, -2) (-2,-2)--(2,2); 
    \draw [fill=blue, opacity=0.2] (0.3,-2)  .. controls (-0.3, -0.1) and (-0.3, 0.1) ..  (0.3,2) -- (3,3) -- (3,-3) -- cycle; 

    \end{tikzpicture}  $\quad$
       \begin{tikzpicture}[scale=0.6]

    \clip (0,0) circle (2);
    \draw [dashed] (0,0) circle (2);
    \draw (-2,2) -- (2, -2) (-2,-2)--(2,2); 
    \draw [fill=blue, opacity=0.2] (-0.3,-2) rectangle (2,2); 
    \end{tikzpicture}
\end{center} 
    \caption{Global propagation of sections from local extendability. }
    \label{fig:local-to-global-extension}
\end{figure}

\bp \label{p:coFF}
Let $\La_\ical$ be a rectilinear skeleton containing zero section.  Then the following conditions are equivalent 
\begin{enumerate}
    \item $SS^L_{Hom}(\pi_* Sh_{\La_\ical}^\dm )  = T^*_{\R^k} \R^k$, 
    \item For any $U \In \R^k$ a convex open set, the co-restriction functor 
$$ \rho_U^L: \pi_* Sh^\dm_{\La_\ical}(U) \to \pi_* Sh^\dm_{\La_\ical}(\R^k) $$
is fully-faithful. 
    \item For any $I \notin \ical$,we have 
    $$  \wb{ \pi(Q_I)} \cap \pi( \z{Leak}(I) )= \wb{ \pi(Q_I)} \cap  \pi(\z{Flood}(I) ), $$ 
    where $\z{Leak}(I)$ is the union of leaks of $Q_I$, and $\z{Flood}(I)$ is the closure of union of the flooded region. 
\end{enumerate}
\ep 
\bpf 
If $A \In \R^k$, we write $\wt A = \pi^{-1}(A)$ for short. 

(1) $\RA$ (2): 
Let $\ccal(U) = Sh^\dm_{\La_\ical}(\wt U)$. 
Without loss of generality, we may assume $U$ is an open convex polytope. Then, we may expand $U$ as increasing family of convex open polytope $U_t$, $t \in [0, \infty)$, such that there are finitely many times $0 < t_1 < t_2 < \cdots < t_k$ when the boundary of the preimage $\pa \wt U_t$ intersects with some strata in $\scal_N$ of $\R^N$ non-transversely. See Figure \ref{fig:local-to-global-extension} for an illustration of how to resolve the critical moment of 'strata-passing' into smaller steps.  By the vanishing of the singular support of microlocal co-restriction functor, we see the co-restriction $\ccal(U_{t_i - \epsilon}) \to \ccal(U_{t_i + \epsilon})$ is fully faithful.  And after the last transition, we also have $\ccal(U_{t_k + \epsilon}) \congto \ccal(\R^k)$. Since $\rho_U^L$ is the composition of fully faithful functors, $\rho_U^L$ is fully-faithful. 

(2) $\RA$ (1): For any non-zero $(x, \xi) \in T^* \R^k$, the microlocal co-restriction 
$$ \rho_{x,\xi}^L: \pi_* Sh^\dm_{\La_\ical}( B_{x,\xi,-}) \to  \pi_* Sh^\dm_{\La_\ical}(  B_{x}) $$
posted composed with a fully-faithful functor $\rho_{B_x}^L$ is another fully-faithful functor $\rho_{B_{x,\xi,-}}^L$, hence $\rho_{x,\xi}^L$ is fully-faithful.

(2) $\RA$ (3): Suppose (3) fails. Since $\z{Leak}(I) \subset \z{Flood}(I)$ (c.f. Proposition \ref{p:leak-flood}), then there exists an $I \notin \ical$, and a point $x \in \wb{ \pi(Q_I)} \cap \pi (\z{Flood}(I))$ and $x \notin \wb{ \pi(Q_I)} \cap  \pi(\z{Leak}(I))$. Let $B_x$ be a small ball around $x$, such that $B_x \cap  ( \wb{ \pi(Q_I)} \cap  \pi(\z{Leak}(I))) = \emptyset$. Then we claim the co-restriction functor $\rho_{B_x}^L$ from $\wt B_x = \pi^{-1}(B_x)$ to $\R^N$ is not fully-faithful, contradicting with (2). Indeed, consider the probe sheaf $F_{I|B_x}$ for $Q_I \cap \wt B_x$. Since there is no leak of $Q_I$ over $B_x$, we have $F_{I|B_x} = P_I \cap \wt B_x$. On the other hand, the co-restriction $\rho_{B_x}^L$ sends $F_{I|B_x}$  to $F_I$. Then consider certain $Q_J \In \z{Flood(I)}$, such that $x \in \wb{ \pi(Q_I)} \cap \wb{ \pi(Q_J)}$, $\wt B_x \cap Q_J \neq \emptyset$. Hence
$$ \Hom(F_{J|B_x}, F_{I|B_x}) = 0, \quad \Hom(\rho_{B_x}^L F_{J|B_x}, \rho_{B_x}^L  F_{I|B_x}) = \Hom(F_J, F_I) \neq 0 $$
thus $\rho_{B_x}^L$ is not fully-faithful, which is a contradiction to the assumption. Hence (3) holds.

(3) $\RA$ (2): We will show that "probe commutes with restriction", i.e. the probe sheaf of $Q_I$ in the restricted region $\wt U:=\pi^{-1}(U)$ is equal to the restriction of the probe $F_I$ to $\wt U$, 
\be F_{I|U} = F_I|_{\wt U} \label{e:probe-restriction-commute} \ee

We first observe that if $I \in \ical$, then the restricted probe sheaf is $F_{I|U} = P_I|_{\wt U}$ since the sheaf $Q_I|_{\wt U}$ can expand non-characteristically to $P_I|_{\wt U}$ and $SS(P_I|_{\wt U}) \In \La_{\ical}|_{\wt U}$. Hence the claim is proven in this case. 

Now we assume $I \notin \ical$. We assume by induction that the cases for $I'$ with $|I'| > |I|$ are proven. We consider two cases. 

{\bf case a:} If there are leaks of $Q_I$ over $U$, i.e., there exists a leak $-\tau_I + \tau_{I'}$ of $Q_I$ and $-\tau_I + \tau_{I'} \cap \wt U \neq \emptyset$, then by Lemma \ref{l:probe-from-leak},  we have
$$ F_{I}  \congto  \uwave{\bigoplus_{J \In (I\cup I')^c, |J|=1} F_{I \sqcup J} } \to \bigoplus_{J \In (I\cup I')^c, |J|=2} F_{I \sqcup J} \to \cdots \to F_{(I')^c} $$ 
and 
$$ F_{I|U}  \congto  \uwave{\bigoplus_{J \In (I\cup I')^c, |J|=1} F_{I \sqcup J} |U} \to \bigoplus_{J \In (I\cup I')^c, |J|=2} F_{I \sqcup J|U} \to \cdots \to F_{(I')^c|U}. $$
Since by induction hypothesis, we have $F_{I \sqcup J |U} \cong F_{I \sqcup J}|_{\wt U}$ for any non-empty $J$, we can restrict the resolution of $F_I$ to $\wt U$ and get the same resolution as  $F_{I|U}$, hence $F_{I|U} = F_I|_{\wt U}$. 

{\bf case b:} If there are no leaks of $Q_I$ over $U$, then by condition (3), there is no flooded quadrant over $U$, i.e., no $Q_J$ such that $Q_J \notin P_I$ and $F_I|_{Q_J} \neq 0$ and $Q_J \cap \wt U \neq \emptyset$. Hence $F_I|_{\wt U} = P_I|_{\wt U}$, which also shows $SS(P_I|_{\wt U}) \In \La_{W_K,v}|_{\wt U}$, thus 
$$ F_{I|U} = P_I|_{\wt U} =  F_I|_{\wt U}. $$
This finishes the proof of the claim. 
\epf 

Since we will only consider skeletons with the above properties, we give it a name. 
\bd 
If $\La_\ical$ is a rectilinear skeleton with zero section, and $\pi: \R^N \to \R^k$ is non-characteristic with respect to $\La_\ical$. We say $\La_\ical$ is {\bf $\pi$-coff} if the microlocal co-restrcition functors are always fully-faithful, i.e. $SS^L_{Hom}(\pi_* Sh_{\La_\ical}^\dm )  = T^*_{\R^k} \R^k$. 
\ed 

Finally, we prove that given a $\pi$-coff skeleton $\La_\ical$, one can extend a constructible sheaf on a fiber on $\pi$ with respect to the restriction of the skeleton to a tubular neighborhood. 

\bp \label{p:coFF-fiber}
Let $\La_\ical$ be a $\pi$-coff skeleton. Then for any $x \in \R^k$, the
co-restriction functor (left-adjoint to restriction)
$$ \rho_x^L:  Sh^\dm(X_b, \La_\ical|_b) \to Sh^\dm(\R^N, \La_\ical) $$
is fully-faithful, where $X_b = \pi^{-1}(b)$, $\La_\ical|_b = \La_\ical|_{X_b}$. 
\ep 
\bpf 
Suffice to check that, for each $I \In [N]$, such that $Q_I \cap X_b \neq \emptyset$, we have the probe sheaf $F_{I|b}$ of region $Q_I \cap X_b$ for skeleton $\La_{\ical}|_b$ equal to the restriction of the global probe sheaf $F_I|_{X_b}$. We proceed as the "(3) $\RA$ (2)" step in the proof of Proposition \ref{p:coFF}, with some modification. 

If $I \in \ical$, then there is nothing to prove, as $F_{I|b} = P_I|_{X_b}=F_I|_{X_b}$. If $I \notin \ical$, then we may prove by induction and assume the case for $I'$ with $|I'| > |I|$ is verified. Since $\pi(Q_I) \ni b$, then $b \in \pi(\z{Leak}(I))$ if and only if $b \in \pi(\z{Flood}(I))$. We have two cases: 

{\bf case a:} $b \in \pi(\z{Leak}(I))$. However unlike the case with open set, this does not imply that $X_b$ contains a transverse disk to the leak of $Q_I$. Suppose there exists a leak $-\tau_I + \tau_J$, such that $X_b$ intersects $\Int(-\tau_I + \tau_J)$ transversely, then the induction step follows as before. Now we consider the case that, for all leak $L = -\tau_I + \tau_J$ of $Q_I$, $\pi(\Int(L))$ does not cover an open neighborhood of $b$. Hence $b \in \pa(\pi(\z{Leak}(I))) \cap \pi(Q_I) = \pa(\pi(\z{Flood}(I))) \cap \pi(Q_I)$, which also means there is all the flooded region $Q_J$ of $Q_I$ satisfies $Q_J \cap X_b = \emptyset$. This is then the same as case (b) below. 

{\bf case b:} If $b \not \in \pi(\z{Leak}(I))$, then $X_b$ does not see any leak and also does not see any flooded quadrant, hence we again have $F_{I|b} = P_I|_{X_b}=F_I|_{X_b}$. 
\epf 

\subsection{Examples for failure of extension of constructible sheaves}
\begin{example}
Let $N=2$, and consider $\ical=\{1,2\}$. We can draw the picture of $\La_\ical$ as following Figure \ref{fig:fail-2d}. Let $U_1 = \{ (x,y) \mid x-y > 1\}$ and $U_2 = \{ (x,y) \mid x-y > -1\}$. It is impossible to extend the sheaf over $U_1$ supported in the blue region to $U_2$. If we apply the co-restriction functor $\rho^L$ along the inclusion $U_1 \into U_2$, then we get the 2d hour-glass sheaf restricted to $U_2$, which changes the stalk at $b$. In other words, the hom space morphism 
$$ 0 \cong Hom_{U_1}(F_b, F_a) \to Hom_{U_2}(\rho^L F_b, \rho^L F_a) \cong \C[-1] $$
is not a quasi-isomorphism. 

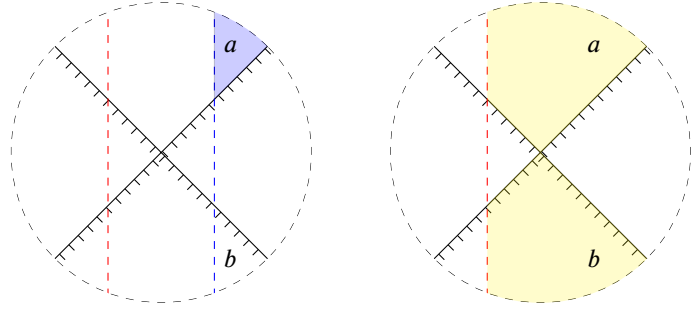
\begin{figure}[h]
    \centering
    \begin{subfigure}[b]{0.3\textwidth}
\begin{tikzpicture}[rotate=45]
\clip (0,0) circle (2);
\draw [dashed] (0,0) circle (2);
\draw [righthairs] (-2,0) -- (2,0) (0,2) -- (0, -2); 
\draw [dashed, blue, shift={(0.5, -0.5)}] (2,2 ) -- (-2,-2); 
\draw [dashed, red, shift={(-0.5, +0.5)}] (2,2 ) -- (-2,-2); 
\draw [fill=blue, opacity=0.2] (1,0) -- (3,0) arc (0:45:2) -- cycle;
\node[anchor=west] at (1.5, 0.5) {$a$};
\node[anchor=west] at (-0.5, -1.5) {$b$};

\end{tikzpicture}              
    \end{subfigure}
    \begin{subfigure}[b]{0.3\textwidth}
\begin{tikzpicture}[rotate=45]
\clip (0,0) circle (2);
\draw [dashed] (0,0) circle (2);
\draw [righthairs] (-2,0) -- (2,0) (0,2) -- (0, -2); 
\draw [dashed, red, shift={(-0.5, +0.5)}] (2,2 ) -- (-2,-2); 
\begin{scope}
\clip (1.5,2.5 ) -- (-2.5,-1.5) -- (0,-5) -- (5,0) -- cycle;
\filldraw [yellow, opacity=0.2] (0,0) rectangle (2,2) (0,0) rectangle (-2,-2);
\end{scope}
\node at (1.5, 0.5) {$a$};
\node at (-0.5, -1.5) {$b$};
\end{tikzpicture}              
    \end{subfigure}
    \caption{Example where co-restriction fails to be fully-faithful. $U_1$ is to the right of the blue line, and $U_2$ is to the right of the red line. Here, the hom space between probe sheaves $\Hom(F_b, F_a)$ are not preserved under co-restriction. }
    \label{fig:fail-2d}
\end{figure}
\end{example}

Here is another example in 3-dimension. 

\begin{example}
Consider the following rectlinear Lagrangian skeleton in $T^* \R^3$, where $\ical = \{1,2, 123\}$. Consider map $\pi: \R^3 \to \R$, where $\pi(x_1, x_2,x_3) = x_1 + x_2 - 2x_3$. We are going to show the restriction of the skeletons over fiber $\pi^{-1}(1)$ and $\pi^{-1}(-1)$. Let $U_1 = \pi^{-1}( (1, \infty))$ and $U_2 = \pi^{-1}( (-1, \infty))$, and point $a = (1, 1, 0.1)$. Then the probe sheaf $F_a$ in $U_1$ is confined in $U_1 \cap Q_\emptyset$, where as the probe sheaf $F_a$ in $U_2$ is leaked into $Q_{12}, Q_{3}, Q_{13}, Q_{23}, Q_{123}$. In particular, $Q_3$ intersects $U_1$, hence $F_a$ flooed into a region that is previously inaccessible in $U_1$. What is shown in the figure is the fiber-wise probe sheaves for the same region (restriction of $Q_\emptyset$, the top wedge), and the sudden leaking indicates that the co-restriction is not fully-faithful. 

\begin{figure}[h]
    \centering
    \begin{subfigure}[b]{0.3\textwidth}
\begin{tikzpicture}
\clip (0,0) circle (2);
\draw [dashed] (0,0) circle (2);
\draw [righthairs] (-2,0) -- (-1,0) -- (2,3); 
\draw [lefthairs] (2,0) -- (1,0) -- (-2,3); 
\fill [blue, opacity=0.2] (-2,3)-- (1,0)  -- (-1,0) -- (2,3);
\fill [blue, opacity=0.2] (-2, 0) -- (2, 0) arc (0: -180: 2); 
\end{tikzpicture}           
\caption{$\La_\ical$ over fiber $\pi^{-1}(-1)$}
    \end{subfigure}
    \begin{subfigure}[b]{0.3\textwidth}
\begin{tikzpicture}
\clip (0,0) circle (2);
\draw [dashed] (0,0) circle (2);
\draw [righthairs] (-2,0) -- (0.5,0) -- (2.5,2); 
\draw [lefthairs] (2,0) -- (-0.5,0) -- (-2.5,2); 
\fill [blue, opacity=0.2]  (2.5,2) -- (0.5,0) -- (-0.5,0) -- (-2.5,2); 
\end{tikzpicture}           
\caption{$\La_\ical$ over fiber $\pi^{-1}(1)$}
    \end{subfigure}
    \caption{Another example in 3 dimensions, where co-restriction is not fully-faithful. The probe sheaf in a slice suffers a sudden leak as one moves from $\pi^{-1}(\epsilon)$ to $\pi^{-1}(-\epsilon)$ for small positive $\epsilon$.}
    \label{fig:fail-3d}
\end{figure}
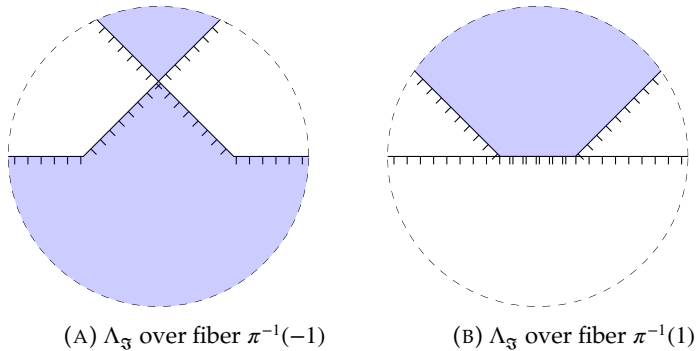
\end{example}

\section{Local Window Skeleton and the Restriction Functor \label{s:SS-Hom}} 
In this section, we will assume quasi-symmetry. 

Recall that the zonotope $\nabla$ is a symmetric polytope $\nabla = (1/2) \mu([0,1]^N)$. For a shift parameter $\delta \in \R^k$, we have the shifted zonotope $\nabla_\delta = \delta + \nabla \In \R^k$ and the shifted window 
$W_\delta = \nabla_\delta \cap \Z^k$. The function $\delta \mapsto W_\delta$ is locally constant and induces a stratification on $\R^k$. We want to describe how the skeleton $\La_{W_\delta}$ changes as $\delta$ varies.

\subsection{Window skeleton near a point}
Let $\La = \La_{W_\delta}$ be some window skeleton, or any skeleton adapted to the grid stratification of $\R^N$.

For any $x \in \R^N$, we consider the specialization $\La_x \In T^*(T_x \R^N)$ of $\La$ at $x$. The specialization $\La_x$ can be obtained in the following way: take a small ball $B_\epsilon(x)$ around $x$ so that $B_\epsilon(x)$ only intersects with strata that contains $x$ in the closure, then restrict $\La$ to $B_\epsilon(x)$, finally identify $B_\epsilon(x)$ with $B_\epsilon(0)$ in $\R^N$ and extrapolate $\La$ from the ball to $\R^N$ as a bi-conical Lagrangian. Here we identified $T_x \R^N$ with $\R^N$ in the obvious way. 

For any $x \in \R^N$, let $\nbhd(x)$ denote the union of strata that contain $x$ in the closure. One can see $\nbhd(x)$ is indeed an open neighborhood of $x$. The specialization $\La_x$ determines the restriction $\La|_{\nbhd(x)}$. 
Let $\wt v \in \Z^N$ and $v = \mu_\Z(\wt v)$. Then the specialization $\La_{W_\delta, \wt v}$ only depends on $v$, since $\La_{W_\delta}$ is invariant under translation by the sub-lattice $\ker(\mu_\Z)$. We thus define $\La_{W_\delta, v} = \La_{W_\delta, \wt v}$, and view it as a Lagrangian on $\R^N$. 

\subsection{Intermediate Level \label{ss:Rm} }
It turns out many discussions would simplify if we consider the factorization of $\mu$
$$ \R^N \xto{\pi} \R^m \xto{q} \R^k, $$
and work with the map $\pi$ first. 

Let $E_1, \cdots, E_m$ be basis of $\R^m$. Define
$\B = (1/2)\pi([0,1]^N).$
For any closed convex polytope $K$, we define
$$ \B_K = K + \B, \quad \wh W_K = \B_K \cap \Z^m, \quad \wt W_K = \pi_\Z^{-1}(\wh W_K). $$
We define the Lagrangian skeleton in $T^* \R^N$ as
$$ \La_{\wh W_K} = \bigcup_{w \in \wt W_K} \La_w. $$

\brem 
If we let $K = q^{-1}(\delta)$, then $\La_{\wh W_K} = \La_{W_\delta}$ defined before. For now, we allow $K$ to be any closed convex polytope to keep the discussion more general and independent of the map $q$. 
\erem 

Let $\wh v \in \Z^m$ and $K \In \R^m$ be a closed convex polytope. 
We first give the combinatorial data that describe the skeleton $\La_{\wh W_K, \wh v}$. Define an index set
$$ \ical(\wh W_K, \wh v) = \{ I \In [N] \mid \La_I \In \La_{\wh W_K, \wh v} \}. $$

\bl \label{l:ical-1}
$I \in \ical(\wh W_K, \wh v)$ if and only if there exists a $w \in \wt  W_K$ and $v \in \pi_\Z^{-1}(\wh v)$, such that $w_i = v_i$ for $i \notin I$ and $w_i < v_i$ for $i \in I$ 
\el 
\bpf
Suffice to consider for all $w \in \Z^N$, how $\La_w$ constribute to skeleton $\La_{w,v}$ near $v$. The contribution is non-empty if and only $w_i \leq v_i$ for all $i$. If we let $I(w,v)=\{i: w_i = v_i\}$, then $\La_{w,v} = \La_{I(w,v)}$. Hence
$$  \La_{\wh W_K, \wh v} =  \bigcup_{w \in \wt W_K} \La_{w,v} = \bigcup_{w \in \wt W_K} \La_{I(w,v)} = \bigcup_{I \in \ical(\wh W_K, \wh v)} \La_{I}. $$
\epf

Next, we relate the index set $\ical(\wh W_K, \wh v)$ to the geometry of $K$. Recall that $e_i$ is a basis of $\R^N$, $e_I = \sum_{i \in I} e_i$, $\tau_I = \cone(e_i, i \in I)$, $\gamma_I = \tau_I^o \cap \Z^N$. 


\bp \label{p:ical-Rm}
$I \in \ical(\wh W_K, \wh v)$ if and only if $K \cap (\wh v + R_{I})  \neq \emptyset$, where 
$$ R_I = \B - \pi(\tau_I) - \pi(e_I). $$
\ep 
\bpf 
Without loss of generality, we only consider the case $\wh v=0$. 

By Lemma \ref{l:ical-1}, $I \in \ical(\wh W_K, 0)$ if and only if $ - \pi(\gamma_I) \cap \wh W_K \neq \emptyset$. Since $\wh W_K = (\B + K) \cap \Z^m$, then 
$( - \pi(\gamma_I)) \cap W_K = (- \pi(\gamma_I)) \cap (\B + K).$

If $x \in (- \pi(\gamma_I))  \cap (\B + K)$, then $x = b + y$ for some $b \in \B, y \in K$. Hence $- x \in \pi(\gamma_I) \In \pi(e_I) + \pi(\tau_I)$, and then $k = x-b \in - \B - \pi(e_I) - \pi(\tau_I) = R_{I}$, where we used $-\B = \B$. Hence $I \in \ical(\wh W_K, 0)$ implies $K \cap R_{I}  \neq \emptyset$. 

If $K \cap  R_{I}  \neq \emptyset$, and let $\delta$ be in the intersection, then $$ \delta \in \B - \pi(\tau_I) - \pi(e_I) $$
hence 
$$ (\delta + \B) \cap (- \pi(\tau_I) - \pi(e_I)) \neq \emptyset. $$
Therefore, in each direction $i \in [m]$, we have 
$$ \delta_i + [-\eta_i/2, \eta_i/2] \cap (- \pi(\tau_{I_i} - \pi(e_{I_i}))) \neq \emptyset. $$
We claim that 
$$ \delta_i + [-\eta_i/2, \eta_i/2] \cap (- \pi(\gamma_{I_i})) \neq \emptyset. $$
Indeed, if $I_i$ is of mixed type, then $- \pi(\tau_{I_i}) - \pi(e_{I_i}) = \R$, and $- \pi(\gamma_{I_i})$ is a shifted sub-lattice of $\Z$ with step size $\leq \eta_i$, hence intersects with $\delta_i + [-\eta_i/2, \eta_i/2]$. Then assume $I_i$ is of pure-type. If $I_i = \emptyset$, then there is nothing to prove. If $I_i \neq \emptyset$, then $v_i - \pi(\gamma_{I_i})$ is a discrete subset in the ray $- \pi(\tau_{I_i} - \pi(e_{I_i}))$, that contains the endpoint $ - \pi(e_{I_i}))$, and has step size  $\leq \eta_i$, hence it would intersect $\delta_i + [-\eta_i/2, \eta_i/2]$. 
\epf 

\subsection{Stratification of the shift parameter space}
For any $\delta \in \R^k, v \in \Z^k$, let $\ical(\delta, v)$ denote the collection of subsets $I \In [N]$, such that $\La_I \In \La_{W_\delta, v}$. Then we have
$ \La_{W_\delta, v} = \cup_{I \in \ical(\delta, v)} \La_I. $
We will consider the stratification on the parameter space $\delta$ induced by the function $\delta \mapsto \ical(\delta, v)$. 

\bp 
$I \in \ical(\delta, v)$ if and only if $\delta \in v + A_I$, where 
\be A_I = q(R_I)=  -\beta_I - C_I + \nabla. \ee
\ep 
\bpf 
Without loss of generality, set $v=0$. Since $\La_{W_\delta} = \La_{\wh W_K}$ for $K = q^{-1}(\delta)$, hence by Proposition \ref{p:ical-Rm}, $\La_I \In \La_{W_\delta, 0}$ if and only if $K \cap R_I \neq \emptyset$, which is equivalent to $\delta \in q(R_I)$.  
\epf 

\bd 
\begin{enumerate}
    \item Let $\scal^\delta$ denote the stratification of $\R^k$ induced by $\delta \mapsto \La_{W_\delta}$. 
    \item For any $v \in \Z^k$, let $\scal^\delta_v$ denote the stratification of $\R^k$ induced by $\delta \mapsto \La_{W_\delta, v}$. 
\end{enumerate}
\ed 
By definition, the strata of $\scal^\delta$ are generated by intersecting $A_{I, v}$ for various $I \in [N]$, $v \in \Z^k$, and the strata of $\scal_v^\delta$ are generated by intersecting $A_{I,v}$ for various $I$ but fixed $v$. 

\subsection{Obstruction of fully-faithfulness for the restriction functor}
Here we consider the singular support of push-forward of hom sheaves, $\mu_* \uhom(F, G)$ for $F, G \in Sh^\dm(\R^N, \La_{W_\delta, v})$, which is bounded by $\mu_*SS_{Hom}(\La_{W_\delta})$ (See section \ref{ss:SS-cat}
). 

Let $\Sigma_{\nabla}$ be the exterior conormal fan of $\nabla$. For a cone $\sigma \in \Sigma_{\nabla}$, let $F_\sigma$ denote the corresponding (closed) face whose exterior conormal is $\sigma$. Let $\Aff(F_\sigma)$ be the affine hull of $F_\sigma$, i.e, the minimal affine linear space containing $F_\sigma$. 

For any nonzero $\xi \in (\R^k)^\vee$, we split $[N]$ into three parts 
$$ [N] = [N]_{\xi, -} \sqcup [N]_{\xi,0} \sqcup [N]_{\xi,+}, \quad [N]_{\xi,s} = \{ i \in [N] \mid \sign(\la \beta_i, \xi \ra) = s \}.  $$
Clearly, the splitting is a locally constant function on $\Sigma_\nabla$, hence we also label $[N]_{\xi, s}$ as $[N]_{\sigma, s}$ if $\Int(\sigma) \ni \xi$.


\bt \label{t:SS-Hom}
(1) For any nonzero $\xi \in (\R^k)^\vee$, let $\sigma$ denote the cone of the exterior normal fan $\Sigma_\nabla$ such that $\Int(\sigma) \ni \xi$. Then $\xi \in \mu_*(SS_{Hom}(\La_{W_\delta, v}))_0$ if and only if $\delta \in v + \Aff(F_\sigma)$. 

(2) For any $v \in \Z^k$, $\delta \in \R^k$, we have
$$ \mu_*(SS_{Hom}(\La_{W_\delta, v})) = T_{\R^k}^*\R^k \quad \cup \,  \bigcup_{\stackrel{ 0 \neq \sigma \in \Sigma_{\nabla}}{ \delta \in v + \Aff(F_\sigma)}} \sigma^\perp \times \sigma. $$ 

(3) For any $\delta \in \R^k$, 
$$ \mu_*(SS_{Hom}(\La_{W_\delta})) = T_{\R^k}^*\R^k \quad \cup \, \bigcup_{\sigma: (\delta + F_\sigma) \cap \Z^k \neq \emptyset} \Aff( \delta + F_\sigma) \times (-\sigma) $$
\et 

\bpf 
(1) Let $\wt \xi = \mu^*(\xi)$, and $s(\xi)=\sign(\wt \xi)$. Suppose $\xi \in \mu_*(SS_{Hom}(\La_{W_\delta}))_0$, then there exists $I, J \in \ical(\delta, v)$, such that $s(\xi) \leq \sign(I, J)$. This means, for any $i \in [N]$, if $\sign(\wt \xi)_i = +$, then $s(I, J)_i = +$, and if  $\sign(\wt \xi)_i = -$, then $s(I, J)_i = -$. Using Lemma \ref{l:ss-uhom-IJ}, this is equivalent to
$$ I^c  \supset [N]_{\xi, +}, \quad \z{and} \quad I \RM J \supset [N]_{\xi,-}, $$
and in turn, equivalent to
$$ [N]_{\xi, -} \In I \In  [N]_{\xi, -} \sqcup [N]_{\xi, 0}, \quad \z{and} \quad J \In  [N]_{\xi, +} \sqcup [N]_{\xi, 0}. $$
Since $\delta - v \in A_I \cap A_J$, and 
$$ A_I \subset -\beta_{[N]_{\xi, -}} - C_{[N]_{\xi, -} \sqcup [N]_{\xi, 0}} + \nabla = -\beta_{[N]_{\xi, -}}/2 - C_{[N]_{\xi, -} \sqcup [N]_{\xi, 0}}$$
$$ A_J \subset - C_{[N]_{\xi, +} \sqcup [N]_{\xi, 0}} + \nabla = \beta_{[N]_{\xi, +}}/2 - C_{[N]_{\xi, +} \sqcup [N]_{\xi, 0}} $$
Since $\beta_{[N]_{\xi, +}} = -\beta_{[N]_{\xi, -}}$, we have 
$$ A_I \cap A_J \subset \beta_{[N]_{\xi, +}}/2 - C_{[N]_{\xi, 0}} = \Aff(F_\sigma).$$
Thus $\delta \in v + \Aff(F_\sigma)$. 

In the other direction, if $\delta \in v + \Aff(F_\sigma)$, we can choose $I = [N]_{\xi, -} \sqcup [N]_{\xi, 0}$ and $J=[N]_{\xi, 0}$, then $A_I \cap A_J = \Aff(F_\sigma)$, hence $\delta \in A_{I, v} \cap A_{J, v}$, i.e. $I, J \in \ical{\delta, v}$. Thus $\wt \xi \in SS(\uhom(P_I, P_J))_0$ and $\xi \in \mu_*(SS_{Hom}(\La_{W_\delta}))$. 

(2) If $\delta - v \in \Aff(F_\sigma)$, then $\sigma^\perp \times \sigma \In \mu_*(SS_{Hom}(\La_{W_\delta}))$. If $\delta - v \in \Aff(F_\sigma)$ and $0 \neq \sigma' \subset \sigma$, then $\Aff(F_\sigma) \In \Aff(F_{\sigma'})$. Finally, since $[N] \in \ical(\delta, v)$, the zero section $T_{\R^k}^* \R^k$ is included. Hence the direction $\supset$ is proven. The other direction can also be easily checked. 

(3) One can verify the statement near each lattice point $v \in \Z^k$, and the rest follows from extrapolation. 
\epf

\section{Local Window Skeleton and the Co-restriction Functor \label{s:SS-Hom-L} }
In this section we prove
\bt \label{t:SSL-Hom-window}
For any $\delta \in \R^k$, we have $$  SS_{Hom}^L(\mu_* Sh^\dm_{\La_{W_\delta}}) = T^*_{\R^k} \R^k. $$
\et 

In fact, we will factorize $\mu$ as $\R^N \xto{\pi} \R^m \xto{q} \R^k$ and prove the following general result about $\pi$. 

\bt \label{t:SSL-Hom-mezzaine}
For any closed convex polytope $K \In \R^m$, we have $$ SS_{Hom}^L(\pi_* Sh^\dm_{\La_{W_\delta}})  = T^*_{\R^m} \R^m, $$ 
where $\wh W_K = (K + \pi([0,1]^N)/2) \cap \Z^m$, $\wt W_K = \pi^{-1}(\wh W_K) \cap \Z^N$ and $\La_{\wh W_K} = \cup_{w \in \wt W_K} \La_w$. 
\et 

\bpf [Proof of Theorem \ref{t:SSL-Hom-window}]
Take $K = q^{-1}(\delta)$ in Theorem \ref{t:SSL-Hom-mezzaine}, we get $\La_{W_\delta} = \La_{\wh W_K}$. Then 
$$ T^*_{\R^k} \R^k \In SS_{Hom}^L(\mu_* Sh^\dm_{\La_{W_\delta}}) \subset q_* SS_{Hom}^L(\pi_* Sh^\dm_{\La_{W_\delta}}) = q_* (T^*_{\R^m} \R^m) = T^*_{\R^k} \R^k. $$
\epf

We will first study the probe sheaf for the localized window skeleton $\La_{\wh W_K, \wh v}$ in detail. We will use notation from Section \ref{ss:Rm}. Since we are mainly going to work with $\pi$, we will drop the hat in $\wh W_K$ and $\wh v$ and hope there is no confusion incurred.

We recall some notation for ease of looking up. $\pi: \R^N \to \R^m$, $[N] = [N_1] \sqcup [N_2] \sqcup \cdots [N_m]$, and for each $i \in [m]$, $[N_i] = [N_i]_+ \sqcup [N_i]_+$. Let $\B = \pi([0,1]^N)/2$, then by quasi-symmetric condition, $\B=\prod_{i=1}^m [-\eta_i/2, \eta_i/2]$. If $I \In [N]$, we have decomposition of $I$ into parts $I_i = I \cap [N_i]$, and $I_i = I_{i,+} \sqcup I_{i,-}$ by intersecting with $[N_i]_\pm$. We say $I_i$ is of mixed type if both $I_{i, \pm}$ are non-empty, and we say $I_i$ is of pure-type if otherwise. We say $I$ is of mixed type if any $I_i$ is of mixed type, and of pure type if all $I_i$ are pure type. If $I_i \In [N_i]$ is of pure-type, we define $\sign(I_i)$ to be $0$ if $I_i =\emptyset$, and $\pm$ if $\emptyset \neq I_i \In [N_i]_\pm$. We also define $[N_i]_0 = \emptyset$. If $I$ is of pure-type, then we define sign vector $\sign(I) = (\sign(I_i))_i$. 

\subsection{Thick Coordinate Cones and Convex Polytope}
Recall that $\Sigma_m$ is the fan in $\R^m$ with $2^m$ quadrants as top dimensional cones, and for each sign vector $s \in \bS_m=\{+,-,0\}^m$ we let $\sigma_s$ denote the corresponding cone in $\Sigma_m$. 

We defined the following closed convex rectilinear regions, which serves as a fattened version of $\sigma_s$, 
$$ \B_s := \B + \sigma_s + \delta_s, \quad \delta_s := (s_1 \eta_1, \cdots, s_m \eta_m). $$ 
We also denote $\B_{s_i}$ as the projection of $\B_s$ to the $i$-th factor.

For any closed convex polytope $K$, we define a sub-poset $\bS_{K} \In \bS_m$,  
$$ \bS_{K} = \{s \in \bS_m \mid \B_s \cap K \neq \emptyset \}. $$

We would like to relate the topology of the geometric realization of the order complex of $\bS_K$ to that of $K$. One possible way is to pick a point $p_s \in K \cap \B_s$ for each $s \in \bS$, then glue the resulting simplices with vertices $p_s$, however this would make the geometric realization non-convex in general. Instead, for each $s \in \bS_K$, we take the constant sheaf $\C_{K_s}$ supported on the closed set $K_s : = K \cap \B_s$.  
Then for each $k$-simplex $[s]=[s_0<s_1<\cdots<s_k] \in \Delta(\bS_K)_k$, we consider the intersection $K_{[s]}=K_{s_0} \cap \cdots \cap K_{s_k}$. Thus, we have a \v{C}ech resolution of $\C_K$, (abusing notation and denote $\C_K$ by $K$, $\C_{K_[s]}$ by $K_{[s]}$),
$$ K \cong \bigoplus_{[s] \in \Delta_0} K_{[s]} \to \bigoplus_{[s] \in \Delta_1} K_{[s]} \to \cdots, $$
where $\Delta = \Delta(\bS_K)$. 

\begin{example}
Throughout this subsection, we will consider the following example,  where $m=2$ and $\eta_1 = 2, \eta_2=2$. The $\B_s$ decomposes $\R^2$ into 9 regions with disjoint interiors. The convex set $K$ intersets with the the following regions: $\B_{++}, \B_{+0}, \B_{+-}, \B_{0-}$.

\begin{figure}[hh]
    \centering
\begin{tikzpicture}[scale=0.7]
\draw [opacity=.1] (-5,-5) grid (5,5);
\draw [fill, opacity=0.3] (-1,-1) rectangle (1,1);
\draw [fill, opacity=0.3, yellow] (1,-1) rectangle (5,1);
\draw [fill, opacity=0.3, yellow] (-1,1) rectangle (-5,-1);
\draw [fill, opacity=0.3, yellow] (-1,1) rectangle (1,5);
\draw [fill, opacity=0.3, yellow] (1,-1) rectangle (-1,-5);
\draw [thick, fill=green,  opacity = 0.3] (3,2) -- (4,0) -- (2,-4) -- (0,-4) -- cycle; 
\draw [thick] (1,-2) -- (1,-4)  (1.5,-1) -- (3.5, -1) (2.5,1) -- (3.5, 1);
\node at (0,-3) {$K_{[0-]}$};
\node at (2,-2) {$K_{[+-]}$};
\node at (3,0) {$K_{[+0]}$};
\node at (3,3) {$K_{[++]}$};
\end{tikzpicture}
    \caption{The intersection pattern of $K$ with $\B_s$. }
    \label{fig:BsKs}
\end{figure}
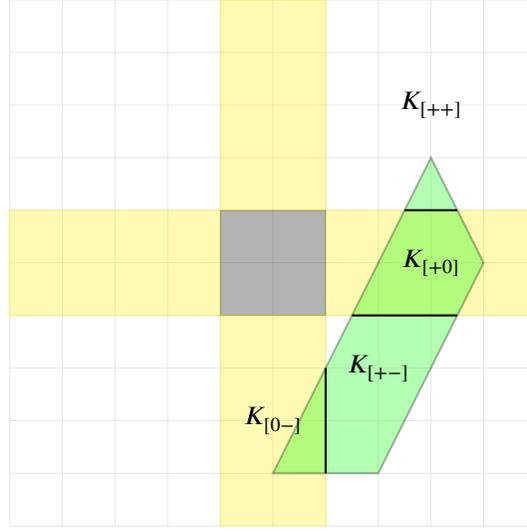

The sub-poset $\bS_K$ is then as following
$$ 
\begin{tikzcd}
(+0) \ar[r] \ar[rd] & (++) \\
(-0) \ar[r] & (+-)
\end{tikzcd}
$$
Then the order complex $\Delta(\bS_K)$ has degree $0$ corresponding to the vertices and degree $1$ corresponding to the arrows 
$$ \bcs 
\Delta(\bS_K)_0 = \{ +0, -0, ++, +- \} \\ 
\Delta(\bS_K)_1 = \{[+0 , ++], [+0 , +-], [0-, +-] \} \\ 
\ecs 
$$ 
The \v{C}ech resolution is then
$$ K_{[++]} \oplus  K_{[+0]} \oplus K_{[+-]} \oplus K_{[0-]} \to K_{[+0 , ++]} \oplus K_{[+0 , +-]} \oplus K_{[0-, +-]}, $$
where for example $K_{[++]} = K \cap \B_{++}$, and $K_{[+0, ++]} = K_{[+0]} \cap K_{[++]}$ is the interface. 
\end{example}

We define the analog of $\z{star}(\sigma_s)$ for the thick cone $\B_s$:
$$ \B_{\geq s} := \bigcup_{s' \geq s} \B_{s'}. $$
\bp \label{p:HomSKsSK}
For any $s \in \bS_m$, we have 
\bee  
\Hom_{\bS_K}((\bS_K)_{\geq s}, \bS_K) &\cong& \bigoplus_{[t] \in \Delta_0} 1_{t_0 \geq s} \cdot \C \to \bigoplus_{[t] \in \Delta_1}  1_{t_1 \geq s} \cdot \C \to \cdots \label{e:SKsSK1}\\
& \cong &  \Hom_{Sh(\R^m)}(\B_{\geq s}, K) \label{e:SKsSK2} \eee
where $\Delta = \Delta(\bS_K)$, and $1_{\z{condition}}$ equals $1$ or $0$ depending on whether the condition is satisfied.  In particular, we have 
$$ \Hom(\bS_K, \bS_K) \cong \C. $$
\ep 
\bpf
The proof follows from Proposition \ref{p:hom-of-subposets}, where we note that if $P_1$ is a sub-poset of $P$ given by $P_1 = P_{\geq a}$ for some element $a \in P$, then for any k-chain $[p]$ in $P$, $[p] \cap P_1 \neq \emptyset$ if and only if $p_k \geq a$. And the last statement follows from $H^*(K) \cong \C$, since $K$ is contractible. 
\epf
\begin{example}
We compute the hom spaces in following two cases:  $s=(+0)$ and $s=(+-)$. See Figure \ref{fig:ex-hom-bs-k}.

(1) $s=(+0)$ case. We first compute using combinatorial complex Eq \eqref{e:SKsSK1}. We have the following cochain complex involving two degrees (first row is degree 0, second row is degree 1, and each node reprsent $\C$.) One can see that this complex is acyclic.
\begin{center}
\begin{tikzcd}
(++) \ar[d] & (+0) \ar[d] \ar[ld]  & (+-) \ar[ld] \ar[d] \\
(+0, ++) &  (+0 , +-) &  (0-, +-)
\end{tikzcd}
\end{center}
We can then compute using Eq. \eqref{e:SKsSK2}. 
$$ \Hom(\B_{\geq s}, K) = \Gamma(K \otimes (\B_{\geq s})^\vee ) \cong H^*(K \cap \B_{\geq s}, K \cap \pa(\B_{\geq s})) \cong 0. $$
The hair on the boundary shows the direction of singular support of the sheaf. 
One can construct a Morse function $f$, such that the downward gradient of $f$ is pointing inward along the inward hair, and outward along the outward hair, and there exists such an $f$ without critical point, hence the cohomology relative to the boundary with outward hair is zero. See \cite[Section 4.5]{NZ}.  

(2) $s=(+-)$ case. The $\B_{\geq s}$ is the closed lower-right quadrant, and the sheaf $K \otimes (\B_{\geq s})^\vee$ is of the following form
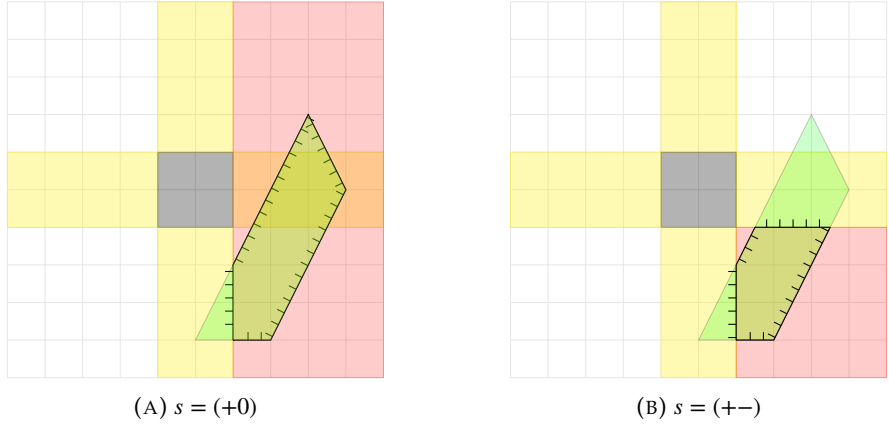
\begin{figure}
    \centering
      \begin{subfigure}[b]{0.4\textwidth}
         \centering
\begin{tikzpicture}[scale=0.5]
\draw [opacity=.1] (-5,-5) grid (5,5);
\draw [fill, opacity=0.3] (-1,-1) rectangle (1,1);
\draw [fill, opacity=0.3, yellow] (1,-1) rectangle (5,1);
\draw [fill, opacity=0.3, yellow] (-1,1) rectangle (-5,-1);
\draw [fill, opacity=0.3, yellow] (-1,1) rectangle (1,5);
\draw [fill, opacity=0.3, yellow] (1,-1) rectangle (-1,-5);
\draw [fill, opacity=0.2, red] (1,-5) rectangle (5,5); 
\draw [righthairs] (1,-2) -- (3,2) -- (4,0) -- (2,-4) -- (1,-4); 
\draw [righthairs] (1,-2) -- (1,-4); 
\draw [fill=green,  opacity = 0.2] (3,2) -- (4,0) -- (2,-4) -- (0,-4) -- cycle; 
\draw [fill=yellow, opacity =0.2] (1,-2) -- (3,2) -- (4,0) -- (2,-4) -- (1,-4);
\end{tikzpicture}
         \caption{$s=(+0)$}
         \label{fig:s=+0}
     \end{subfigure}
 \begin{subfigure}[b]{0.4\textwidth}
         \centering
\begin{tikzpicture}[scale=0.5]
\draw [opacity=.1] (-5,-5) grid (5,5);
\draw [fill, opacity=0.3] (-1,-1) rectangle (1,1);
\draw [fill, opacity=0.3, yellow] (1,-1) rectangle (5,1);
\draw [fill, opacity=0.3, yellow] (-1,1) rectangle (-5,-1);
\draw [fill, opacity=0.3, yellow] (-1,1) rectangle (1,5);
\draw [fill, opacity=0.3, yellow] (1,-1) rectangle (-1,-5);
\draw [fill, opacity=0.2, red] (1,-1) rectangle (5,-5); 
\draw [fill=green,  opacity = 0.2] (3,2) -- (4,0) -- (2,-4) -- (0,-4) -- cycle; 
\draw [fill=yellow, opacity =0.2] (1.5,-1) -- (3.5,-1) -- (2,-4) -- (1,-4) -- (1,-2);
\draw [righthairs] (1,-2) -- (1.5,-1) (3.5,-1)  -- (2,-4) -- (1,-4); 
\draw [righthairs] (1,-2) -- (1,-4) (3.5,-1) -- (1.5,-1); 
\end{tikzpicture}
         \caption{$s=(+-)$}
         \label{fig:s=+-}
     \end{subfigure}
    \caption{Computation of $\Hom(\B_{\geq s}, K)$, $\B_{\geq s}$ is shaded red, and $K$ is shaded in green.}
    \label{fig:ex-hom-bs-k}
\end{figure}

hence we have
$$ \Hom(\B_{\geq s}, K) \cong H^*(K \cap \B_{\geq s}, K \cap \pa(\B_{\geq s})) \cong \C[-1] $$. 
\end{example}

We also record a lemma that will be used later. Recall $\z{star} \sigma_s = \cup_{\sigma_{s'} \supset \sigma_s} \sigma_{s'}$. 
\bl  \label{l:Hom-Bs-K}
If for some $s \neq 0$ and for some convex set $K$, we have $\Hom_{\bS_K}((\bS_K)_{\geq s}, \bS_K) \neq 0$, then $(K-\delta_s/2) \cap -(\B_{\geq s} - \delta_s/2) =\emptyset$, or equivalently $ K \cap (-\z{star}\sigma_s + \B) = \emptyset$. 
\el 

\bpf 

We can compute $\Hom(\B_{\geq s}, K)$ by pushforward along the map $\pi_s = \pi_{\sigma_s}: \R^m \to \span(\sigma_s) = \R^{|s|}$. Indeed, $\B_{\geq s} =(\pi_s)^* (\pi_s)_* \B_{\geq s}$, hence 
$$ \Hom_{\bS_K}((\bS_K)_{\geq s}, \bS_K) \cong \Hom_{\R^m} (\B_{\geq s}, K) \cong  \Hom_{\R^{|s|}} ( (\pi_s)_*\B_{\geq s}, (\pi_s)_* K)) \cong \Hom_{\R^{|s|}} ( \pi_s \B_{\geq s}, \pi_s K) $$
where the middle term is to push-forward the constant sheaves, and the last term is claiming the push-forward sheaves equal to the constant sheaves on the image of $\pi_s$, since the fibers of $\pi_s$ on $\B_{\geq s}$ and $K$ are convex thus contractible.

Let $B = \pi_s (\B_{\geq s} - \delta_s/2)$, $C = \pi_s (K- \delta_s/2)$, where the shift is to make sure the quadrant $B$ has vertex at $0$. Then from previous arguments we have 
\be \label{e:nonzero} H^*(B \cap C, \pa B \cap C) \cong \Hom(\C_B, \C_C) \cong \Hom_{\R^{|s|}} ( \pi_s \B_{\geq s}, \pi_s K) \neq 0. \ee 

The claim is equivalent to $-B \cap C = \emptyset$. Suppose not, then let $x \in -B \cap C$, and let 
$$ S_x = \{ y \in \pa B \cap C \mid \z{ the segment $\wb{xy}$  intersects $B$ only at $y$} \},$$ 
then both $B \cap C$  and $\pa B \cap C$ deformation retract to $S_x$ along rays emitted from $x$, hence they are all contractible. Hence the relative cohomology $H^*(B \cap C, \pa B \cap C) = 0$, which contradicts with Eq \eqref{e:nonzero}. 
\epf

\subsection{Sign Pattern and Probe Sheaf}
Recall that for $v \in \Z^m$, we have the localized window skeleton 
$\La_{W_K, v}$ and index set $\ical(K, v)$ containing those subset $I\In [N]$ such that $\La_I \In \La_{W_K, v}$. From Proposition \ref{p:ical-Rm}, we know if $I \in \ical(K, v)$ if and only if $(K + \pi(e_I) + \pi(\tau_I) ) \cap (v + \B) \neq \emptyset$. 

From the general discussion of probe sheaf in Section \ref{ss:probe-general}, in particular Proposition \ref{p:probe-as-holim}, we can express $F_I$ as a homotopy limit, or more concretely
$$ F_I: = F_{I, \ical(K,v)} \cong  \bigoplus_{[J] \in \Delta_0} P_{J_0} \to \bigoplus_{[J] \in \Delta_1 } P_{J_1} \to \cdots, $$
where $\Delta$ is the order complex of the poset $\ical(K, v)_{\geq I}$.  We can also compute which quadrant is in the support of $F_I$ by computing the hom of sub-posets, c.f Proposition \ref{p:hom-FI-FJ}. 

Our goal in this section is to find a simpler presentation of $F_I$, using the sign pattern of some convex set $K_{I,v}$. 

For any $I \In [N]$, let $$ K(I) = K + \pi(e_I) + \pi(\tau_I). $$ 

\bd \label{d:SKIV}
For any $I \In [N], v \in \Z^m$ and convex polytope $K \In \R^k$, we define the sign pattern $\bS_{K,I,v} \In \bS_m$, where $s \in \bS_{K, I, v}$ if and only if
\begin{enumerate}
    \item $K(I) \cap (v + \B_s) \neq \emptyset$. 
    \item For each $i \in [m]$, if $I_{i, +} \neq \emptyset$ then $s_i \neq -$ , and if  $I_{i, -} \neq \emptyset$ then $s_i \neq +$.
\end{enumerate}
We also define $\Sigma_{K, I, v} = \{\sigma_s \mid s \in \bS_{K, I, v}\}$. 
\ed 
\brem 
The condition (2) is equivalent to the following list 
\begin{itemize}
    \item If $I_i = \emptyset$, then $s_i \in \{-,+,0\}$. 
    \item If $\emptyset \neq I_i \In [N_i]_+$, then $s_i \in \{-, 0\}$. \item If $\emptyset \neq I_i \In [N_i]_-$, then $s_i \in \{+, 0\}$.
    \item If $I$ is of mixed type, then $s_i = 0$. 
\end{itemize}
If we drop condition (2), one can still define the simplified resolution of $F_I$ (c.f Eq \eqref{e:FI-better}),  but it would involve sheaves not admissible for the skeleton $\La_{W_K, v}$. 
\erem 

\bp \label{p:KIv}
There exists a closed convex polytope $K_{I,v}$, such that 
$$ s \in \bS_{K, I, v} \LRA K_{I,v} \cap \B_s \neq \emptyset. $$
Furthermore, $K_{I,v}$ satisfies the property that
$$  K_{I,v} \cap \B_s \neq \emptyset \RA K_{I,v} \cap \Int(\B_s) \neq \emptyset. $$
\ep 
\bpf
We define $K_{I,v, \epsilon} = (K(I)-v) \cap B_{I,\epsilon}$ for small enough positive $\epsilon$, where 
$$ B_{I, \epsilon} = \bigcap_{i=1}^m \bigcap_{s_i=+,-} B_{I_i, s_i, \epsilon}, \quad B_{i, s_i, \epsilon} = \bcs \{x \in \R^m \mid  s_i x_i \geq  -\eta_i/2 + \epsilon  \} & I_{i,s_i} \neq \emptyset \\ \R^m & I_{i,s_i} = \emptyset \ecs.  $$
We will consider small enough $\epsilon$ and drop $\epsilon$ in the notation $K_{I,v, \epsilon}$.  

For the direction $\LA$, if $K_{I,v} \cap (v + \B_s) \neq \emptyset$, then $K(I) \cap (v + \B_s) \neq \emptyset$ since $K(I) \supset K_{I,v}$, hence cnodition (1) is satisfied. For each $i \in [m], s_i' \in \{+,-\}$, if $I_{i,s_i'} \neq \emptyset$, then $s_i \neq -s'_i$ since otherwise $B_{I_i, s_i,  \epsilon} \cap \B_{s_i} = \emptyset$, hence cnodition (2) is satisfied. Thus $\LA$ is proven. 

For the direction $\RA$, if $\sigma_s \in \Sigma_{K, I, v}$, we just need to prove that $(K(I)-v) \cap \B_s \cap \Int(B_{I, 0}) \neq \emptyset$, since $\Int(B_{I,0}) = \cup_{\epsilon >0} B_{I,\epsilon}$. Set $v=0$ for simplicity. Since $\B_s \In B_{I,0}$, we have $(K(I)-v) \cap \B_s \cap B_{I,0} \neq \emptyset$. Let $x \in K(I) \cap \B_s$, then $x + \epsilon \sum_{i, s'_i: x \in \pa B_{i,s'_i,0}} s'_i E_i$ will be in the interior of $B_{I,0}$ and remain in $K(I)$.

We may replace $K_{I,v}$ by  $K_{I,v} + [-\delta, \delta]^m$ for small enough $\delta$ to achieve the second property. 
\epf

\bl \label{l:PIJ}
If $s \in - \bS_{K, I, s}$,  and $J \In [N]$ is of pure-type $s$, then $J \cap I = \emptyset$ and $I \sqcup J \in \ical(K, v)$. 
\el 
\bpf 
From condition (2) in the definition of $\Sigma_{K, I, s}$, we see $I \cap J = \emptyset$. Also note that $K(I) \cap (v - \B_s) \neq \emptyset$, and $\pi(\tau_J) = \sigma_s$, 
$$ K(I) \cap (v + \B + \sigma_{-s} + \delta_{-s}) \neq \emptyset  \LRA (K(I) +  \pi(\tau_J) + \delta_s)  \cap (v + \B) \neq \emptyset. $$
and
\[ K(I \sqcup J) = K + \pi(\tau_I + \tau_J) + \pi(e_I) + \pi(e_{J}) \supset K + \pi(\tau_I + \tau_J) + \pi(e_I) + \delta_{s} =K(I) +  \pi(\tau_J) + \delta_s. \]
Hence $K(I \sqcup J) \cap (v + \B) \neq \emptyset$, and $I \sqcup J \in \ical(K, v)$. 
\epf 

\bp 
The probe sheaf $F_I$ is quasi-isomorphic to the following complex 
\be \label{e:FI-better} \wt F_I := \uwave{\bigoplus_{[\sigma] \in \Delta_0} F_{I, \sigma_0}} \to \bigoplus_{[\sigma] \in \Delta_1} F_{I, \sigma_1} \to \cdots , \quad \Delta = \Delta(-\Sigma_{K, I, v}), \quad [\sigma]=(\sigma_0< \cdots < \sigma_k) \in \Delta_k. 
\ee
where 
\begin{itemize}
    \item For each $\sigma \in - \Sigma_{K,I,v}$, $F_{I, \sigma}$ is a sheaf in $Sh(\prod_{i=1}^m \R^{[N_i]}, \La_{W_K, v})$
    $$ 
F_{I, \sigma} = F_{I_1, s_1} \boxtimes F_{I_1, s_1} \boxtimes \cdots \boxtimes F_{I_m, s_m}, \quad \sigma = \sigma_s, \; s= (s_1, \cdots, s_m), 
$$
and for each $i \in [m]$,  $F_{I_i, s_i}$ is a sheaf in $Sh(\R^{[N_i]_+} \times \R^{[N_i]_-})$
\begin{center}
\begin{tabular}{c|c|c|c|c}
 $F_{I_i, s_i}$ & $I_i = \emptyset$ & $\emptyset \neq I_i = I_{i,+}$ & $\emptyset \neq I_i = I_{i,-}$ & I is mixed \\
\hline
$s_i=0$ & $P_{\emptyset} \boxtimes P_{\emptyset}$ & $P_{I_{i,+}} \boxtimes P_{\emptyset}$ & $P_{\emptyset} \boxtimes P_{I_{i,-}}$ & $P_{I_{i,+}} \boxtimes P_{I_{i,-}}$ \\
\hline
$s_i=+$ & $\hourglass_{[N_i]_+} \boxtimes P_{\emptyset}$ & ----- & $\hourglass_{[N_i]_+} \boxtimes P_{I_{i,-}}$ & ----- \\
\hline
$s_i=-$ & $P_{\emptyset} \boxtimes \hourglass_{[N_i]_-}$ & $P_{I_{i,+}} \boxtimes \hourglass_{[N_i]_-}$ & ----- & ----- \\
\hline
\end{tabular}
\end{center}
  \item For each $\sigma < \sigma'$, the morphism $F_{I, \sigma'} \to F_{I,\sigma}$ is a product of $F_{I_i,s_i'} \to F_{I_i, s_i}$ for $s_i' < s_i$, using the inclusion $P_\emptyset \into \hourglass$ and identity maps. 
\end{itemize}
\ep 
\bpf
By the Koszul resolution of the hourglass sheaves $\hourglass$ (Section \ref{ss:hourglass}), we note that $F_{I,\sigma_s}$ can be resolved using sheaves of type $P_{I \sqcup J}$ where $J$ is pure of type $s$. Since $\sigma_s \in - \Sigma_{K, I, v}$, by Lemma \ref{l:PIJ}, we have $SS(P_{I\sqcup J}) \In \La_{W_K, v}$, hence $SS(F_{I,\sigma_s}) \In \La_{W_K, v}$. 

Next, suffice to prove that for any $J \in \ical(K, v)$, we have 
we have 
$$ \Hom(\wt F_I, P_J) = \bcs 
0 & \z{ if } I \not \In J \\
\C & \z{ if } I \In J 
\ecs. 
$$

If $I \not \In J$, then for any $I' \supset I$ we have $I' \not \In J$.  Since $\wt F_I$ can be resolved with $P_{I'}$ with $I' \supset I$, we have $\Hom(F_I, P_J) = 0$. 

If $I \In J$ and $J \in \ical(W_K,v)$, then $K(J) \cap (v + \B) \neq \emptyset$. Let $J' = J \RM I$, then $K(J) = K(I) + \pi(\tau_{J'}) + \pi(e_{J'})$, and $K(J) \cap (v + \B) \neq \emptyset$, which is equivalent to $K(I) \cap (v + \B - \pi(\tau_{J'}) - \pi(e_{J'})) \neq \emptyset$. Since $(v + \B - \pi(\tau_{J'}) - \pi(e_{J'})) \In v + \B - \pi(\tau_{J'})$,  hence 
$$ K(I) \cap (v + \B - \pi(\tau_{J'})) \neq \emptyset. $$

Let $\Sigma_{J'} : = \{ \sigma \in \Sigma_m \mid \sigma \in -\pi(\tau_{J'})\}$, then $\Sigma_{J'}$ closed under taking sub-cone. We then define 
$$ \Sigma^{J'}_{K,I,v} = \Sigma_{J'} \cap \Sigma_{K,I,v}. $$

\bl 
Let $K_{I,v}$ be given as in Proposition \ref{p:KIv}, and let $K^{J'}_{I,v} = K_{I,v} \cap (\B - \pi(\tau_{J'}))$, then $$ \Sigma^{J'}_{K,I,v} = \{ \sigma_s \mid  K^{J'}_{I,v} \cap \Int \B_s \neq \emptyset\}. $$ 
\el 
\bpf 
We note that
$$ K_{I,v} \cap (\B - \pi(\tau_{J'})) \cap \Int \B_s \neq \emptyset \LRA \bcs K_{I,v} \cap \Int \B_s \neq \emptyset \\  \Int \B_s \In \B - \pi(\tau_{J'}) \ecs \LRA \bcs \sigma_s \in \Sigma_{K,I,v} \\ \sigma_s \in  \Sigma_{J'}. \ecs $$
\epf

\bl 
For any $\sigma \in -\Sigma_{K, I, v}$, we have
$$ \Hom(F_{I,\sigma}, P_J) = \bcs \C & \z{if } \sigma \in -\Sigma^{J'}_{K,I,v} \\
0 & \z{otherwise} 
\ecs. $$
\el 
\bpf 
Write $\Sigma^{J'}_{K,I,v}$ as $\Sigma'$ for short. 
If $\sigma_s \in -\Sigma'$, then suffice to check that for each $i \in [m]$, $\Hom(F_{I_i, s_i}, P_{J_i}) = \C$. If $F_{I_i, s_i} = P_{I_i}$, then since $I_i \In J_i$, this is verified. If $F_{I_i, s_i} = \hourglass_{[N_i]_+} \boxtimes P_{I_i}$ in the case of  $s_i=+$ and $I_i \In [N_i]_-$, then $J'_{i,+}  \neq \emptyset$, otherwise $(v + \B - \pi(\tau_{J'}))) \cap (v + \B_{-s}) = \emptyset$ by examining the projection to the $i$-th direction, thus contradicting with the requirement of $\Sigma'$. Thus
$$ \Hom(F_{I_i, s_i}, P_{J_i}) = \Hom(\hourglass_{[N_i]_+}, P_{J'_{i,+}}) \otimes \Hom(P_{I_i}, P_{J_{i,-}}) \cong \C. $$
The case for $s_i=-$ and  $I_i \In [N_i]_+$ is similar. 

If $\sigma_s \notin -\Sigma_{K, I, v} \RM -\Sigma'$, then  $\sigma_s \not \In \pi(\tau_{J'})$, thus there exists one $i \in [m]$, such that 
$\sigma_{s_i} \not \In \pi(\tau_{J'_i})$. Then $J'_i$ is not of mixed type, since that would make $\pi(\tau_{J'_i}) = \R$, and $s_i \neq 0$ since that would make $\sigma_{s_i}=0$. In all cases, $J'_{i, s_i} = \emptyset$. Take $s_i=+$ for example, then 
$$ \Hom(F_{I_i, s_i}, P_{J_i}) = \Hom(\hourglass_{[N_i]_+}, P_{J'_{i,+}}) \otimes \Hom(P_{I_i}, P_{J_{i,-}}) = 0 \otimes \C = 0, $$
where we used $\Hom(\hourglass, P_\emptyset) = 0$. The other case $s_i=-$ is similar.
\epf

\bl 
The geometric realization of $\Sigma^{J'}_{K,I,v}$ is contractible, hence its homology and cohomology are $\C$ in degree $0$.  
\el 
\bpf
If $\pi(\tau_J') = \R^m$ we let $K' = K^{J'}_{I,v}$, otherwise we pick a unit vector $v$ in the interior of $-\pi(\tau_J')$ and let $K' = \epsilon v + K^{J'}_{I,v}$. Then, one may check that $\Sigma^{J'}_{K,I,v} = \{ \sigma_s \mid K' \cap \B_s \neq \emptyset \}$ and $K'$ is a convex closed set, the result then follows from Proposition \ref{p:HomSKsSK}. 
\epf

Now we are ready to prove that if $I \In J$ and $J \in \ical(K,v)$, then $\Hom(\wt F_I, P_J) = \C$. Define $\Delta' = \Delta(-\Sigma^{J'}_{K,I,v})$. We have
$$ \Hom(\wt F_I, P_J) = \bigoplus_{[\sigma] \in \Delta'_0} \C \gets \bigoplus_{[\sigma] \in \Delta'_1} \C \gets \cdots = H_* (-\Sigma^{J'}_{K,I,v}) \cong \C.$$
\epf

Next, we use the better presentation of $F_I$ Eq \eqref{e:FI-better} to find the support of $F_I$. 
For any sign vector $s$, we define a maximal subset of pure-type $s$
$$ [N]_s = \sqcup_{i=1}^m [N_i]_{s_i} $$
where we set $[N_i]_0 =\emptyset$. 
We also define
$$ \deg(s) = \sum_{i, s_i \neq 0} (N_{i,s_i} -1), \quad N_{i,s} = |[N_i]_s|.$$

\bl 
For any $J \In [N]$ and any $\sigma_s \in - \Sigma_{K, I, v}$, 
$$ \Hom(Q_J, F_{I,\sigma}) = \bcs \C[-\deg(s')] &\z{ if $J \RM I = [N]_{s'}$ for some $s' \leq s$}  \\
0 & \z{ else }\ecs 
$$ 
\el 
\bpf 
This can be verified using the table of the definition of $F_{I_i, s_i}$. 

Let $J' = J \RM I$. Assume $\Hom(Q_J, F_{I,\sigma})$ is non-zero. If $s_i = 0$, then $F_{I_i,s_i} = P_{I_i}$, then one needs to have $J'_i = \emptyset$. If $s_i = +$, then $F_{I_i, s_i} = \hourglass_{[N_i]_+} \boxtimes P_{I_{i,-}}$, hence $J'_i = \emptyset$ or $[N_i]_+$. Similarly, if $s_i=-$, then $J'_i = \emptyset$ or $[N_i]_-$. Conversely, suppose $J'$ satisfies the condition, then note that for any $n \geq 1$ the hourglass sheaf $\hourglass_{[N]}$ is supported only on $Q_\emptyset$ with stalk $\C$ and $Q_{[n']}$ with stalk $\C[-n'+1]$, we have the desired value of the hom. 
\epf 


\bp \label{p:hom-QJ-FI-1}
For any $I, J \In [N]$, $\Hom(Q_J, F_I)$ is zero unless $J \RM I = [N]_s$ for some sign vector $s$, in which case
$$  \Hom(Q_J, F_I) \cong \Hom( (-\Sigma_{K, I, v})_{\geq \sigma_{s}} , -\Sigma_{K, I,v})[-\deg(s)] \cong \Hom(\B_{\geq s}, -K_{I,v})[-\deg(s)]. $$
In particular, if $J \In I$, then $ \Hom(Q_J, F_I)  = \C$. 

\ep 
\bpf 
If $J \RM I \neq [N]_s$ for some sign vector $s$, then all the terms in $\Hom(Q_J, F_I)$ vanish, hence we have the first claim.  

If $J \RM I = [N]_s$, then let $\Delta = \Delta(-\Sigma_{K,I,v})$, 
$$ \Hom(Q_J, F_I) = \bigoplus_{[\sigma] \in \Delta_0} 1_{\sigma_s \In \sigma_0}\C[-\deg(s)]  \to \bigoplus_{[\sigma] \in \Delta_1} 1_{\sigma_s \In \sigma_1} \C[-\deg(s)+1] \to \cdots $$ 
$$= \Hom( (-\Sigma_{K, I, v})_{\geq \sigma_{s}} , -\Sigma_{K, I,v})[-\deg(s)] $$
where $1_{...}$ is the indicator function, valued in $1$ or $0$. The next equality comes from Proposition \ref{p:HomSKsSK}. 

If $J \In I$, then $J \RM I = \emptyset$ and $s=0$, hence $\Hom(\B_{\geq 0}, -K_{I,v}) \cong H^*(-K_{I,v})) = \C$. 
\epf 

\subsection{Leaks and flooded quadrants for localized window skeleton}
Let $\Sigma \In \Sigma_m$ be a sub-poset. 
We define the upper closure and lower closure of $\Sigma$ in $\Sigma_m$ as
$$ \ceil{\Sigma} = \bigcup_{\sigma \in \Sigma}  (\Sigma_m)_{\geq \sigma}, \quad \floor{\Sigma} = \bigcup_{\sigma \in \Sigma}  (\Sigma_m)_{\leq \sigma}. $$
We also denote the underlying geometric representation of $\Sigma$ as
$$ | \Sigma | = \bigcup_{\sigma \in \Sigma} \Int(\sigma). $$
The reason we only use the relative interior of each cone in $\Sigma$ is that it does not lose any information of $\Sigma$. 

If $A \In \R^m$, we denote $\Sigma_A = \{ \sigma \mid \sigma \In A\}$ as the fan generated by $A$. 

Recall that a leak for $Q_I$ in $\La_\ical$ is of the form $-\tau_I + \tau_J$ for some $J$ disjoint from $I$ and satisfies $I \cup J' \notin \ical$ for any $J' \In J$. We say $J$ {\bf labels a leak for $Q_I$ (in $\La_\ical$)}.  We now give a bound for the $\pi(\tau_J)$. 

\bp 
(1) If $J$ labels a leak for $Q_I$,  then for any cone $\sigma_s \In \pi(\tau_J)$, we have  $\sigma_s \notin \ceil{-\Sigma_{K,I,v}}. $

(2) Conversely, if $\sigma_s$ is a cone such that $\sigma_s \notin \ceil{-\Sigma_{K,I,v}}$, then for any $J$ of pure-type $s$ and disjoint from $I$, $J$ labels a leak for $Q_I$. 
\ep 
\bpf
(1) We prove by contradiction. Suppose that $\sigma_s \in -\Sigma_{K,I,v}$ and $\sigma_s \In \pi(\tau_J)$, then we can lift $\sigma_s$ to a pure-type subset $J' \In J$. Then by Lemma \ref{l:PIJ}, $I \cup J' \in \ical(K,v)$ contradicting with the requirement of a leak (Lemma \ref{l:leak2}).

(2) It suffices to prove that, for any $\sigma_s \notin \ceil{-\Sigma_{K,I,v}}$ and any $J$ of pure type $s$ and disjoint from $I$, we have $I \sqcup J \notin \ical(K,v)$. For simplicity of notation, set $v=0$. We let $\Sigma_{K(I)} = \{\sigma_s \mid \B_s \cap K(I) \neq \emptyset\}$, then $\Sigma_{K,I,v} \In \Sigma_{K(I)}$. We may check that $\ceil{\Sigma_{K,I,v}} = \ceil{\Sigma_{K(I)}}$.  Then we have
\bea  I \cup J \notin \ical(K,v)  &\LRA&  K + \pi(\tau_I + \tau_J) + \pi(e_I + e_J) \cap \B = \emptyset \\
& \LA & (K(I) + \sigma_{s}) \cap \B = \emptyset \\
& \LRA & K(I) \cap (\B-\sigma_{s}) = \emptyset \\
& \LRA & \Sigma_{K(I)} \cap \Sigma_{\sigma_{-s}} = \emptyset \\
& \LRA & \ceil{\Sigma_{K(I)} } \cap  \{ \sigma_{-s} \} = \emptyset \\
& \LRA &  \sigma_{-s} \notin  \ceil{\Sigma_{K,I,v} } 
\eea 
\epf

\bc \label{c:no-leak-over-U}
For any subset $I \notin \ical(K,v)$ and open set $U \In \R^m$,  the following are equivalent: 
\begin{itemize}
    \item there is no leak of $Q_I$ over $U$, i.e. for all leak $L = -\tau_I + \tau_J$ of $Q_I$,  $\pi(L) \cap U = \emptyset$
    \item $U + \pi(\tau_I) \In | \ceil{-\Sigma_{K,I,v}} |.$
\end{itemize}
\ec
\bpf 
There is no leak of $Q_I$ over $U$ if and only if for all $s$ such that $\sigma_s \notin \ceil{-\Sigma_{K,I,v}}$, we have $(\pi(-\tau_I) + \sigma_s ) \cap U = \emptyset$, which is in turn equivalent to 
\bea 
& \LRA & \sigma_s  \cap (U + \pi(\tau_I)) = \emptyset, \quad \forall \sigma_s \notin \ceil{-\Sigma_{K,I,v}} \\
& \LRA & U + \pi(\tau_I) \In | \ceil{-\Sigma_{K,I,v}} |.
\eea 
\epf 

Next, we  bound the $\pi$-image of flooded quadrants of $F_I$. If $Q_J$ is flooded by $F_I$, then by Proposition \ref{p:hom-QJ-FI-1}, we have
$$ J = I' \sqcup J', \quad I' = I \cap J, \quad J' = J \RM I = [N]_s $$
for some sign vector $s$. 

\bp \label{p:bound-flooded-quadrant}
With the above notation, $Q_J$ is a flooded quadrant of $F_I$ only if 
$$  \ceil{ -\Sigma_{K,I,v}}  \cap \ceil{ \{\sigma_{-s}\}} = \emptyset $$
\ep  
\bpf 
$Q_J$ is a flooded quadrant if and only if 
$\Hom(\B_{\geq s}, -K_{I,v}) \neq 0$. 
From Lemma \ref{l:Hom-Bs-K}, a necessary condition is 
\bea 
\Hom(\B_{\geq s}, -K_{I,v}) \neq 0 &\RA& -K_{I,v} \cap (\B - \z{star} \sigma_s) = \emptyset \\
&\LRA& -\Sigma_{K,I,v} \cap \floor{\ceil{\{\sigma_{-s}\}}} = \emptyset \\
&\LRA&  \ceil{ -\Sigma_{K,I,v}}  \cap \ceil{ \{\sigma_{-s}\}} = \emptyset
\eea 
\epf


\bp \label{p:no-leak-no-flood}
For any subset $I \notin \ical(K,v)$, we have 
 $$  \wb{ \pi(Q_I)} \cap \pi( \z{Leak}(I) )= \wb{ \pi(Q_I)} \cap  \pi(\z{Flood}(I) ), $$ 
    where $\z{Leak}(I)$ is the union of leaks of $Q_I$, and $\z{Flood}(I)$ is the closure of the union of flooded regions. 
\ep 
\bpf 
Since the $\In$ direction is automatic by Proposition \ref{p:leak-flood},  it suffice to show that for any open subset $U \In \R^m$, such that $U \cap \pi(Q_I) \neq \emptyset$,  if there is no leak of $Q_I$ over $U$, then there is no flooded quadrant of $F_I$ over $U$. In other words, for any flooded quadrant $Q_J$, $\pi(Q_J) \cap U = \emptyset$. 

We partition $[m]$ into three parts, 
$$ [m] = [m]_{0} \sqcup [m]_{\pm} \sqcup [m]_{+-},  $$
where 
$$ \bcs 
[m]_0 = \{ i \in [m]:  I_i = \emptyset \}  \\
[m]_{\pm} = \{ i \in [m]: I_i \neq \emptyset, I_i \In [N_i]_+ \z{ or }  I_i \In [N_i]_- \} \\
[m]_{+-} = i \in [m]: I_i \cap [N_i]_+ \neq \emptyset \z{ and }  I_i \cap [N_i]_- \neq \emptyset \}
\ecs.
$$
For $\bu \in \{0, \pm , +-\}$, we define subspaces $\R^{[m]_\bu}$, and projections $\pi_{[m]_\bu}: \R^{[m]} \to \R^{[m]_\bu}$. 

Since there is no leak of $Q_I$ over $U$, we have $U + \pi(\tau_I) \In | \ceil{ -\Sigma_{K,I,v}} |$. Since both sides are translation invariant by $\R^{[m]_{+-}}$, we may quotient out $\R^{[m]_{+-}}$, and assume $[m]_{+-}=\emptyset$. Now that $\pi(\tau_I)$ is a proper cone, we note that $| \ceil{ -\Sigma_{K,I,v}} |$ is invariant by translation $-\pi(\tau_I)$, and $U + \pi(\tau_I))$ invariant by translation $\pi(\tau_I)$, hence we have
\be U + \pi(\tau_I)) \In | \ceil{ -\Sigma_{K,I,v} \cap \R^{[m]_0}  } |, \quad \z{ where } -\Sigma_{K,I,v} \cap \R^{[m]_0} = \{\sigma \in -\Sigma_{K,I,v} \mid \sigma \In \R^{[m]_0}\}.  \label{e:no-leak-no-flood}\ee
This also shows that $[m]_0 \neq \emptyset$, since otherwise $-\Sigma_{K,I,v} \cap \R^{[m]_0} \In \{\sigma_0\}$ the origin, and  $\sigma_0 \notin -\Sigma_{K,I,v}$ since $I \notin \ical(K,v)$. 

We also note that Eq \eqref{e:no-leak-no-flood} is equivalent to 
\be  \pi_{[m]_0}(U) \In  | \ceil{ \pi_{[m]_0}(-\Sigma_{K,I,v} \cap \R^{[m]_0})  }|. \label{e:no-leak-2} \ee

Suppose $Q_J$ is a flooded quadrant, with $J = I' \sqcup J'$, $I'= J \cap I$ and $J' = J \RM I = [N]_s$, then by Proposition \ref{p:bound-flooded-quadrant}, we have $\ceil{ -\Sigma_{K,I,v}}  \cap \ceil{ \{\sigma_{-s}\}} = \emptyset $. If $I_i = \emptyset$, then $J_i = [N_i]_{s_i}$, hence under $\pi_i: \R^{N_i} \to \R$, we have 
$$ \pi_i(Q_{[N_i]_{s_i}}) = \bcs \R & s_i = 0 \\ 
\R_{\leq 0} & s_i = + \\
\R_{\geq 0} & s_i = - 
\ecs \quad = \quad  \z{star} (\sigma_{-s_i})
$$ 
Hence we have
\bea  
\pi(Q_J) \cap U = \emptyset &\LRA & \pi(-\tau_{I'} - \tau_{[N]_s} + \tau_{J^c}) \cap U = \emptyset  \\
& \LA & \pi_{[m]_0} \pi(-\tau_{I'} - \tau_{[N]_s} + \tau_{J^c})) \cap \pi_{[m]_0}(U) = \emptyset \\
& \LRA & (-\z{star} (\pi_{[m]_0} \sigma_s)) \cap \pi_{[m]_0}(U) = \emptyset \\
& \LA & (-\z{star} (\pi_{[m]_0}  \sigma_s)) \cap  | \ceil{ \pi_{[m]_0}(-\Sigma_{K,I,v} \cap \R^{[m]_0})} |   = \emptyset \\
& \LRA & \ceil{ \{\pi_{[m]_0} \sigma_{-s}\} } \cap \ceil{ \pi_{[m]_0}(-\Sigma_{K,I,v} \cap \R^{[m]_0})} = \emptyset
\eea 
where in the next to last step, we used Eq \eqref{e:no-leak-2}. 

We note that, for any cone $\sigma_s \in \Sigma_m$, 
$$ \pi_{[m]_0} \ceil{ \{ \sigma_s \} } = \{ \pi_{[m]_0}(\sigma') \mid \sigma' \geq \sigma_s \} = \{ \sigma' \in \Sigma_{m_0} \mid \sigma' \geq \pi_{[m]_0} \sigma_s \} = \ceil{ \{ \pi_{[m]_0} \sigma_s\} }. $$
and similarly
$$ \ceil{ \pi_{[m]_0}(-\Sigma_{K,I,v} \cap \R^{[m]_0})} =  \pi_{[m]_0} \ceil{  -\Sigma_{K,I,v} \cap \R^{[m]_0}}. $$
Hence, we have
\bea  && \ceil{ \{\pi_{[m]_0} \sigma_{-s}\} } \cap \ceil{ \pi_{[m]_0}(-\Sigma_{K,I,v} \cap \R^{[m]_0})} = \emptyset \\
& \LRA & \pi_{[m]_0} \ceil{ \{\sigma_{-s}\} } \cap  \pi_{[m]_0}  \ceil{ -\Sigma_{K,I,v} \cap \R^{[m]_0}} = \emptyset \\
& \LRA &  \ceil{ \{\sigma_{-s}\} } \cap  \pi_{[m]_0}^{-1} \pi_{[m]_0}  \ceil{ -\Sigma_{K,I,v} \cap \R^{[m]_0}} = \emptyset \\
& \LRA &  \ceil{ \{\sigma_{-s}\} } \cap   \ceil{ -\Sigma_{K,I,v} \cap \R^{[m]_0}} = \emptyset \\
& \LA &  \ceil{ \{\sigma_{-s}\} } \cap   \ceil{ -\Sigma_{K,I,v} } = \emptyset 
\eea 
where in the next to last step, we used that $ \ceil{ -\Sigma_{K,I,v} \cap \R^{[m]_0}} $ is translation-invariant by $\ker(\pi_{[m]_0})$, to get $ \pi_{[m]_0}^{-1} \pi_{[m]_0}  \ceil{ -\Sigma_{K,I,v} \cap \R^{[m]_0}}  =  \ceil{ -\Sigma_{K,I,v} \cap \R^{[m]_0}} $. 
\epf 

\subsection{Proof of Theorem \ref{t:SSL-Hom-mezzaine}}
Since the claim depends on the local behavior of $\La_{W_K}$, we only need to check locally near each lattice point $\wt v \in \Z^N$. Let $v = \pi(\wt v) \in \Z^m$, and denote as before $\La_{W_K, v}$ the specialization of $\La_{W_K}$ to $\wt v$. It suffices to prove that 
$$ SS_{Hom}^L(\pi_* Sh^\dm_{\La_{W_K, v}})  \In T^*_{\R^m} \R^m. $$ 
From Proposition \ref{p:no-leak-no-flood}, we see $\La_{\ical(K,v)}$ satisfies condition (3) in Proposition \ref{p:coFF}, hence satisfies condition (1) there.

\section{Global Window Skeleton and Sheaf-theoretic Parallel Transport\label{s:global} }
In this section, we mainly consider the global window skeleton $\La_\delta = \La_{W_\delta} \In T^*\R^N$ for a fixed shift parameter $\delta \in \R^k$. We will show that for generic $\delta$, the family of skeleton $\La_{\delta, l} = \La_\delta|_l$ parametrized by $l \in \R^k$ forms a non-characteristic deformation. It is more interesting when $\delta$ is non-generic, when the constructible sheaf of category $\ccal_\delta = Sh^\dm_{\La_{W_\delta}}$ has singular support on the discriminant locus (not just on the boundary, but also in the interior). Roughly speaking, the jumping loci is on the affine hull of the faces of $\nabla_\delta$ that contains lattice points, and the 'jump amount' is associated to the semi-orthogonal decomposition. 

We establish some notation that will be used throughout this section. Let $\Sigma_\nabla$ be the exterior conormal fan of $\nabla$, and for a cone $\sigma \in \Sigma_\nabla$, let $F_\sigma$ denote the corresponding face in $\nabla$ whose exterior conormal is $\sigma$. 
Define the tangent cone of $\nabla$ to the face $F_\sigma$ as
\be \label{eq:CF} C_{F_\sigma, \nabla} = \R_{\geq 0} \cdot (\nabla - F_\sigma). \ee
Define the index subset 
\be \label{eq:Is}  I_\sigma = \{ i \in [N] \mid \beta_i \notin C_{F_\sigma, \nabla} \}, \ee
and the shifted tangent cone
$$ C_{\delta, \sigma} = \delta + F_\sigma + C_{F_\sigma, \nabla} $$ adjacent to $\delta + F_\sigma$.

For $\wt v \in \Z^N$, we define $L_{\wt v} = \C_{\wt v + \R_{> 0}^N}$. 

\subsection{Generic shifted zonotope}

\bt \label{t:l-indep}
If $\delta \in \R^k$ is generic, then for any $l \in \R^k$, the restriction to the fiber $X_l = \mu^{-1}(l)$ induces an equivalence of categories
$$ \rho_{\delta, l}: Sh^\dm(\R^N, \La_{W_\delta}) \to Sh^\dm(X_l, \La_{W_\delta}|_{l}) $$
\et
\bpf 
Let $\ccal_\delta = Sh^\dm(\R^N, \La_{W_\delta})$ and $\ccal_{\delta, l} =  Sh^\dm(X_l, \La_{W_\delta}|_{l})$. To show the equivalence, it suffices to show that $\rho_{\delta, l}$ and its left-adjoint is fully-faithful.

First, we show $\rho_{\delta, l}$ is fully-faithful. For any two objects $F, G$ in $\ccal_\delta$, we have $\mu_* \uhom(F, G)$ being a local system on $\R^m$. Hence 
$$ \Gamma(\R^N, \uhom(F,G)) = \Gamma(\mu^{-1}(l), \rho_{\delta, l} \uhom(F,G)) = \Gamma(\mu^{-1}(l), \uhom(\rho_{\delta, l}  F, \rho_{\delta, l} G)) $$
where the last step is because $SS(F), SS(G), SS(\uhom(F,G))$ are non-characteristic with respect to $\mu$. 

Next, we show the left-adjoint $\rho_{\delta, l}^L$ is fully-faithful. 
We may first co-restrict from $X_l$ to a tubular neighborhood of it. Let $B_l$ be a small enough ball around $l$, then by Proposition \ref{p:coFF-fiber}, we have the co-restriction functor 
$$ \rho_{\delta, l, B_l}^L: Sh^\dm(X_l, \La_{W_\delta}|_{l}) \to Sh^\dm(\mu^{-1}(B_l), \La_{W_\delta}|_{\mu^{-1}(B_l)}) $$
being fully-faithful. By further composition with the co-restriction from $B_l$ to $\R^k$, we see $\rho_{\delta, l}^L$ is fully-faithful.
\epf 

From this theorem, we can already deduce the main theorem in the introduction, which we restate here. 
\bt 
The universal window skeleton $\La$ over the punctured base $\bcal^o = \R^k_\delta \times \R^k_l \RM \dcal$ defines a local system of categories whose value over $(\delta, l)$ is $Sh^\dm(\R^n, \La_{\delta, l})$, where we identified $\R^n=\R^{N-k}$ with the fiber over $\wt \mu^{-1}((\delta, l))$. 
\et 
\bpf 
By the construction of $\La$, $\La_{\delta, l}$ is locally constant in $\delta$. By Theorem \ref{t:l-indep}, for generic $\delta$ (hence for all $\delta$), $\ccal_{\delta, l}$ is independent of $l$. Hence we have a local system of categories as claimed. More concretely, given an embedded smooth path in $\bcal^o$, we may straighten it to a piecewise linear path with each segment constant in $\delta$ or $l$, then the parallel transport along the constant $l$ segment is the identity, and the parallel transport along the constant $\delta$ path is by co-restriction then restriction. 
\epf 

\subsection{Window objects as generators}
Let $\La_{full} \In T^*\R^N$ be the full skeleton
$$ \La_{full} = \Z^N + SS(\C_{\R_{>0}^N}),$$
which can also be described as the equivariant FLTZ skeleton for $\C^N$. 
For any $\delta \in \R^k$, let $\nabla_\delta, W_\delta, \wt W_\delta$ be given as before, then we define the equivariant 'window subcategory'
$$ \acal_\delta = \la \{  L_w  \mid w \in \wt W_\delta\} \ra \In Sh^\dm(\R^N, \La_{full})$$
as the full triangulated subcategory generated by the window objects. This is a direct translation of the B-model window subcategory. We show that the window subcategory is precisely the constructible sheaves on $\R^N$ admissible for the window skeleton. 

\bt\label{t:generation}
The canonical fully-faithful embedding $\acal \to Sh^\dm(\R^N, \La_\delta)$ is an equivalence of categories. 
\et 
\bpf 
We only need to prove that $\{L_w  \mid w \in \wt W_\delta\}$ generate the category $Sh^\dm(\R^N, \La_\delta)$. Let $Q_w = w + (0,1)^N$ be the open cell, and $F_w$ be the probe sheaf for $Q_w$ for skeleton $\La_\delta$. It suffices to prove that $F_w$ can be generated using window objects $\{L_w  \mid w \in \wt W_\delta\}$. 

We will prove this by induction.  First, for any $r \in \R, r \geq 1$ we defined the rescaled zonotope and windows 
$$ \nabla_{\delta, r} = \delta + r \nabla,  \quad W_{\delta, r} = \Z^k \cap \nabla_{\delta, r}, \quad \wt W_{\delta, r} = \mu_\Z^{-1} (W_{\delta, r}).$$ 
As $r$ increases, only for $r$ in a sequence of $1=r_0 < r_1 < r_2 \cdots$ does $W_{\delta, r}$ change. 
Assume for $r \leq r_k$, $F_w$ with $w \in \wt W_{\delta, r}$ is generate by the window objects. Then the case for $k=0$, $r \leq r_0=1$ is clear, since $F_w = L_w$ for $w \in \wt W_{\delta, 1}$.  We now prove the hypothesis for $r \leq r_{k+1}$. Suppose $w \in \wt W_{\delta, r_{k+1}} \RM \wt W_{\delta, r_{k}}$, then $v=\mu(w)$ is on the boundary of $\nabla_{\delta, r}$.
Let $F_\sigma$ be the minimal face of $\nabla$ such that $\delta + r F_\sigma$ contains $v$, and $C_{F_\sigma, \nabla}, I_\sigma$ be defined as in Eqs. \eqref{eq:CF} and \eqref{eq:Is}. Then, $\tau_{I_\sigma^c}$ is a leak for $Q_\emptyset$ in the local skeleton $\La_{W_\delta, v}$, since there is no points in $W_\delta$ that is contained in $v - \tau_{I_\sigma^c}$. Thus, we have the acyclic complex
$$ F_w \to \bigoplus_{I \In I_\sigma, |I|=1}  F_{w - e_I} \to \bigoplus_{I \In I_\sigma, |I|=2}  F_{w - e_I} \to \cdots $$
Since for any $I \In I_\sigma$, $|I| > 0$,  we have 
$$ \mu(w-e_I) \in (v + \sum_{i \in I_\sigma} [0, -\beta_i]) \RM \{v\} \cap \Z^k \In \Int(\nabla_{\delta, r_{k+1}}) \cap \Z^k \In W_{\delta, r_{k}} $$
hence $F_w$ can be generated by $F_{w'}$ with $w' \in \wt W_{\delta, r_{k}}$ hence by induction hypothesis can be generated by the window objects $L_w$ for  $w \in \wt W_{\delta, 1}$.  
\epf 

\brem 
The induction is in a similar fashion as \cite{halpern2020combinatorial} and  \cite{vspenko2017non}, using rescaled zonotope and induction on the radius. 
\erem

\subsection{Microlocal Stalk and Jumping Loci \label{s:microlocal-stalk}} 
Next, we consider the case where $\delta$ is fixed and non-generic. Recall that in general
$$ SS(\mu_* Sh^\dm_{\La_\delta}) =  SS_{Hom} (\mu_* Sh^\dm_{\La_\delta}) \cup   SS_{Hom}^L (\mu_* Sh^\dm_{\La_\delta}). $$
By Theorem \ref{t:SSL-Hom-window}, we have $SS_{Hom}^L (\mu_* Sh^\dm_{\La_\delta})$ is just the zero section. By Theorem \ref{t:SS-Hom}, we have $SS_{Hom} (\mu_* Sh^\dm_{\La_\delta}) = \mu_* SS_{Hom} (Sh^\dm_{\La_\delta})$ is given explicitly. Here we describe the jump associated to $0 \neq (x, \xi) \in SS(\mu_* Sh^\dm_{\La_\delta})$.

\bp \label{p:microstalk-1}
Let $(x, \xi) \in T^* \R^k$ be a non-zero covector in the smooth part of the singular support $SS(\mu_* Sh^\dm_{\La_\delta})$. 

(1) Then there is a unique face $F_\sigma$ of the zonotope $\nabla$ with exterior conormal the cone $\sigma$, such that $x \in \Aff(\delta + F_\sigma)$ and $\xi \in \Int(-\sigma)$. 

(2) Let $B_x$ and $B_{x,\xi,-}$ be given by 
$ B_x = \{x' \mid |x'-x|<\epsilon\}, \quad B_{x,\xi,-} = \{x' \in B_x \mid \la x'-x, \xi \ra < 0 \},$
where $\epsilon$ is small enough. Denote the sheaf of categories $\mu_* Sh^\dm(\La_\delta)$ as $\ccal_\delta$, then we have fully-faithful co-restriction
$$ \rho^L_{x,\xi}: \ccal_\delta(B_{x,\xi,-}) \to \ccal_\delta(B_{x}).$$
If we identity $\ccal_\delta(B_{x,\xi,-})$ as the full subcategory in $\ccal(B_x)$, then we have a semi-orthogonal decomposition
$$ \ccal_\delta(B_{x}) = \la \ccal_\delta(B_{x,\xi,-})^\perp, \ccal_\delta(B_{x,\xi,-}) \ra, $$
where $\ccal_\delta(B_{x,\xi,-})^\perp$ is the full subcategory of  $\ccal(B_x)$ consisting of sheaves that vanish under restriction to $B_{x,\xi,-}$.

\ep
\bpf 
(1) The geometric statement follows from the vanishing (away from zero-section) of $SS^L_{Hom}(\mu_* Sh^\dm_{\La_\delta})$ (Theorem \ref{t:SSL-Hom-window}), and the explicit description of $SS_{Hom}(\mu_* Sh^\dm_{\La_\delta})$ (Theorem \ref{t:SS-Hom} (3)). 

(2) If $F \in \ccal_\delta(B_{x,\xi,-})^\perp$, then for any $G \in \ccal_\delta(B_x)$, we have
$$ 0 = \Hom(\rho^L_{x,\xi} G, F) = \Hom(G, \rho_{x,\xi} F) $$
Hence $\rho_{x,\xi} F$ has to vanish. 

\epf

In fact, more is true. These microlocal stalks or vanishing cycles $\ccal_\delta(B_{x,\xi,-})^\perp$ form a constructible sheaf of categories as $x$ varies on $\Aff(\delta + F_\sigma)$. In the following, we fix a face $F_\sigma$, such that $(\delta + F_\sigma) \cap \Z^k \neq \emptyset$. Our goal is to describe the sheaf of categories $\ccal_\ds$ on the affine space $V_{\delta, \sigma}$, such that for any $x \in V_{\delta, \sigma}$ and $\xi \in \Int(-\sigma)$, and $B_x, B_{x,\xi,-}$ as in Proposition \ref{p:microstalk-1}, we have $\ccal_{\delta, \sigma}(B_x) = \ccal_{\delta}(B_{x,\xi,-})^\perp$. 

Define the face zonotope, and its affine hull $$ \nabla_{\delta, \sigma} = \delta + F_\sigma, \quad V_{\delta, \sigma} = \Aff(\nabla_{\delta, \sigma}). $$
Define a partition of $[N] = I_{\sigma, +} \sqcup I_{\sigma, 0} \sqcup I_{\sigma, -}$, by fixing any element $\xi$ in the interior of $\sigma$, and let
$$I_{\sigma, s} := \{ i \in [N] \mid \sign(\la \beta_i, \xi \ra) = s \} \z{ for } s= +,0,-. $$ 
Define an affine lattice and a lattice that acts on it
$$ V^\Z_{\delta, \sigma} = V_{\delta, \sigma} \cap \Z^k, \quad U_{\sigma}^\Z = \Z \cdot \{\beta_i: i \in I_{\sigma,0}\}, \quad U_{\sigma}^\R = \R \cdot \{\beta_i: i \in I_{\sigma,0}\}.$$
Upstairs in $\R^N$ and $\Z^N$, we define
$$ \wt V_{\delta, \sigma} = \mu^{-1}(V_{\delta, \sigma} ), \quad \wt V^\Z_{\delta, \sigma} =  \mu_\Z^{-1}(V^\Z_{\delta, \sigma} )$$
and $\Z^{I_{\sigma, 0}}$ acts on $\wt V^\Z_{\delta, \sigma}$ by translation. 

Let $\R^{I_{\sigma,0}}$ act on $\wt V_\ds$ by translation, and for any $\wt v\in \wt V^\Z_\ds$, consider the orbit $O_{\wt v} = \wt v + \R^{I_{\sigma,0}}$. Let $O = O_{\wt v}$, and define the following partial conormal
$$ \La_{O} := O \times (\R^\vee_{<0})^{I_{\sigma,+} \sqcup I_{\sigma,-}} \In \R^N \times (\R^\vee)^N = T^* \R^N. $$

\bl 
For any $\R^{I_{\sigma,0}}$ orbit $O$ in $\wt V_\ds$ that passes through a lattice point $\wt V^\Z_\ds$, the Lagrangian $\La_{O} \In \La_{\delta}$. 
\el 
\bpf 
Suffice to check at every lattice point $\wt w \in O \cap \Z^N$, that the specialization $\La_{\delta, \wt w}$ contains the $\La_{I_{\sigma, 0}}$. This is equivalent to 
\bea 
&& \wt w - e_{I_{\sigma, 0}} - \tau_{I_{\sigma, 0}} \cap \wt W_\delta \neq \emptyset \\
& \LRA & w - U_\sigma \cap W_\delta \neq \emptyset \\
& \LRA & V_\ds \cap W_\delta \neq \emptyset
\eea 
which holds since $\nabla_\ds \cap \Z^k \neq \emptyset$. \epf 

We are interested in the restriction of $\La_\delta$ to a small Weinstein neighborhood $\Omega_O$ of $\La_O$, or equivalently the specialization $\La_\delta|_{\La_O} \In T^* (\La_O)$. Then only Lagrangian component of the form $A \times B$ in $\La_\delta$ will contribute, where $A \In O \In \R^N$ and $(\R^\vee)^N \supset B \supset (\R^\vee_{<0})^{I_{\sigma,+} \sqcup I_{\sigma,-}}$. 

The following is a description of the specialization $\La_\delta|_{\La_O}$. 
\bp 
Let $O_\Z = O \cap \Z^N$, consider the restriction $\mu_{O, \Z}$ of $\mu_\Z$ to $O_\Z$ with image $\wb O_\Z$. Then $\wb O_\Z$ is a $U_\sigma^\Z$ orbit in $V^\Z_\ds$. Let $\La_{\nabla_\ds, O} \In T^* O$ be the window skeleton associated to the map $\mu_{O, \Z}: O_\Z \to \wb O_\Z$ and the shifted zonotope $\nabla_\ds$. Then we have 
$$ \La_\delta|_{\La_O} \cong \La_{\nabla_\ds, O} \times (\R^\vee_{<0})^{I_{\sigma,+} \sqcup I_{\sigma,-}}. $$
In other word, $\La_\delta|_{\La_O}$ is $\La_{\nabla_\ds, O} $ up to 'stabilization' by multiplying the zero-section of $T^* (\R^\vee_{<0})^{I_{\sigma,+} \sqcup I_{\sigma,-}}$. 
\ep 
\bpf 
Suffice to verify this locally at each lattice point $\wt v \in O$. Consider the specialization to point $\wt v$,  $\La_{\delta, \wt v}$, then $\La_O$ specializes to $\La_{I_{\sigma,0}}$.  We note that only $\La_I$ with $I \In I_{\sigma, 0}$ will contribute to the specialization to $\La_{I_{\sigma,0}}$. The condition that $I \In I_{\sigma,0}$ and $\La_I \In \La_{\delta, \wt v}$ is equivalent to 
\bea 
 ( \mu(\wt v) - \beta_{I} - \mu(\tau_{I}) )\cap \wt W_\delta \neq \emptyset 
& \LRA & ( \mu(\wt v) - \beta_{I} - \mu(\tau_{I})) \cap (\nabla_\ds \cap \Z^k) \neq \emptyset 
\eea 
which is equivalent to $\La_I$ (in $T^* O$) is in the specialization $\La_{\nabla_\ds, O, \wt v}$. 
\epf 

\bl 
Let $\wt v \in \wt V^\Z_\ds \In \Z^N$ be a lattice point, and let $\La_{\delta, \wt v}$ and $\La_{\nabla_{\delta, \sigma}, O, \wt v}$ be specialization of $\La_{\delta}$ and $\La_{\nabla_{\delta, \sigma}, O}$  at $\wt v$. Let $F$ be a constructible sheaf in $Sh(\R^N, \La_{\delta, \wt v})$, such that $\mu(\supp(F)) \In C_{F_\sigma, \nabla}$. Then there exists a constructible sheaf $F_0$ in $Sh(T_{\wt v} O, \La_{\nabla_{\delta, \sigma}, O, \wt v})$, such that 
$$ F = (\wt v_+ +  (\R_{\leq 0})^{I_{\sigma, +}} ) \boxtimes (\wt v_- + \R_{> 0}^{I_{\sigma,-}} ) \boxtimes F_0.  $$
Conversely, for any $F_0$ in $Sh(T_{\wt v} O, \La_{\nabla_{\delta, \sigma}, O, \wt v})$, the above construction for $F$ satisfies the support condition that $\mu(\supp(F)) \In C_{F_\sigma, \nabla}$. 
\el 
\bpf 
We note that any sheaf in $Sh^\dm(\R^N, \La_N)$ such that $\mu(\supp(F)) \In C_{F_\sigma, \nabla}$ has to be of the form $(\R_{\leq 0})^{I_{\sigma, +}} \boxtimes \R_{> 0}^{I_{\sigma,-}} \boxtimes F_0$ for some sheaf $F_0$ on $T_{\wt v} O \cong \R^{I_{\sigma, 0}}$. Now suffice to show that $SS(F_0) \In \La_{\nabla_{\delta, \sigma}, O, \wt v}$. Let 
$$ \ical_0 = \{ I \In I_{\sigma,0} \mid \La_I \In \La_{\nabla_{\delta, \sigma}, O, \wt v} \}, \quad \ical = \{ I \In [N] \mid \La_I \In \La_{\delta, \wt v} \}.  $$
First we note that for any $I_+ \In I_{\sigma, +}$ and $I_0 \in \ical_0$, we have $I_+ \sqcup I_0 \in \ical$. Indeed, $\nabla_\delta + \beta_{I_+ \sqcup I_0} + \mu(\tau_{I_+ \sqcup I_0})$ contains $v = \mu(\wt v)$. Conversely, if there is an $I_0 \In I_{\sigma_0}$, such that for any $I_+ \In I_{\sigma, +}$ we have $I_+ \sqcup I_0 \in \ical$, then $I_0 \in \ical_0$. Indeed, taking $I_+ = \emptyset$ will show. 
\epf 

For any $\R^{I_{\sigma,0}}$ orbit $O$ in $\wt V_\ds$ that passes through a lattice point $\wt V^\Z_\ds$, we define a sheaf of categories on $O$, 
$$ \ccal_{\ds, O}: = Sh^\dm_{\La_{\nabla_\ds, O}}. $$

And we define the sheaf of vanishing cycles $\ccal_{\ds}$ along $V_\ds \In \R^k$ as a sub-sheaf of $\ccal_\delta$, such that for any $x \in V_\ds$ and small ball $B_x$ around $x$, we have
$$ \ccal_{\ds}(B_x \cap V_\ds) = \z{ full subcategory of $\ccal_\delta(B_x)$ consisting of object $F$ such that $\mu(\supp(F)) \In C_\ds \cap B_x$. } $$ 

\def\riz{\R^{I_{\sigma,0}}} 

For each $\riz$ orbit $O$ that passes a lattice point, 
we define the {\bf microlocal stalk functor} between the two sheaves on $V_\ds$: 
$$ \psi_O: \ccal_\ds \to  \ccal_{\ds, O}. $$
Here we abuse notation and identify the sheaf $\mu_* \ccal_{\ds, O}$ on $V_\ds$ as the sheaf on $\ccal_{\ds, O}$ on $O$ since $V_\ds \cong O$ by $\mu$. 
For any $x \in V_\ds$ and small ball $B_x$ around $x$, and for any $F \in \ccal_\ds(B_x \cap V_\ds) = Sh^\dm(\mu^{-1}(B_x), \La_\delta)$, we define
$$ \psi_O (F) =  p_* \psi_{\La_O}(F) \in Sh^\dm(O|_{B_x}, \La_{\nabla_\ds, O}) $$
where the specialization $\psi_{\La_O}(F)$ values in $Sh^\dm(\La_O|_{B_x}, \La_\delta|_{\La_O})$ and $p: \La_O \to O$ is the projection map, quotient out the stabilization factor. Equivalently, up to a degree shift, $\psi_O(F)$ is the restriction of $F$ to $O + \epsilon (e_{I_{\sigma, -}} - e_{I_{\sigma, +}})$.

We also define the fully-faithful embedding
$$ \iota_O: \mu_* \ccal_{\ds, O} \to \ccal_\ds $$
where in the setup as above, we send $F_O \in Sh^\dm(O|_{B_x}, \La_{\nabla_\ds, O})$ to the germ of $(\wt x_+ + \R_{\leq 0}^{I_{\sigma, +}})[-|I_{\sigma, +}|] \boxtimes (\wt x_- + \R_{>0}^{I_{\sigma, -}}) \boxtimes F_O$ near $O$, where $\wt x$ is the unique $\mu$-preimage of $x$ to $O$, and $\wt x_\pm$ are the components of $\wt x$ in the decomposition $\R^N \cong  \R^{I_{\sigma, +}} \times \R^{I_{\sigma, -}} \times \R^{I_{\sigma, 0}}$. 

\bt \label{t:decomposition}
Taking microlocal stalk along all orbits $O$ is an equivalence of categories
$$ \psi_\ds = \bigoplus_O \psi_O: \ccal_\ds \to  \bigoplus_{O}  \ccal_{\ds, O}, $$
In other words, for any convex open $U \in V_\ds$, $F \in \ccal_\ds(U)$, we have
$$ F \cong \sum_O \iota_O (\psi_O(F)) $$
where the summation is over each $\riz$ orbit $O$ that passes a lattice point in $\wt V^\Z_\ds$. 
\et 
\bpf 
Indeed, for any convex open $U \In V_\ds$, if $F \in \ccal_\ds(U)$, then $F$ is a germ of sheaves on $\R^N$ whose support is contained in the disjoint union of neighborhood of lattice $\riz$ orbit $O$. In particular, since $\mu(\supp(F)) \In C_\ds$, along any $O$ orbit and any $\wt x \in O \cap \mu^{-1}(U)$, $F$ has a common factor in the $\R^{I_{\sigma, +}} \times \R^{I_{\sigma, -}}$ direction, namely the constant sheaf supported on $\R_{\leq 0}^{I_{\sigma, +}} \times \R_{>0}^{I_{\sigma, -}}$ translated by $(\wt x_+, \wt x_-)$, which is a constant shift depending only on $O$. The specialization $\psi_O(F)$ amounts to factoring out this common factor. 
Hence the total $\Psi_\ds$ is an equivalence.
\epf

Finally, we take global sections and describe some global generators using window objects. 

We need to introduce certain global vanishing cycles associated to the data $(\delta, \sigma, \wt v)$, where $\delta \in \R^k$, $\sigma$ is a cone in the exterior conormal fan of the zonotope $\nabla$ corresponding to the face $F_\sigma$, and $\wt v \in \mu_\Z^{-1}( (\delta + F_\sigma) \cap \Z^k )$. 
Then we define
\be \label{e:L-sigma-v} L_{\sigma, \wt v} = \uwave{L_{\wt v}} \to \bigoplus_{I \In I_\sigma, |I|=1} L_{\wt v - e_I} \to \bigoplus_{I \In I_\sigma, |I|=2} L_{\wt v - e_I} \to \cdots.
\ee

\bl \label{l:Koszul-supp}
$L_{\sigma, \wt v}$ is an object in $Sh^\dm(\R^N, \La_{\delta})$.  Up to a degree shift, it is quasi-isomorphic the product of constant sheaf
$$ L_{\sigma, \wt v} [|I_{\sigma,+}|] \cong (\wt v_+ +  (-1, 0]^{I_{\sigma, +}} ) \boxtimes (\wt v_- + \R_{> 0}^{I_{\sigma,-}} ) \boxtimes (\wt v_0 + \R_{> 0}^{I_{\sigma,0}} ).  $$
And its support and $\mu$-image are
$$\supp(L_{\sigma, \wt v}) = \wt v + [-1,0]^{I_\sigma} + \tau_{I_\sigma^c}, \quad \mu(\supp(L_{\sigma, \wt v})) = \delta + F_\sigma + C_{F_\sigma, \nabla}. $$ 
\el 
\bpf 
For the first statement, suffice to note that for any $I \In I_\sigma$, $\mu(\wt v - e_I) = v - \beta_I$ is contained in $\delta + \nabla$, hence $L_{\wt v - e_I}$ is in $Sh^\dm(\R^N, \La_{\delta})$. For the support, we note that in the $\R^{I_\sigma}$ factor, we have a cube in $[-1,0]^{I_\sigma}$, and for the remaining factor, it is  $\R_{\geq 0}^{I_{\sigma}^c }$. Then its image under $\mu$ is 
$$ \mu(\supp(L_{\sigma, \wt v}) ) = v + \sum_{i \in I_\sigma} [0, -\beta_i] + \sum_{i \notin I_\sigma} \R_{\geq 0} \beta_i = v + C_{F_\sigma, \nabla} = \delta + F_\sigma + C_{F_\sigma, \nabla} = C_\ds. $$ 
\epf 

\bl \label{l:ess-surj}
For any convex open set $U \In V_\ds$, 
the restriction functor $C_\ds(V_\ds) \to C_\ds(U)$ is essentially surjective. 
\el 
\bpf 
We define the extension functor using co-restriction along the $O$ direction. 
$$ C_\ds(U) \into C_\ds(\R^k), \quad F \mapsto \sum_O \iota_O \circ \rho_{U_O}^L \circ \psi_O(F), $$
where $U_O = \mu^{-1}(U) \cap O$, and $\rho_{U_O}^L: \ccal_{\ds, O}(U_O) \to \ccal_{\ds, O}(O)$ is the co-restriction. 
\epf 

\bp 
The collection of global vanishing cycles $\{ L_{\sigma, \wt v} \mid \wt v \in \mu_\Z^{-1}( (\delta + F_\sigma) \cap \Z^k ) \}$ generates $\ccal_\ds(V_\ds)$. Furthermore, for any convex open $U \In V_\ds$, the restriction of these generators to 
$U$ generates $\ccal_\ds(U)$. 
\ep 
\bpf 
The first statement follows from the decomposition Theorem \ref{t:decomposition}, and the global generation Theorem \ref{t:generation} for window skeletons $\La_{\nabla_\ds, O}$. The second statement follows from the essential surjectivity of restriction in Lemma \ref{l:ess-surj}. 
\epf

Since the singular support in $SS(\ccal_\delta)$ has 'hair' pointing inward towards the zonotope , namely for any non-zero $(x,\xi) \in SS(\ccal_\delta)$, we have $\la x-\delta, \xi \ra < 0 $, the restriction to the stalk at a point $x \in \Int(\nabla_\delta)$ in the interior of the zonotope is an equivalence
$$\rho_x: \ccal_\delta(\R^k) \to \ccal_\delta(B_x), $$
where $B_x$ is a small enough ball at $x$. Furthermore, since there exists fully-faithful co-restriction from the category of constructible sheaves on the fiber to the germ of the fiber (Proposition \ref{p:coFF-fiber}), and restriction is also fully-faithful since $x$ is in the interior of $\nabla_\delta$, we have $\ccal_\delta(B_x) \cong Sh^\dm(\mu^{-1}(x), \La_\delta|_x)$. 

If $\gamma: [0,1] \to \R^k$ is a smooth path in $\R^k$, such that whenever $\gamma$ crosses a jumping locus in $\R^k$ for $\La_\delta$, $\dot \gamma$ points in the direction of the hair, $\la \dot \gamma, \xi \ra > 0$, and the endpoints of $\gamma$ are not in the jumping loci, then we have an fully-faithful embedding 
$$ \rho^L: \ccal_\delta(\gamma(0)) \to \ccal_\delta(\gamma(1)) $$
which is independent of the choice of $\gamma$, since $\rho^L$ is the composition of co-restriction to $\R^k$ and restriction to $\gamma(1)$. And there is a semi-orthogonal decomposition that depends on $\gamma$ 
$$ \ccal_\delta(\gamma(1)) = \la T_k, T_{k-1}, \cdots, T_1,   \ccal_\delta(\gamma(0)) \ra $$
where $T_i = \ccal_{\delta, \sigma(\gamma(t_k))}(\gamma(t_k))$, and $0 < t_1 < t_2 <  \cdots < t_k < 1$ are the jumping loci.

\bibliographystyle{amsalpha}
\bibliography{HMS}

\end{document}